      \documentclass[a4paper,11pt]{article}

\usepackage[T2A]{fontenc}
\usepackage[english]{babel}

\usepackage{latexsym}
\usepackage{amssymb,amsmath,amsfonts}

\usepackage[all]{xy}
\usepackage{graphics}
\usepackage{epic}
\usepackage{eepic}
\usepackage{epsfig}

\newcommand{\torsp}{\mathbb T^\nu}
\newcommand{\param}{\mathbb R^l}
\newcommand{\RR}{\mathbb R}
\newcommand{\ZZ}{\mathbb Z}

\newcommand{\lam}{\lambda}

\newcommand{\const}{{\mbox{\rm const}}}
\newcommand{\supp}{{\mbox{\rm supp}}\,}
\newcommand{\m}{\mathfrak m}
\newcommand{\A}{{\it E}}  

\renewcommand{\v}{\bar a} 
\newcommand{\w}{\bar d}
\newcommand{\x}{\bar x}

\newtheorem{theorem}{Theorem}[section]

\newtheorem{definition}[theorem]{Definition}
\newtheorem{statement}[theorem]{Statement}
\newtheorem{note}[theorem]{Remark}

\title{Classification of Singularities and Bifurcations of Critical Points
 of Even Functions}
\author{E.A.~Kudryavtseva\thanks{Moscow State University.
e-mail: eakudr@mech.math.msu.su},\quad  E.L.~Lakshtanov
\thanks{Moscow State University.
e-mail: lakshtanov@rambler.ru}}
\date{}
\begin{document}

 \maketitle

\begin{abstract}
Singularities of even smooth functions are studied. A classification of singular
points which appear in typical parametric families of even functions with at
most five parameters is given. Bifurcations of singular points near a caustic
value of the parameter are also studied. A determinant for singularity types
and conditions for versal deformations are given in terms of partial
derivatives (not requiring a preliminary reduction to a canonical form).
\end{abstract}

\section{Introduction}
In this work, we consider families of smooth even functions
$\omega_\lam(k) : \torsp \rightarrow \mathbb R$ on a torus $\torsp=(S^1)^\nu$
where $S^1=[0,2\pi]/0\sim2\pi$ and the parameter $\lam$ takes values in some
domain $\param$:
 $$
\refstepcounter{theorem}
\label{fam}
 \omega_\lam(k)=\omega_\lam(-k), \quad k \in \torsp,\ 
 \lam \in \param.
\eqno (\thetheorem)
 $$
We will study bifurcations of critical points for smooth families of
even smooth functions, their singularities, and, in particular, obtain
a classification of degenerate critical points for generic families
of germs of even functions, which are determined in a neighbourhood of zero
(or in a neighbourhood of a fixed point of the involution $k\mapsto -k$),
with the number of parameters $l \leq 5$. Moreover, we will show how to
determine types of critical points of even functions which appear in typical
parametric families with $l\le 4$ parameters, and formulate conditions on
their deformations to be typical (more precisely, versal), in terms of partial
derivatives (not requiring a preliminary reduction to a canonical form).

The study of properties of families of even functions appears in problems
about planetary systems with satellites, integrable Hamiltonian systems
with 2 degrees of freedom, problems of spectral analysis of stochastic operators
of multiparticle systems, and problems involving two (quasi)particles
interaction.
It was shown in~\cite{LakshtanovMinlos, MMinfP, MamatovMinlos, Sinay} how
the type of singularities of a family $\omega_\lam(k)$ affect the spectrum
structure. An extremely important task is therefore to describe
{\it stable families} of even functions in a small neighbourhood of a
singular point (they are also called {\it families in generic position}
or {\it typical families}), i.e.\ families which do not change their properties
under small perturbations (for the precise definition see~\ref {note:vozmTh}).

We tried to make this work readable
for those readers who are not familiar with singularity theory,
as well as to demonstrate, using simple examples, how to use results
of~\cite{Arnold} which are very important for applications.

The authors are grateful to S.M.~Gusein-Zade, M.E.~Kazaryan, G.~Wassermann, and
V.M.~Zakalyukin for useful discussions.

\refstepcounter{theorem}
\subsection{Critical points of even functions on a torus $\torsp$}
The automorphism $k\mapsto -k$ of the torus has $2^\nu$ fixed points
\[
(k_1,\ldots,k_\nu), \quad k_i \in \{0,\pi \} \subset S^1, \quad i=1,\ldots,\nu.
\]
These points are obviously critical. We will call them {\it basic}.
Observe that, if a non-basic point $k \in \torsp$ is critical, then the point
$-k$ is also critical. Thus, non-basic critical points appear in pairs
$(k,-k)$, moreover the degeneracy type of these points is the same.
We will call such points {\it twin} or {\it additional}.

Recall~\cite[v.~2, 6.1.$\Gamma$]{Arnold} that a value of the parameter $\lam$
is called
{\it caustic} if the function $\omega_\lam(k)$ admits a degenerate critical
point. It is obvious that, if some value $\lam=\lam_0$ is not caustic, the
number of critical points and their types remain the same in a neighbourhood of
$\lam_0$. Thus, appearance and disappearance of critical points is related
with the passage of the parameter $\lam$ across hypersurfaces
(i.e.\ surfaces of codimension one) of caustic values. For a detailed study
of these processes, let us describe the types of degenerate critical points,
which are non-removable for the entire family of functions. In particular,
this is necessary for study of oscillating integral asymptotics, see below.

\section{Classification of singularities of even smooth functions}
Thus, for even functions $\omega_\lam(k)$, there are two types of critical
points. In a neighbourhood of a basic critical point, the function
$\omega_\lam(k)$ is even, while, in a neighbourhood of an additional critical
point, it is an arbitrary function (i.e.\ a function which is not assumed to be
symmetric relative to the critical point).
In the latter case, we can use the classification from~\cite{Arnold}.

\refstepcounter{theorem}
\subsection{A survey of the classification of critical points of arbitrary
functions due to~\cite[v.~1, 11.2]{Arnold}} 

We present here a table of germs of smooth functions at zero up to the
following transformations: adding a constant, adding a nondegenerate quadratic
form in the remaining variables (the number of these variables equals the rank
of the second differential of the function at zero), and a smooth change of
variables leaving the origin fixed. The integer $c$ in the table is the {\it
codimension of the class of singularities} in the space of function germs at
zero. It equals the minimal number $l \geq c$ of parameters of a family,
such that a critical point of the class under consideration is non-removable
under small perturbations of the family. The integer $\mu$ is the
{\it multiplicity} of the critical point $0$~(see~\cite[v.~1, 6.3] {Arnold}
and~\ref {def:mu}). For all known types of singularities, it equals the
{\it codimension} of the singularity in the space of function germs at zero
plus~$1$ (see~\ref {note:vozmTh} and~\ref {note:strat}).

\begin{note} \label {note:osob} \rm
Recall the definitions of the notions {\it singularity} and {\it class of
singularities} from~\cite {Arnold} mentioned above. If two functions can be
obtained from the same function by the above mentioned transformations, one
says that these functions have the same {\it type of critical point} or the
same {\it type of singularity} at zero, and the germs of these functions at
zero are {\it $R^+-$equivalent}
(more precisely, stably $R^+-$equivalent~\cite[v.~1,~II]{Arnold}).
The notion of $R-$equivalence is similarly determined, where adding a constant
is forbidden, see~\cite{Arnold} and~\ref {def:RO}.
An $R^+-$equivalence class of germs, i.e.\ the space of all germs with a
given type of singularity, is also called a {\it singularity}.
Singularities which can be joined to each other by a smooth path in the space
of singularities having the same multiplicity $\mu$ form a
{\it class of singularities}, more precisely a {\it class of
$\mu-$equivalent singularities}, or a {\it $\mu=\const$ stratum},
see~\ref{note:strat} and~\cite [v.~1,~15.0.1]{Arnold}. A class of
singularities is (at least for all known finite-multiple singularities)
a smooth parametric family of singularities. For such a parametric family,
the parameters are called {\it moduli}, with the number $m$
of parameters being called the {\it modality} of the singularity
(for a more general definition of modality see the introduction
of~\cite [v.~1,~II]{Arnold}). The equality $\mu=c+m+1$ is valid, where $m$
is the modality and $c$ is the codimension of the class of singularities,
see~\ref{note:strat} and~\cite [v.~1,~15.0.4]{Arnold}.
See also the examples after Remark~\ref{note:strat}.
 \end{note}

We remark that the singularities listed in Table~(\ref {obKl}) do not have
moduli ($m=0$), thus they are {\it simple}~\cite [v.~1,~II]{Arnold}.
For generic families of functions with $l \leq 5$ parameters, there are no
critical points apart from those which are equivalent to the singularities
listed in Table~(\ref {obKl}).

 $$
\refstepcounter{theorem} \label{obKl}
\begin{tabular}{|c|c|c|c|c|c|}
\hline
Type of singularity & Normal form & Restrictions & $\mu$ & $c$ & $\beta$ \\
\hline \phantom{$I^{I^{I^I}}$}\!\!\!\!\!\!\!\!\!\!
$A^\pm_k$ & $\pm x^{k+1}$      & $k\geq 1$ & $k$ & $k-1$ & $\frac{k-1}{2k+2}$\\ \phantom{$I^{I^{I^I}}$}\!\!\!\!\!\!\!\!\!\!
$D^\pm_k$ & $x^2y \pm y^{k-1}$ & $k\geq 4$ & $k$ & $k-1$ & $\frac{k-2}{2k-2}$\\ \phantom{$I^{I^{I^I}}_{j_j}$}\!\!\!\!\!\!\!\!\!\!
$E^\pm_6$ & $x^3 \pm y^{4}$    & $-$       & $6$ & $5$   & $\frac{5}{12}$    \\
 \hline
\end{tabular}
 \eqno (\thetheorem)
 $$
Seven of these singularities, with $c\le4$, are known in catastrophe theory
as the ``seven of Thom''.
Unimodular singularities ($m=1$) appear for $c=6$,
however we do not consider them here.

\begin{note} \rm
The integer $\beta$ is the {\it singularity exponent} of a given critical
point. It can be determined by means of the asymptotics of the oscillating
integral
\[
\int_{U(0)} e^{itf(x)} \varphi(x)\,dx_1 \ldots dx_\nu \sim {\const \cdot
}{t^{\beta-\nu/2}}, \quad t \rightarrow \infty .
\]
Here, the integration is taken over a small neighbourhood of the point
$x=0 \in \mathbb R^\nu$, which is a critical point of the phase $f(x)$
with singularity type under consideration, and the amplitude $\varphi(x)$
does not vanish in $U(0)$. (For more details see~\cite[v.~2, p.~134]{Arnold}.)
\end{note}

\begin{note} \rm
The listed types of critical points never correspond to a
(local) extremum (i.e.\ a point of local minimum or local maximum) apart from
the types $A_{2k-1}$, $k \geq 1$. Moreover, Table~(\ref {obKl}) implies that for
generic one-parameter families, an additional extreme point never degenerates.
V.A.~Vasiliev composed a table of germs of smooth functions in a neighbourhood
of minimum points, which occur in generic families of functions with $l\le16$
parameters, see~\cite[v.~1,~17.2]{Arnold}.
\end{note}

\refstepcounter{theorem}
\subsection{Classification of critical points of even functions}

We present here a table of germs of smooth even functions at zero up to the
following transformations: adding a constant, adding a nondegenerate quadratic
form in the remaining variables, and a smooth odd change of variables.
The integer $c_e$ in the table is the {\it even codimension of an even class
of singularities} and equals the minimal number $l \geq c_e$ of parameters of
a family, such that a critical point of the class under consideration is
non-removable under small perturbations of the family of even functions. The
integer $\mu_e$ is the {\it even multiplicity} of the function at zero
(see~\ref {def:mu}). It equals the even codimension of the singularity plus~$1$.

The notions {\it $R_O-$equivalence} of even germs (see~\ref{def:RO}),
{\it even singularity}, {\it even class of singularities}, {\it even type}
of a critical point, {\it even modality} $m_e$, and {\it even moduli}
are defined similarly to~\ref {note:osob}. The equality
$\mu_e=c_e+m_e+1$ is valid.

Even singularities of the series $A_{e,k}$ do not have moduli ($m_e=0$).
Other even singularities listed in Table~(\ref {evKl}) are unimodal
($m_e=1$). Actually, all finite-even-multiple even singularities but
$A_{e,k}$ have moduli ($m_e\ge1$)~\cite[Theorem~4.3]{Wass}.
For generic families of even functions with $l \leq 5$ parameters, there are no
even singularities at zero apart from those which are equivalent to the even
singularities listed in Table~(\ref {evKl}).

 $$
\refstepcounter{theorem}
\label{evKl}
\begin{tabular}{|c|c|c|c|c|c|}
\hline
Even class of & Normal form & Restrictions & $\mu_e$ & $c_e$ & $\beta$\\
singularities &             &              &         &       &        \\
\hline \phantom{$I^{I^{I^I}}$}\!\!\!\!\!\!\!\!\!\!
$A_{e,k}^\pm \subset A_{2k-1}^\pm$   & $\pm x^{2k}$
              & $k\geq 1$                       & $k$ & $k-1$ & $\frac{k-1}{2k}$\\ \phantom{$I^{I^{I^I}}$}\!\!\!\!\!\!\!\!\!\!
$X_{e,5}^{\pm\pm} \subset X_9^{\pm\pm}$ & $\pm x^4 +ax^2y^2 \pm y^4$
              & $a^2\neq 4$, if $++$ or $--$ & $5$   & 3       & 1/2 \\ \phantom{$I^{I^{I^I}}$}\!\!\!\!\!\!\!\!\!\!
$X_{e,r+3}^{\pm\pm} \subset X_{2r+5}^{\pm\pm}$ & $\pm x^4 \pm x^2 y^2 + a y^{2r}$
              & $a \neq 0, \quad r\geq 3$       & $r+3$ & $r+1$   & 1/2 \\ \phantom{$I^{I^{I^I}}$}\!\!\!\!\!\!\!\!\!\!
$Y_{e,r,s}^{\pm\pm} \subset Y_{2r,2s}^{\pm\pm}$ & $\pm x^{2r} + ax^2 y^2 \pm y^{2s}$
              & $a \neq 0, \quad r,s\geq 3$     & $r+s+1$ & $r+s-1$ & 1/2 \\ \phantom{$I^{I^{I^I}}$}\!\!\!\!\!\!\!\!\!\!
$\widetilde{Y}_{e,r}^\pm\subset\widetilde{Y}_{2r}^\pm$ &$\pm(x^2+y^2)^2+ay^{2r}$
              & $a > 0, \quad r\geq 3$          & $2r+1$  & $2r-1$  & 1/2 \\ \phantom{$I^{I^{I^I}}_{j_j}$}\!\!\!\!\!\!\!\!\!\!
$Z_{e,7}^\pm\subset Z_{13}^\pm$ & $x^3y \pm y^6 + axy^5$
              & $-$                             & $7$     & $5$     & 5/9 \\
 \hline
\end{tabular}
 \eqno (\thetheorem)
 $$
Table~(\ref {evKl}) is an extension of the table of even germs of even
multiplicity $\mu_e\le5$ from the diploma work of M.~Beer,
see~\cite[Satz~5.6]{Beer} or~\cite [Theorem~4.1]{Wass}.

\begin{note} \rm
A singularity $X_{e,5}$ is a (local) minimum only if its normal form is
$x^4+ax^2y^2+y^4$ where $a>-2$, $a\neq 2$.
Singularities $X_{e,r+3}$ and $Y_{e,r,s}$ are (local) minima only if the
normal form is $x^4 + x^2 y^2 + a y^{2r}$, or resp.\
$x^{2r} + ax^2 y^2 + y^{2s}$ with $a>0$.
A singularity $\widetilde{Y}_{e,r}$ is a (local) minimum only if the normal
form is $(x^2+y^2)^2+ay^{2r}$.
\end{note}

The result formulated above is proved in Section~\ref{sec:vers}.
Section~\ref {sec:versal} gives the definition and properties of even versal
deformations of germs, and explains their relation to generic parametric
families of smooth even functions.
The reader who is mainly interested in a description of bifurcations of critical
points for generic families of smooth even functions can pass to
Section~\ref{sec:bif}.

\section{Germs of even functions} \label{sec:vers}

\refstepcounter{theorem}
\subsection{Equivalence of germs}

Let $\A$ be the algebra of germs at zero of $C^\infty-$functions
$f:\RR^\nu\to\RR$, and let $\m \subset \A$ be the subalgebra of germs of
functions, whose value at zero vanishes.
Further $\A_e \subset \A$ is the subalgebra of germs of even functions, and
$\m_e=\m \cap \A_e$ is the maximal ideal of the subalgebra $\A_e$.
\\
Let $f \in \m_e$. Then all partial derivatives
$\frac{\partial f}{\partial x_i}$ are odd functions. Observe that
linear combinations of the form $h_1 \frac {\partial f}{\partial x_1}  +
\ldots + h_\nu \frac {\partial f}{\partial x_\nu}$ with
coefficients $h_k \in \A$, which are odd functions, form an ideal
$I^e_{\nabla f} \subset \m_e$ of the algebra $\A_e$.
This ideal is called the {\it even gradient ideal}, or the
{\it even Jacobian ideal}, of the germ $f$.

\begin{definition} \rm
The {\it even local algebra of the gradient map of $f$ at zero}
is the quotient algebra of the algebra of even germs, having a
vanishing value at zero, by the ideal which is generated by the
components of $\nabla f$:
\[
Q^e_{\nabla f}=\m_e \slash I^e_{\nabla f}.
\]
\end{definition}
{\bf Example.} Let $\nu=1$ and $f(x)=x^6$. Then
$Q^e_{\nabla f}=\{ a_2 x^2 + a_3 x^4 \}$.
(The divided powers algebra with degrees less than 5 and a vanishing
free term.)

\begin{definition} \label {def:mu} \rm
The integer $\mu_e = \dim Q^e_{\nabla f}+1$ is called
the {\it even local multiplicity} of the gradient map of
$f \in \A_e$ at zero, or
the {\it even multiplicity} of the critical point $0$ of the function
$f \in \A_e$
(compare~\cite[v.~1, 6.3]{Arnold}).
If $\mu_e < \infty$ then the critical point $0$ of the even function $f$
is called a {\it finite-even-multiple} critical point.
\end{definition}

Let $G$ be the group of germs of smooth odd changes of variables in $\RR^\nu$,
i.e.\ germs of diffeomorphisms $h:$
\[
h(0)=0, \quad h(x)=-h(-x).
\]
The change of variables by means of $h\in G$ obviously transforms even
functions to even ones, and odd functions to odd ones. Besides, under
such a change of variables, the ideal $I^e_{\nabla f}$ and the algebra
$Q^e_{\nabla f}$ are equivariantly transformed:
$I^e_{\nabla(h^*f)}=h^*I^e_{\nabla f}$ and
$Q^e_{\nabla(h^*f)}=h^*Q^e_{\nabla f}$. Therefore, the even multiplicity $\mu_e$
of the critical point $0$ is preserved under odd changes of variables in
$\RR^\nu$.

\begin{definition} \label {def:RO}
\rm
Elements $f_1,f_2 \in \A_e$ which admit a change of variables $h\in G$ such that
\[
f_1 \circ h^{-1} \equiv f_2
\]
are called {\it $R_O-$equivalent}, i.e.\ {\it right-odd equivalent}.
\end{definition}

Recall that (not necessarily even) function germs $f_1,f_2$ are called
$R-$equivalent if the equality $f_1 \circ h^{-1}\equiv f_2$ is valid for some,
not necessarily odd, diffeomorphism $h$ with $h(0)=0$.

The following analogues of the Morse Lemma, the splitting lemma, and
the Tougeron theorem about a finite-determination show that
the class of odd changes of variables is wide enough:

\begin{theorem}\label{Mors}
{\bf (<<Morse lemma>> for even functions)} In a neighbourhood
of the nondegenerate critical point $0$, an even function
is $R_O-$equivalent to the sum of a quadratic form and a constant.
\end{theorem}

The above theorem is a partial case of the so called equivariant Morse lemma,
see~\cite[v.~1,~17.3]{Arnold}.

\begin{theorem}\label{th:paramMors}
{\bf (Splitting lemma for even functions)}
In a neighbourhood of the critical point $0$ of co-rank $k\le\nu$,
an even function is $R_O-$equivalent to a function having the form
$\psi(x_1,\dots,x_k)+Q(x_{k+1},\dots,x_\nu)$ where $\psi$ is an even
function in $k$ variables, and $Q$ is a nondegenerate quadratic form in
$\nu-k$ variables.
\end{theorem}

Theorem~\ref {th:paramMors} generalizes Theorem~\ref {Mors}. It is proved by
means of the <<parametric Morse lemma>> for even functions.
Theorem~\ref {th:paramMors} is a partial case of the $\ZZ_2-$splitting
lemma~\cite{Beer}, which is a partial case of the equivariant
splitting lemma~\cite[Lemma 2.1]{Wass}.
For arbitrary smooth functions, a similar assertion is proved by J.~Mather in
his unpublished notes about $R-$equivalence. For the case $k=1$, a proof is
given in the beginning of the proof of the lemma from~\cite[v.~1, 9.6]{Arnold}.

\begin{theorem}\label{Tuazhron}
{\bf (Tougeron theorem about finite-determination of an even germ)}
In a neighbourhood of the finite-even-multiple critical point $0$,
an even function is $R_O-$equivalent to a polynomial (namely, to its Taylor
polynomial of degree $2\mu_e$ at $0$ where $\mu_e$ is the even multiplicity).
\end{theorem}

Actually, this theorem admits the following generalization:
one can choose the neighbourhood of $0$ in such a way that any even
function which is close enough to the function under consideration is also
$R_O-$equivalent in this neighbourhood to a polynomial of degree $2\mu_e$
(namely, to its Taylor polynomial of degree $2\mu_e$ at $0$).

A proof of this theorem (as well as the generalization from above) is analogous
to the proof of the Tougeron theorem about finite-determination of a usual
germ (for a proof by J.~Mather see~\cite [v.~1,~6.3 and~6.4]{Arnold}),
see~\cite {Poenaru} and~\cite {Beer}.
\\{\bf Example.} For a nondegenerate critical point ($\mu_e=1$),
an even function is $R_O-$equivalent to its Taylor polynomial of degree~2
(<<Morse lemma>>). This example shows that the degree $2\mu_e$ can not be
replaced by a lesser one.

\refstepcounter{theorem}
\subsection{Classification of singularities of even germs}

In order to obtain a classification of singularities of even smooth functions,
we can not unfortunately just consider all even normal forms of singularities
of the usual smooth functions, since we do not know whether the $R-$equivalence
of even germs with a vanishing value at $0$ implies their $R_O-$equivalence.
But we can transfer the corresponding classification assertions for the
singularities of usual germs (see~\cite[chap.~2]{Arnold}) to our case.

\begin{theorem} \label {th:cork}
The space of germs at $0$ of even functions in $\nu$ variables with the
critical point $0$ of co-rank $k \leq \nu$ has codimension $k(k+1)/2$ in
the space of germs of even functions at $0$.
(The co-rank of a critical point is the co-rank of the second differential of the
function at this point.)
\end{theorem}

In particular, in generic families of even functions with at most~5 parameters,
only singularities of co-rank at most~2 occur. Thus, due to
Theorem~\ref {th:paramMors}, in order to classify singularities and even
singularities of functions which occur in generic families with at most~5
parameters, one may restrict oneself to functions in two variables.

A proof of the above theorem follows from the fact that the codimension of
the set of matrices of co-rank $k$ in the space of all symmetric matrices
of order $\nu$ equals $k(k+1)/2$, see~\cite[v.~1,~2.2 and~11.1]{Arnold}.
Here one should take into account that the second differential
of an even function at zero can be an arbitrary symmetric matrix.

\begin{theorem}\label{odnokr}
Let $0$ be a finite-even-multiple critical point of co-rank~$1$ of an even
function. Then, in a neighbourhood of $0$, the function is $R_O-$equivalent to
the function $c+x^{2\mu_e}+Q$ where $c$ is a constant, $\mu_e$ is the even
multiplicity, and $Q$ is a nondegenerate quadratic form in the remaining
variables.
\end{theorem}

A proof is obtained by means of Theorem~\ref {th:paramMors}. It is analogous
to the proof of the similar assertion for arbitrary smooth functions, see the
lemma from~\cite[v.~1, 9.6]{Arnold} and Statement~\ref {st:defe}.

Thus, in order to describe all singularities in generic (i.e.\ stable,
see~\ref {note:vozmTh}) families of even functions with $l \leq 5$ parameters,
it remains to classify even singularities of co-rank $2$ of even functions in
two variables.
Let us show that any even singularity of interest has a nonvanishing
$4-$jet (i.e.\ its Taylor series contains a nonvanishing term of order $4$
or less).

\begin{theorem}
The space of germs at $0$ of even functions in two variables
with a vanishing $4-$jet has codimension $8$ in the space of germs
of even functions with critical value $0$ at the point $0$.
\end{theorem}

A proof of the above theorem follows from the invariance of the order of
the lowest term under a diffeomorphic change of variables.
Moreover, one should observe that the condition about the absence
of the second order terms is defined by 3 equations, see Theorem~\ref {th:cork},
while the condition about the absence of the fourth order terms is defined by 5
equations.

\refstepcounter{theorem}
\subsection{Even singularities of co-rank 2 with vanishing 2-jet and
nonvanishing 4-jet} \label{subsec:4forms}

Consider a real-valued 4-form in two variables
 $$
A(x,y)=a_1 x^4 + a_2 x^3y + a_3 x^2 y^2 + a_4 x y^3 + a_5 y^4,
\quad a_i \in \mathbb R ,
 $$
where at least one of the coefficients $a_i$ does not vanish.
Zeros of the equation $A(x,y)=0$ define four lines on the complex plane $(x,y)$.
Thus the form $A$ can be represented as the product of four linear forms
$\widetilde{A}_i=u_ix+w_iy$ with complex coefficients, $i=1,2,3,4$.
Combining the linear forms into pairs, one can obtain a representation
 $$
A(x,y)=A_1(x,y)A_2(x,y)
 $$
where $A_i(x,y)$, $i=1,2$ are quadratic forms with real coefficients.
We obtain the following cases of canonical forms of real-valued 4-forms in two
variables.
\medskip
\\ 1.
If the four lines $\widetilde{A}_i(x,y)=0$ are pairwise different, or
the lines of one of the pairs are not real and the lines of the other pair are
different, then the quadratic forms $A_i$ are nondegenerate. Moreover, they
can be simultaneously reduced to diagonal form by a real linear change of
variables (after a change of the grouping of the lines into pairs if necessary).
This reduces the form $A(x,y)$ to the form
 $$
A(x,y)=\pm x^4+ax^2y^2 \pm y^4, \quad a \in \mathbb R.
 $$
2.
If exactly two of the real lines coincide, for example,
$\widetilde{A}_1(x,y)=0 \Leftrightarrow \widetilde{A}_2(x,y)=0$,
then one of the quadratic forms $A_i(x,y)$ is nondegenerate,
and the second one is degenerate. Moreover, they can be simultaneously reduced
to diagonal form by a real linear change of variables. This reduces the form
$A(x,y)$ to the form
 $$
A(x,y) = x^2(\pm x^2\pm y^2).
 $$
3. If the four lines are real and the lines of any pair coincide, while the
lines of different pairs are different (after a change of the grouping of the
lines into pairs if necessary), for example
$\widetilde{A}_1(x,y)=0 \Leftrightarrow \widetilde{A}_2(x,y)=0$ and
$\widetilde{A}_3(x,y)=0 \Leftrightarrow \widetilde{A}_4(x,y)=0$, then,
after a real linear change of variables, the form $A(x,y)$ looks like
 $$
A(x,y)=\pm x^2y^2.
 $$
4. If exactly three (real) lines pairwise coincide then,
after a real linear change of variables, the form $A(x,y)$ looks like
  $$
A(x,y)=x^3y.
  $$
5. If all four (real) lines coincide then,
after a real linear change of variables, the form $A(x,y)$ obviously looks like
 $$
A(x,y) = \pm x^4.
 $$
We remark that all changes of variables from above are linear and therefore are
{\it odd}.

Let us notice at once that the case~5 is not of interest to us, since
the codimension of the subspace of even germs with a $4-$jet (i.e.\ the Taylor
polynomial of degree 4) of such type equals~6. Thus, such a singularity can not
appear in generic families of even functions with $l \leq 5$ parameters.
In fact, the subspace of such germs is defined by three equations
on the quadratic terms and three equations on the fourth order terms.

Consider the case 1. If $|a|=2$ and the signes $\pm\pm$ in the form $A(x,y)$
equal $++$ or $--$ then $A(x,y)=\pm(x^2 \pm y^2)^2$. Therefore, the form
$A(x,y)$ has the form $\pm x^2y^2$ (case~3), or $\pm(x^2+y^2)^2$
(the latter subcase will be considered after case~3).

\begin{theorem}\label{stFour}  {\bf (Even singularity $X_{e,5}^{\pm\pm}$)}
Let the Taylor polynomial of the $4$th order of an even function $f(x,y)$
at the point $0$ have the form $\pm x^4+ax^2y^2\pm y^4$ where $|a|\neq 2$
in the cases $++$ and $--$.
Then the germ of $f$ at zero is $R_O-$equivalent to its Taylor polynomial
of the $4$th degree.
\end{theorem}

A proof to this theorem is analogous to the proof of a similar assertion
for arbitrary smooth functions, see the theorem from~\cite[v.~1,~12.6]{Arnold}
about a normal form of semiquasihomogeneous singularities. Here one should use
the Tougeron Theorem~\ref{Tuazhron} for even functions.

Let us analyze case 2. We will use the Newton diagram (Fig.~1) which
corresponds to the Taylor series $\sum a_{p,q} x^p y^q$ of $f$ and is
defined as follows. This diagram represents the support $\supp f$ consisting of
those integer points $(q,p)\in\mathbb R^2$
which are exponents of the monomials contributing to the series with
nonvanishing coefficients.

\begin{figure}
\unitlength = 3mm
\begin{center}
\begin{picture}(13,11)(0,0) 
\put(0.1,0.1) {\vector(1,0){14}}
\put(0.1,0.1) {\vector(0,1){11}}
\put(14,-.7) {$q$}
\put(-.5,-1) {0}
\put(4,-1) {2}
\put(-.7,11) {$p$}
\put(-.7,4) {2}
\put(-.7,8) {4}
\multiput(-.1,-.1)(4,0) {3} {\multiput(0,0)(0,4) {3} {\,$\cdot$}}
\multiput(1.9,1.9)(4,0) {3} {\multiput(0,0)(0,4) {3} {\,$\cdot$}}
\qbezier[25](4.7,2.2)(5.1,2.3)(5.53,2.84) 
\put(5.53,2.83){\vector(2,3){0.2}}        
\thicklines
\path(0.2,8)(4.2,4)   
\path(4.2,4)(12.2,0)  
\put(7.5,5.5){$\Omega$} 
\put(-.1,7.8){$\bullet$} \put(3.9,3.8){$\bullet$}
\end{picture}
\end{center}

\vskip 0.5truecm
\centerline{{\bf Fig.~1:}
The ruler-method for an even germ $f$ with
$(0,4),\, (2,2)\in\supp f\subset \Omega$}
\end{figure}

The {\it method of a ruler turning} by Newton consists of the following steps.
At first, one draws a line (a <<ruler>>) through the exponent of a marked
monomial (which is $\pm x^2y^2$ in our case) such that this line separates the
origin from the unmarked points of the support. Later, one turns the ruler
around the marked exponent (counter-clockwise in our case) until it meets the
exponent of another monomial present.
\\ Subcase 1. The ruler meets just one point $(2k,0)$ of the $q$-axis. In this
case, one proves that all other points of the support do not affect the
$R_O-$equivalence class of the germ $f$. Namely, the germ reduces to
the form $\pm x^4+ax^2y^2\pm y^{2k}$ with $a\ne0$ (the even singularity
$X_{e,k+3}$), see the corollary from~\cite[v.~1,~12.7]{Arnold}.
Moreover, it follows directly from the proof (the cross-words method) that this
reduction can be fulfilled by means of an odd change of variables. Here one
should use the Tougeron Theorem~\ref{Tuazhron} for even functions.
\\ Subcase 2. The ruler meets two points $(k+1,1)$ and $(2k,0)$ which
correspond to terms $xy^{k+1}$ and $y^{2k}$ with $k$ even.
In this case, let us consider the polynomial which is defined by the monomial
$\pm x^4$ and the points falling on the ruler:
 $$
Ax^4 + Bx^2y^2 + C xy^{k+1} + Dy^{2k} =
Ax^4 + y^2 \left [ Bx^2 + Cxy^{k-1} + Dy^{2(k-1)} \right ] .
 $$
We will call this polynomial the {\it mean part} of the germ being investigated.
It is easy to observe that, by means of an {\it odd} change of variables
$x'=x+\lam y^{k-1}$, one can turn the coefficient $C$ into zero.
Now the subcase~2 is reduced to the subcase~1.

Case 3 is analyzed similarly to case~2 (the ruler turns clockwise at
first and counter-clockwise later). Here, by means of an odd change of
variables, the germ reduces to the form $\pm x^{2k}+ax^2y^2\pm y^{2l}$ with
$a\ne0$ (the even singularity $Y_{e,k,l}$).

Let us return to case 1. It remains to analyze the subcase where the 4-jet
of the even germ $f$ reduces to the form $\pm(x^2+y^2)^2$. Let us associate to
the even Taylor series $\pm(x^2+y^2)^2 + \sum_{p+q\ge 6} a_{p,q} x^p y^q$ of $f$
the series of complex numbers $a_3,a_4,\dots$, which are defined as
$a_k=\sum_{p=0}^{2k} i^p a_{p,2k-p}$. In particular,
$a_3=a_{0,6}-a_{2,4}+a_{4,2}-a_{6,0}+i(a_{1,5}-a_{3,3}+a_{5,1})$.
Using the Tougeron Theorem~\ref{Tuazhron} for even functions,
it is easy to check that, by means of an odd change of variables, the germ
reduces to the form $\pm(x^2+y^2)^2 + ay^6$ with $a=|a_3|$ as soon as
$a_3\ne0$ (the even singularity $\tilde Y_{e,3}$).
Moreover, one easily shows that there exists a sequence of polynomials
$\tilde a_3=a_3,\tilde a_4,\dots$ in $a_{p,q}$ with complex coefficients such
that, by means of an odd change of variables, the germ $f$ reduces to the form
$\pm(x^2+y^2)^2+ay^{2k}$ with $a=|\tilde a_k|$ as soon as $k$ is the index of
the first nonvanishing element of this sequence (the even singularity
$\tilde Y_{e,k}$). Here the polynomial $\tilde a_k-a_k$ explicitly depends only
on $a_{p,q}$, $6\le p+q<2k$.
These properties of the sequence $a_3,\tilde a_4,\tilde a_5,\dots$ can be
considered as an analogue of the method of a ruler turning, see above.
Let us notice that we actually do not need to consider the case $a_3=0$,
since the codimension of the corresponding subspace of the space of even
singularities having value $0$ at the origin equals~7. Here, one considers the
subspace consisting of those singularities which satisfy the following
condition:
the 4-jet of the singularity reduces to the form $\pm(x^2+y^2)^2$ by means of
an odd change of variables and, in the obtained variables, $a_3=0$.

Case 4.
It is easy to show that, by means of an odd change of variables, the 6-jet
of the germ reduces to the form $x^3y+a_1y^6+a_2xy^5$. Moreover, the germ itself
reduces to the form $x^3y\pm y^6+axy^5$ if $a_1\ne0$ (the even singularity
$Z_{e,7}$).
Consider the space of even germs in two variables such that the Taylor
polynomial of the fourth order reduces to the form $x^3y$ by means of an odd
change of variables, and the coefficient at the term $y^6$ (with respect to
these variables) vanishes.
It is easy to check that the validity of this condition does not depend on the
choice of variables. The codimension of this subspace in the space of all
even germs taking value $0$ at the origin equals~6.
In fact, this subspace is defined by
three equations on the 2nd order terms of the Taylor series,
two equations on the 4th order terms, and
one equation on the coefficient at the term $y^6$.

\section {Even deformations of even functions} \label {sec:versal}

\begin{definition} \label{def:evdef} \rm
An {\it even deformation} of an even germ $f\in \m_e$ with a {\it base} $\Lambda
\subset \mathbb R^l$ is the germ of a smooth map
$F:(\mathbb R^\nu \times \Lambda,0) \rightarrow (\mathbb R,0)$
such that $F(x,0)\equiv f(x)$ and $F(\cdot,\lam) \in \m_e$ for an
arbitrary value of the parameter $\lam \in \Lambda$.
\end{definition}

\begin{definition} \rm
An even deformation $F'$ is {\it $R_O-$equivalent} to an even deformation $F$ if
\[
F'(x,\lam)\equiv F(g(x,\lam),\lam)
\]
where $g:(\mathbb R^\nu \times\mathbb R^l,0) \rightarrow (\mathbb
R^\nu,0)$ is a smooth germ which is odd with respect to $x$ and satisfies the
condition $g(x,0)=x$.
\end{definition}

\begin{definition} \rm
A deformation $F'$ is {\it induced} from $F$ if
\[
F'(x,\lam')\equiv F(x,\varphi(\lam'))
\]
where $\varphi : (\mathbb R^l,0) \rightarrow (\mathbb R^l,0)$ is the germ of a
smooth map between bases.
\end{definition}

\begin{definition} \label{def:ROvers} \rm
An even deformation $F$ of an even germ $f \in \m_e$ is called {\it
$R_O-$versal} if any even deformation of $f$ is $R_O-$equivalent to
a deformation which is induced from $F$. Therefore, an even deformation
$F$ of $f$ is $R_O-$versal if any even deformation $F'$ of $f$ can be
represented as
 $$
\refstepcounter{theorem}
\label{indT}
F'(x,\lam')\equiv F(g(x,\lam'),\varphi(\lam')), \quad g(x,0)=x,\
\varphi(0)=0,\ g(-x,\lam)=-g(x,\lam).
 \eqno (\thetheorem)
 $$
\end{definition}

With regard to~\ref {subsec:vozmTh}, the following theorems provide
a means for the investigation of singularities of even germs and caustic
surfaces.

\begin{theorem}\label{verMO}
Any $R_O-$versal even deformation of a finite-even-multiple even germ $f$
has at least $\mu_e-1$ parameters, where $\mu_e$ is the even multiplicity
of $f$. An even deformation
\[
F(x,\lam)=f(x)+\lam_1 e_1(x)+\ldots + \lam_{\mu_e-1}
e_{\mu_e-1}(x)
\]
is $R_O-$versal if the germs of functions $e_k$ at zero define a basis
of the linear space $Q^e_{\nabla f}$.
\end{theorem}

\begin{theorem}\label{verMT}
Any $l-$parameter $R_O-$versal even deformation of an even germ $f$ is
$R_O-$equivalent to an even deformation which is induced from any other
$R_O-$versal even deformation with $l$ parameters by means of a diffeomorphism
of bases.
\end{theorem}

Proofs are similar to the proofs of theorems
from~\cite [v.~1,~8.2 and~8.3]{Arnold} and~\cite [v.~1, 8.5]{Arnold}, resp.
Actually, Theorems~\ref {Tuazhron},~\ref {verMO}, and~\ref {verMT} are partial
cases of the corresponding equivariant theorems. The equivariant theorems can
be proved in the same way as their non-equivariant analogues, since the groups
of germs of equivariant changes of variables in $\RR^\nu$ satisfy the
conditions of J.~Damon~\cite[3.2.4~and~3.2.5]{Itogi}.

\refstepcounter{theorem}
\subsection{Stable families of even functions} \label{subsec:vozmTh}

Let us make an important remark relating the theory of germs deformations
and the perturbation theory for parametric families of functions.

Let $U$ be a domain of $\RR^\nu$ which contains $0$ and is invariant under
the transformation $x\mapsto-x$, and let $\Lambda$ be a domain of $\RR^l$ which
contains $0$. A smooth function $F\colon\, U\times\Lambda\to\RR$ is called
a {\it smooth family of even functions on $U$ with the base $\Lambda$} if
$F(\cdot,\lam) \in \m_e$ for any value of the parameter $\lam\in\Lambda$,
compare~\ref {def:evdef}.
Any representative $F\colon\,(U\times\Lambda,0)\to(\RR,0)$ of an even
deformation of an even germ $f\in\m_e$ is obviously a smooth family of even
functions.

\begin{note} \label{note:vozmTh} \rm
The {\it $R-$versality of a deformation} is equivalent
to the stronger property {\it $R-$stability of a deformation},
see Remark~3 from~\cite[v.~1, 8.4]{Arnold}
(compare the definition from~\cite[v.~1, 21.1]{Arnold}).
Similarly, the {\it $R_O-$versality of an even deformation} is equivalent
to the stronger property {\it $R_O-$stability of an even deformation}:
any representative $F\colon\,(U\times\Lambda,0)\to(\RR,0)$ of an even
deformation $F$ admits a neighbourhood $E$ in the space of smooth families
of smooth even functions on $U$ with the base $\Lambda$
(i.e.\ smooth functions on $U\times\Lambda$ satisfying the conditions
$F(-x,\lam)=F(x,\lam)$ and $F(0,\lam)=0$) such that,
for any family $F'\colon\,U\times\Lambda\to\RR$ from $E$,
there exists a point $0'$ such that the germ $F'$ at $0'$ defines
a deformation of a germ $f'$, $R_O-$equivalent to the deformation $F$,
where $f'$ is $R_O-$equivalent to $f$.
Moreover, if the family $F'$ is close enough to $F$, the point $0'$
can be chosen arbitrarily close to $0$ and the $R_O-$equivalence
can be chosen arbitrarily close to the identity.
\end{note}

The versality and the stability of a deformation are apparently equivalent
to the stronger property {\it strict stability}, which is formulated as
follows in the case of even deformations:
any representative $F\colon\,(U\times\Lambda,0)\to(\RR,0)$ of an even
deformation $F$ admits a neighbourhood $E$ in the space of smooth families
of smooth even functions on $U$ with the base $\Lambda$
and a neighbourhood $U'\times\Lambda'$ of the origin in $U\times\Lambda$
such that $U'$ is invariant under the transformation
$(x,\lam)\mapsto(-x,\lam)$ and,
for any family $F'\colon\,U\times\Lambda\to\RR$ from $E$,
there exist mappings $g\colon\, U'\times\Lambda'\to U\times\Lambda$ and
$\varphi\colon\,\Lambda'\to\Lambda$, close to the identities, such that
\[
F(x,\lam)\equiv F'(g(x,\lam),\varphi(\lam)),\quad g(-x,\lam)=-g(x,\lam)
\]
in the domain $U'\times\Lambda'$. We remark that the conditions
$g(x,0)=x$ and $\varphi(0)=0$ are not required, compare~(\ref{indT}).

The property $R^+-$stability of an $R^+-$versal deformation (resp.\
$R_O-$stability of an $R_O-$versal deformation) implies that a singularity
of multiplicity (resp.\ even multiplicity) $\mu$ is non-removable for generic
$l-$parameter families of functions with $l\ge\mu-1$.
Moreover, the caustic surfaces for such families have the same local structure
as the caustic surfaces for $l-$parameter $R^+-$versal
(resp.\ $R_O-$versal) deformations.

\begin{note} \label{note:strat} \rm
$R^+-$versal deformations of all known types of finite-multiple singularities
are {\it $R^+-$stratified} in the following sense.
If $F\colon\,(U\times\Lambda)\to(\RR,0)$ is a representative of an
$R^+-$versal deformation of a germ $f\in\A$
then there exists a neighbourhood $U'\times\Lambda'$ of the origin in
$U\times\Lambda$ with the following properties.
For any $(x,\lam)\in U'\times\Lambda'\setminus\{0\}$, the multiplicity of the
germ $F(\cdot,\lam)$ at $x$ is not larger than the multiplicity $\mu$ of the
germ $f$ at $0$.
Furthermore, the points $(x,\lam)\in U'\times\Lambda'$ corresponding to germs
which are $R^+-$equivalent to the germ $f$ form a smooth submanifold $M_1$
of codimension $\mu-1$, while the points $(x,\lam)\in U'\times\Lambda'$
corresponding to germs of multiplicity $\mu$ form a smooth submanifold
$M_2\subset M_1$.
Moreover, the restriction of the projection $(x,\lam)\mapsto\lam$ to any of
the submanifolds $M_1$ and $M_2$ is a diffeomorphism onto its image.
Here the image of the submanifold $M_2$ is called a {\it $\mu=\const$ stratum},
its codimension $c$ is called the {\it codimension} of the singularity of $f$,
and the integer $m=\mu-c-1$ is called the {\it modality} of the singularity of
$f$.

Similarly, $R_O-$versal deformations of all known types of finite-even-multiple
singularities of even germs are {\it $R_O-$stratified} in the following sense.
If $F\colon\,(U\times\Lambda)\to(\RR,0)$ is a representative of an
$R_O-$versal deformation of an even germ $f\in\m_e$ then
there exists a neighbourhood $U'\times\Lambda'$ of the origin in
$U\times\Lambda$ such that $U'$ is invariant under the transformation
$(x,\lam)\mapsto(-x,\lam)$ and the following conditions are fulfilled.
For any $\lam\in\Lambda'$,
the even multiplicity of the germ $F(\cdot,\lam)$ at $0$ is not larger than
the even multiplicity $\mu_e$ of the germ $f$ at $0$.
Moreover, the points $\lam\in\Lambda'$ corresponding to even germs
which are $R_O-$equivalent to the even germ $f$ form a smooth submanifold
of codimension $\mu_e-1$, while the points $\lam\in\Lambda'$
corresponding to even germs of even multiplicity $\mu_e$ form a smooth
submanifold called an {\it even $\mu_e=\const$ stratum}.
Here the codimension $c_e$ of the even stratum is called the {\it even
codimension} of the even singularity of $f$, and the integer $m_e=\mu_e-c_e-1$
is called the {\it even modality} of the even singularity of $f$.
\end{note}

If a singularity of multiplicity $\mu$ admits an $R^+-$versal deformation
which is $R^+-$stratified and strictly $R^+-$stable then the class of
singularities containing this singularity appears in typical $l-$parameter
families of functions if and only if $l\ge c=\mu-m-1$.

Similarly, if an even singularity of even multiplicity $\mu_e$ admits an
$R_O-$versal deformation which is $R_O-$stratified and strictly $R_O-$stable
then the even class of singularities containing this singularity appears in
typical $l-$parameter families of even functions if and only if
$l\ge c_e=\mu_e-m_e-1$.
\\
{\bf Example.} Consider an even germ $f \in \m_e$ at the origin with a
singularity of type $A_{e,4}$, see Table~(\ref{evKl}).
Theorem~\ref {verMO} implies that,
as an even $R_O$-versal deformation of $f$, one can take
\[
F(x,\lam)=\pm x^8+\lam_1 x^6 + \lam_2 x^4 + \lam_3 x^2.
\]
Thus, an even germ of type $A_{e,4}$ becomes non-removable under small
perturbations of a family if the dimension $l$ of the parametric space is
at least three. Moreover, the surface of caustic values corresponding
to singularities of this type has codimension $l-3$,
see~\ref {note:strat}.
\\ {\bf Example.} Consider an even germ $f \in \m_e$ at the origin with the
even singularity $X_{e,5}^{++}$, see Table~(\ref{evKl}).
Theorem~\ref {verMO} implies that,
as an even $R_O$-versal deformation of $f$, one can take
\[
F(x,y;\lam)= x^4+ax^2y^2 + y^4 + \lam_1 x^2 + \lam_2 xy + \lam_3 y^2 +
\lam_4 x^2y^2 .
\]
Thus, even germs of class $X_{e,5}^{++}$ become non-removable under small
perturbations of an even family if the dimension $l$ of the parametric space is
at least four. Observe however that if $\lam_1=\lam_2=\lam_3=0$ then, for
small values of $\lam_4$, the even germ $F(\cdot;\lam)$ at zero has a
singularity of class $X_{e,5}^{++}$, although it is not equivalent to its
initial form $f=F(\cdot;0)$. Thus, if we want an even singularity of class
$X_{e,5}^{++}$ to be preserved under arbitrary small perturbations of an even
family, it is enough to require that the dimension $l$ of the parametric space
is at least three. Such singularities are called {\it modal} (in our case,
the even singularity of $f$ is unimodal, since it contains one even module $a$).
This term is used for (even) germs whose arbitrary small neighbourhood in the
space of (even) germs is covered by a finite number of $m-$parameter families
of (even) singularities. The minimal integer $m$ with this property is called
the (even) {\it modality} of the singularity~\cite [v.~1,~II]{Arnold}.

\section{Bifurcations of critical points of even functions in typical
families} \label{sec:bif}

Consider a family $\omega_\lam(k) : \torsp \rightarrow \mathbb R$
satisfying~(\ref{fam}).
Values of the parameter $\lam$ run through some domain in $\param$.
We will study {\it generic} families of functions $\omega_\lam(k)$,
that is such families that all critical points of $\omega_\lam(\cdot)$ are
non-removable for the entire family under small perturbations of the family.

For an arbitrary dimension $l$ of the parametric space,
there are only three types of caustic surfaces of codimension one.
\\ 1. A caustic hypersurface of type $A_{e,2}^+\subset A_3$ or
$A_{e,2}^-\subset A_3$. It corresponds to a singularity at a basic critical
point. Under a passage of the parameter across this hypersurface, a pair
$(k,-k)$ of additional critical points separates from the basic critical point.
\\2. A caustic hypersurface of type $A_2$. It corresponds to a singularity at a
twin critical point. Under a passage of the parameter across this hypersurface,
two pairs of additional critical points (existing on one side of the
hypersurface) merge into one pair and disappear (i.e.\ are absent on the other
side of the hypersurface).

Below, we will give a detailed description of the behaviour (bifurcations)
of critical points near the caustics for small (physical) dimensions $l=1,2,3$
of the parametric space. Bifurcations for the codimensions $l=1,2$ are described
in the works~\cite{MamatovMinlos, LakshtanovMinlos}.

In Sections~\ref {sec:def} and~\ref {sec:cond:versal}, we will show how to
determine types of critical points of functions which appear in typical
parametric families with a small number $l=1,2,3,4$ of parameters, and
formulate conditions on their deformations to be typical
(more precisely, versal).

\refstepcounter{theorem}
\subsection{Typical one-parameter families}

Although we considered even germs taking value 0 at the origin, the values
of an even function at the basic critical points freely depend on $\lam$. For
convenience, we will assume below that the values of even functions at
the basic critical points vanish.

There exists a finite number of caustic values $\lam_1,\ldots,\lam_p$,
for which a {\it basic} critical point degenerates. Due to Table~(\ref{evKl})
and Statement~\ref {st:defe}, the singularity types at these points are
$A_{e,2}$.
There exists a finite number of caustic values $\bar\lam_1,\ldots,\bar\lam_q$,
for which {\it additional} (twin) critical points degenerate. Due to
Table~(\ref{obKl}) and Statement~\ref {st:def}, the singularity type at
these points is $A_2$.

\begin{note} \rm
It is easy to observe that an additional critical point which is a local
extremum (i.e.\ a point of local minimum or local maximum) does not degenerate.
\end{note}

By analyzing (even) versal deformations for degenerate critical points, we
understand what happens under a passage of the parameter across a caustic value:
\\ {\bf 1.} ($A_{e,2}^\pm$)
An even versal deformation of a function in a neighbourhood of a degenerate
{\it basic} critical point $0$ has the form
$F(k,\lam)=\pm k_1^4+(\lam-\lam_i) k_1^2 \pm k_2^2 \pm \dots \pm k_\nu^2$,
see Statement~\ref {st:versale}.
Under a passage of the parameter across the caustic value $\lam_i$, a pair
$(k,-k)$ of additional critical points separates from the basic critical point.
\\ {\bf 2.} ($A_2$)
An $R^+-$versal deformation of a function in a neighbourhood of a degenerate
{\it additional} critical point $\bar k=(\bar k_1,\dots,\bar k_\nu)$
has the form $F(k,\lam)=
(k_1-\bar k_1)^3+(\lam-\bar\lam_i)(k_1-\bar k_1)\pm k_2^2\pm\dots\pm k_\nu^2$,
see Statement~\ref {st:versal}.
Under a passage of the parameter across the caustic value $\bar\lam_i$,
two twin critical points $(\bar k',-\bar k')$ and $(\bar k'',-\bar k'')$ which
exist for values of $\lam$ lying on one side of $\bar\lam_i$ (namely, for
$\lam<\bar\lam_i$) merge into a twin critical point $(\bar k,-\bar k)$
and disappear for $\lam$ lying on the other side of $\bar\lam_i$ (namely,
for $\lam>\bar\lam_i$).
Bifurcations of level lines of a function for the case $\nu=2$ are illustrated
in Fig.~45~\cite[v.~1, 9.6]{Arnold}.

\refstepcounter{theorem}
\subsection{Typical two-parameter families} 

For $l=2$, all caustic values form a one-dimensional (not necessarily connected)
curve which has one of the following forms in a neighbourhood of any of its
points:
\\ {\bf 1.} ($A_2$)
A one-dimensional curve corresponding to a degeneracy of an {\it additional}
critical point with singularity type $A_2$.
Under a passage of the parameter across this curve, two pairs of critical
points merge into one pair $(k_0,-k_0)$ and disappear.
\\ {\bf 2.} ($A_{e,2}^\pm$)
A one-dimensional curve corresponding to a degeneracy of a {\it basic}
critical point with singularity type $A_{e,2}\subset A_3$.
Under a passage of the parameter across this curve, a pair
$(k,-k)$ of additional critical points separates from the basic critical point.
\\ {\bf 3.} ($A_3^\pm$)
In a neighbourhood of a caustic value $\overline{\lam}$ corresponding to a
degeneracy of an {\it additional} critical point with singularity type $A_3$,
the set of caustic values consists of two curves $\Gamma_1,\Gamma_2$ of type
$A_2$ (corresponding to a degeneracy of different twin critical points) which
adjoin the point $\overline{\lam}$ (a <<point of regression>>, or a
<<cusp point>>) and are tangent to each other at this point with contact
order $3/2$ (a <<semicubical parabola>>), see Fig.~2.

\begin{figure}
\unitlength = 5mm
\begin{center}
\begin{picture}(23,5.5)(-4,-3.5) 
\put(0,0) {\vector(1,0){2}}
\put(0,0) {\vector(0,1){2}}
\put(2.2,-.5) {$\lam_1$}
\put(.3,1.8) {$\lam_2$}
\thicklines
\qbezier[100](-4,2)(-1.5,0.1)(0,0) \put(-4.7,1.5){$\Gamma_1$}
\qbezier[100](-4,-2)(-1.5,-0.1)(0,0) \put(-4.7,-1.6){$\Gamma_2$}
\put(-2,1){\tiny$A_2$}
\put(-2,-1.2){\tiny$A_2$}
\put(-.2,-.5){\tiny$A_3$} \put(-.18,-.17){$\bullet$}

\put(0,-3.5){\large a}
\put(15,0) {
\thicklines
\qbezier[100](-4,2)(-1.5,0.1)(0,0) 
\qbezier[100](-4,-2)(-1.5,-0.1)(0,0) 
\thinlines
\qbezier[25](-5,0)(-2.5,0)(0,0)
\put(-2,1){\tiny$A_2$}
\put(-2,-1.2){\tiny$A_2$}
\put(-.2,-.5){\tiny$A_3$} \put(-.18,-.17){$\bullet$}

\put(.5,1.5){I}
\put(2.2,0){\circle{1.5}
  \qbezier[70](-.5,.55)(0,-.8)(.5,.55)
 }
\put(-3,-.3){II}
\put(-4.3,1){\circle{1.5}
  \qbezier[70](-.5,.55)(-.3,-.3)(-.15,0)
  \qbezier[70](-.15,0)(0,.5)(.1,-.2)
  \qbezier[70](.5,.55)(.3,-.8)(.1,-.2)
 }
\put(-4.3,-1){\circle{1.5}
  \qbezier[70](-.5,.55)(-.3,-.8)(-.1,-.2)
  \qbezier[70](-.1,-.2)(0,.5)(.15,0)
  \qbezier[70](.5,.55)(.3,-.3)(.15,0)
 }
\put(0,-3.5){\large b}
}
\end{picture}
\end{center}

\vskip 0.5truecm
\centerline{
\begin{tabular}{ll} \bf Fig.~2
& {\bf a:} A semicubical cusp point --
           a typical singularity $A_3^\pm$ of caustic curves \\
& {\bf b:} Bifurcations of functions at a singularity $A_3^+$
           in typical families of functions
\end{tabular}
}
\end{figure}

For $\nu=1$ and $\overline{\lam}=0$, an $R^+-$versal deformation of a function
in a neighbourhood of a degenerate additional critical point $\bar k$ of type
$A_3^+$ has the form
$F(k,\lam)=(k-\bar k)^4 +\lam_1(k-\bar k)^2+\lam_2(k-\bar k)$.
The caustic curve is defined by the equation
$\lam_2^2=(\frac{2}{3})^3(-\lam_1)^3$.

Let us describe the behaviour of critical points for those values of the
parameter $\lam$ which are close to the point $\overline{\lam}$:
In the domain $I$, there is one twin critical point $(k_1,-k_1)$.
Under a passage across the curve $\Gamma_1$ from the domain $I$ to
the domain $II$, two twin critical points $(k_1',-k_1')$ and $(k_2',-k_2')$
arise.
Under a passage across the curve $\Gamma_2$ from the domain $II$ to
the domain $I$, the points $(k_1,-k_1)$ and $(k_1',-k_1')$ merge and disappear,
and the point $(k_2',-k_2')$ is renamed $(k_1,-k_1)$.
\\ {\bf 4.} ($A_{e,3}^\pm$)
In a neighbourhood of a caustic value $\overline{\lam}$ corresponding
to a degenerate {\it basic} critical point of type $A_{e,3}\subset A_5$,
the set of caustic values consists of two curves $\Gamma_1,\Gamma_2$ of type
$A_{e,2}\subset A_3$ (corresponding to a degeneracy of the basic critical point)
which adjoin $\overline{\lam}$ and form a regular curve
$\Gamma_1\cup\overline{\lam}\cup\Gamma_2$, as well as one curve
$\Gamma'$ of type $A_2$ (corresponding to a degeneracy of a twin critical point)
which also adjoin $\overline{\lam}$ and is tangent to the curve $\Gamma_2$
with the contact order $2$, see Fig.~3.

\begin{figure}
\unitlength = 5mm
\begin{center}
\begin{picture}(27,8)(-6,-3) 
\put(0,0) {\vector(1,0){2}}
\put(0,0) {\vector(0,1){2}}
\put(1.8,.2) {$\lam_1$}
\put(.3,1.8) {$\lam_2$}
\thicklines
\path(-5,0)(5,0) \put(-5.5,.3){$\Gamma_2$} \put(5,0.2){$\Gamma_1$}
\qbezier[150](-3,4)(-1.5,0.1)(0,0) \put(-3.7,3.8){$\Gamma'$}
\put(-1.9,2){\tiny$A_2$}
\put(-3,-.5){\tiny$A_{e,2}$}
\put(3,-0.5){\tiny$A_{e,2}$}
\put(-.2,-.6){\tiny$A_{e,3}$} \put(-.18,-.17){$\bullet$}

\put(0,-3.5){\large a}
\put(15,0) {
\thicklines
\path(-5,0)(5,0)
\qbezier[150](-3,4)(-1.5,0.1)(0,0)
\thinlines
\qbezier[30](-5,2.5)(-2,0.05)(0,0)
\put(-1.9,2){\tiny$A_2$}
\put(-3,-.5){\tiny$A_{e,2}$}
\put(3,-0.5){\tiny$A_{e,2}$}
\put(-.2,-.6){\tiny$A_{e,3}$} \put(-.18,-.17){$\bullet$}

\put(2.5,2){I}
\put(1,2){\circle{1.5}
  \path(-.7,0)(.7,0)
  \put(0,.5){\vector(0,1){.2}} \path(0,-.6)(0,.6)
  \qbezier[70](-.5,.55)(0,-.5)(.5,.55)
 }
\put(1.2,-1.5){II}
\put(-1,-1.5){\circle{1.5}
  \path(-.7,0)(.7,0)
  \put(0,.5){\vector(0,1){.2}} \path(0,-.6)(0,.6)
  \qbezier[70](-.5,.55)(-.3,-.5)(-.15,-.2)
  \qbezier[70](-.15,-.2)(0,.2)(.15,-.2)
  \qbezier[70](.5,.55)(.3,-.5)(.15,-.2)
 }
\put(-3,.6){III}
\put(-5,1){\circle{1.5}
  \path(-.7,0)(.7,0)
  \put(0,.5){\vector(0,1){.2}} \path(0,-.6)(0,.6)
  \qbezier[80](-.5,.55)(-.4,-.5)(-.25,0)
  \qbezier[50](-.25,0)(-.18,.3)(-.1,.1)
  \qbezier[50](-.1,.1)(0,-.1)(.1,.1)
  \qbezier[50](.25,0)(.18,.3)(.1,.1)
  \qbezier[80](.5,.55)(.4,-.5)(.25,0)
 }
\put(-3.7,2.9){\circle{1.5}
  \path(-.7,0)(.7,0)
  \put(0,.5){\vector(0,1){.2}} \path(0,-.6)(0,.6)
  \qbezier[50](-.5,.55)(-.4,0)(-.3,.2)
  \qbezier[50](-.3,.2)(-.2,.4)(-.12,.2)
  \qbezier[50](-.12,.2)(0,-.2)(.12,.2)
  \qbezier[50](.3,.2)(.2,.4)(.12,.2)
  \qbezier[50](.5,.55)(.4,0)(.3,.2)
 }
\put(0,-3.5){\large b}
}
\end{picture}
\end{center}

\vskip 0.5truecm
\centerline{
\begin{tabular}{ll} \bf Fig.~3
& {\bf a:} A typical singularity $A_{e,3}^\pm$ of caustic curves \\
& {\bf b:} Bifurcations of functions at a singularity $A_{e,3}^+$
           in typical families of even functions
\end{tabular}
}
\end{figure}

For $\nu=1$ and $\overline{\lam}=0$, an $R^+-$versal deformation of a function
in a neighbourhood of a degenerate basic critical point $0$ of type
$A_{e,3}^+$ has the form $F(k,\lam)=k^6+\lam_1k^4+\lam_2k^2$.
The branches $\Gamma_1\cup\overline{\lam}\cup\Gamma_2$ and $\Gamma'$
of the caustic curve are defined by the equations $\lam_2=0$ and resp.\
$\lam_2=\frac{1}{3}\lam_1^2$, $\lam_1<0$.

Let us describe the behaviour of critical points for those values of the
parameter $\lam$ which are close to the point $\overline{\lam}$:
In the domain $I$, there are no twin critical points.
Under a passage across the caustic $\Gamma_1$ to the domain $II$,
a twin critical point $(k_1,-k_1)$ separates from the basic critical point.
Then, under a passage across the caustic $\Gamma_2$ to the domain $III$,
one more twin point $(k_2,-k_2)$ separates.
Finally, under a passage of $\lam$ from the domain $III$ to the domain $I$
across the caustic $\Gamma'$, both twin points merge and disappear.

Of course, except for the singularities listed above, transversal intersections of
different branches of caustics are also possible:
\\ {\bf 5, 6, 7.} ($A_2+A_2$, $A_2+A_{e,2}^\pm$, $A_{e,2}^\pm+A_{e,2}^\pm$)
In a neighbourhood of a caustic value $\overline{\lam}$ corresponding
to a degeneracy of {\it two different} (basic or additional) critical points
with singularity types $A_2$ or $A_{e,2}^\pm$,
the set of caustic values is the union of two curves $\Gamma_1,\Gamma_2$
(corresponding to a degeneracy of different critical points) which
intersect each other transversally at the point $\overline{\lam}$.

\refstepcounter{theorem}
\subsection{Typical three-parameter families}

All caustic values form a two-dimensional (not necessarily connected) surface
which has one of the following forms in a neighbourhood of any of its points:
\\ {\bf 1.} ($A_{e,2}^\pm$)
A surface corresponding to a degeneracy of a {\it basic}
critical point with singularity type $A_{e,2}\subset A_3$.
Under a passage of the parameter across this surface, a twin critical point
$(k,-k)$ separates from the basic critical point.
\\ {\bf 2.} ($A_2$)
A surface corresponding to a degeneracy of an {\it additional}
critical point with singularity type $A_2$.
Under a passage of the parameter across this surface, two twin critical
points $(k_1,-k_1)$ and $(k_2,-k_2)$ merge into one twin point $(k_0,-k_0)$
and disappear.
\\ {\bf 3.} ($A_3^\pm$)
Two caustic surfaces of type $A_2$ corresponding to a degeneracy of an
{\it additional} point and adjoining a curve of type $A_3$ (an <<edge of
regression>>, or a <<cuspidal edge>>, see Fig.~55~\cite[v.~1,~21.3]{Arnold}).
Bifurcations of critical points happen similarly to those in the case of
two-parameter even families, see~($A_3$).
\\ {\bf 4.} ($A_{e,3}^\pm$)
Three caustic surfaces of types $A_2,A_{e,2},A_{e,2}\subset A_3$
adjoining a common curve of type $A_{e,3}\subset A_5$ which corresponds to
a degeneracy of a {\it basic} point.
Bifurcations of critical points happen similarly to those in the case of
two-parameter even families, see~($A_{e,3}$).
\\ {\bf 5.} ($A_4$)
In a neighbourhood of a caustic value $\overline{\lam}$ corresponding to a
degeneracy of an {\it additional} critical point with singularity type $A_4$,
the set of caustic values has the form shown on
Fig.~56~\cite[v.~1,~21.3]{Arnold} (the <<swallow-tail>>).
Bifurcations of functions are shown in Fig.~28~\cite {Arnold2},
see also Fig.~55 and~97~\cite {Arnold2}).
\\ {\bf 6.} ($D^+_4$)
In a neighbourhood of a caustic value $\overline{\lam}$ corresponding to a
degeneracy of an {\it additional} critical point with singularity type $D^+_4$,
the set of caustic values has the form shown in
Fig.~57~\cite[v.~1,~21.3]{Arnold} (the <<purse>>).
\\ {\bf 7.} ($D^-_4$)
In a neighbourhood of a caustic value $\overline{\lam}$ corresponding to a
degeneracy of an {\it additional} critical point with singularity type $D^-_4$,
the set of caustic values has the form shown in
Fig.~57~\cite[v.~1,~21.3]{Arnold} (the <<pyramid>>).
\\ {\bf 8.} ($A_{e,4}^\pm$)
In a neighbourhood of a caustic value $\overline{\lam}$ corresponding
to a degeneracy of a {\it basic} critical point with singularity type
$A_{e,4}\subset A_7$, the set of caustic values is formed by the following
surfaces (Fig.~4, see also Fig.~62~\cite{Arnold}):

\begin{figure}
\unitlength = 5mm
\begin{center}
\begin{picture}(8,8)(-4,-4) 
\thinlines
\put(0,0) {\vector(2,1){1.5}}
\put(0,0) {\vector(-2,1){1.5}}
\qbezier[12](0,0)(0,-.75)(0,-1.5) \put(.02,-1.5){\vector(0,-1){.2}}
\put(1.2,0.9) {\tiny$\lam_2$}
\put(-1.3,0.8) {\tiny$\lam_1$}
\put(0,-2) {\tiny$\lam_3$}
\path(-3,1.5)(3,-1.5)
\path(-6,0)(0,3) \path(0,3)(4.1,.95) 
\qbezier[9](4.2,.9)(4.7,.65)(5.2,.4) \path(5.2,.4)(6,0) 
\path(6,0)(0,-3) \path(0,-3)(-6,0) 
\qbezier[25](-3,1.5)(-4,1)(-4.5,-.7) \put(-3.5,1.8){$l_1^-$}
\qbezier[90](-4.5,-.7)(-4.9,-2)(-5,-5)
\qbezier[220](-5,-5)(0,-4.7)(5.2,-5) \put(-4.6,-4.4){$\Gamma_1'$}
\qbezier[140](5.2,-5)(5,-2)(5.2,2)
\qbezier[95](5.2,2)(5,-.5)(3,-1.48) \put(3,-2.1){$l_1^+$}
\qbezier[35](0,0)(1.2,-.6)(4.2,-.5)
\qbezier[40](4.2,-.5)(4.4,-.5)(5.1,-.4) \put(5.3,-.9){$l$}
\thicklines
\qbezier[1000](5.2,2)(2,-1)(0,0) \put(5.1,2.3){$l_1'$}
\put(-.8,2){$\Gamma_1$}
\put(-5.3,-.2){$\Gamma_2$}
\put(2.8,-0.75){$\Gamma_2'$}
\put(-2.8,.7){\tiny$A_{e,3}$}
\put(1.1,-1.2){\tiny$A_{e,3}$}
\put(0.2,-4.5){\tiny$A_2$}
\put(3.8,1.5){\tiny$A_3$}
\put(1.2,1.7){\tiny$A_{e,2}$}
\put(-3,-1.4){\tiny$A_{e,2}$}
\put(-.55,.35){\tiny$A_{e,4}$} \put(-.18,-.17){$\bullet$}
\put(-1,4){II}
\put(5.7,-2){I}
\put(6.1,.3){III} \path(6,.5)(4.95,0.1)
\put(-5.7,-1.7){IV}
\end{picture}
\end{center}

\vskip 0.5truecm
\centerline{{\bf Fig.~4:}
A typical singularity $A_{e,4}^\pm$ of caustic surfaces}
\end{figure}

a) Two surfaces $\Gamma_1,\Gamma_2$ corresponding to a degeneracy of a
{\it basic} critical point. These surfaces adjoin each other on a curve $l_1$
which passes through the point $\overline{\lam}$. Here the degeneracy type
of the basic point of this curve equals $A_{e,3}\subset A_5$ (except the
point $\overline{\lam}$).
The point $\overline{\lam}$ splits the curve $l_1$ into two curves which will
be denoted by $l_1^+$ and $l_1^-$, resp.

b) A surface $\Gamma_1'$ corresponding to a degeneracy of an {\it additional}
critical point. The surface $\Gamma_1'$ adjoins the surfaces
$\Gamma_1$,$\Gamma_2$ on the curve $l_1^-$,
as well as intersects the surface $\Gamma_1$ in a curve $l$ which is tangent
to $l_1^+$ at the point $\overline{\lam}$.
A degeneracy of a basic point of type $A_{e,2}\subset A_3$ and
a degeneracy of an additional point of type $A_2$
correspond to values $\lam \in l$ (notation $A_{e,2}+A_2$).

c) A surface $\Gamma_2'$ corresponding to a degeneracy of type $A_2$ of an
{\it additional} critical point. The surface $\Gamma_2'$ adjoins the surfaces
$\Gamma_1$,$\Gamma_2$ on the curve $l_1^+$. Moreover the contiguity of
$\Gamma_1'$ and $\Gamma_2'$ to the surfaces $\Gamma_1,\Gamma_2$ takes place
on different sides. The surfaces $\Gamma_1'$ and $\Gamma_2'$ adjoin each other
on the curve $l_1'$ which is tangent to $l_1$ at the point $\overline{\lam}$ and
corresponds to a degeneracy of type $A_3$ of an additional critical point.

For $\nu=1$ and $\overline{\lam}=0$, an even versal deformation of a function in
a neighbourhood of a degenerate basic critical point $0$ of type $A_{e,4}^+$ has
the form $F(k,\lam)=k^8+\lam_1k^6+\lam_2k^4+\lam_3k^2$.
The branches $\Gamma_1\cup l_1\cup\Gamma_2$ and $\Gamma_1'$, $\Gamma_2'$
of the caustic surface are respectively defined by the equations
 $$
\begin{array}{lll}
\lam_3=0, && \\
8\lam_3+\lam_1(-4\lam_2+\lam_1^2)=(-\frac{8}{3}\lam_2+\lam_1^2)^{3/2} &
\mbox{for} &
\mbox{$\lam_2<0$ or $0\le\lam_2\le\frac{3}{8}\lam_1^2$, $\lam_1<0$}, \\
8\lam_3+\lam_1(-4\lam_2+\lam_1^2)=-(-\frac{8}{3}\lam_2+\lam_1^2)^{3/2} &
\mbox{for} &
\mbox{$0<\lam_2\le\frac{3}{8}\lam_1^2$, $\lam_1<0$.}
\end{array}
 $$

\begin{figure}
\unitlength = 5mm
\begin{center}
\begin{picture}(30,18.5)(0,-13.5) 
\path(0,-5)(0,5)
\path(0,5)(20,5)
\path(20,5)(20,-5)
\path(20,-5)(0,-5) 
\put(0,0){\vector(1,0){1.2}} \put(1.2,.2){\tiny$\varphi$}
\put(0,0){\vector(0,-1){1.2}} \put(-.7,-1.2){\tiny$\lam_3$} 
\thicklines 
\qbezier[70](5,0)(6.5,0)(8.5,-2) \put(5,0){\put(-.18,-.17){$\bullet$}}
\qbezier[130](8.5,-2)(10.6,-4.4)(14,-3.7)
\qbezier[130](14,-3.7)(18.2,-2.5)(18,3) \put(18,3){\put(-.18,-.17){$\bullet$}}
\qbezier[70](18,3)(17.5,0)(15,0) \put(15,0){\put(-.18,-.17){$\bullet$}}
\path(0,0)(20,0)

\put(4.5,.4){\tiny$A_{e,3}$} \put(10,.3){\tiny$A_{e,2}$} \put(2.5,-.5){\tiny$A_{e,2}$}
\put(14.5,.4){\tiny$A_{e,3}$} \put(10,-4.1){\tiny$A_2$} \put(17.7,3.4){\tiny$A_3$}
\put(18.2,1.3){\tiny$A_2$}

\put(-.5,-12.6){\tiny$A_{e,3}$} \put(14,-13.5){\tiny$A_{e,2}$}
\put(29.5,0){\tiny$A_2$} \put(28.5,-13.5){\tiny$A_2+A_{e,2}$} \put(13.7,-9.6){\tiny$A_2$}

\thinlines 
\put(2,-3){I}
\put(5,-3){\circle{1.5}
  \path(-.7,0)(.7,0)
  \put(0,.5){\vector(0,1){.2}} \path(0,-.6)(0,.6)
  \qbezier[70](-.5,.55)(0,-.5)(.5,.55)
 }
\put(4,3){II}
\put(8,2.5){\circle{1.5}
  \path(-.7,0)(.7,0)
  \put(0,.5){\vector(0,1){.2}} \path(0,-.6)(0,.6)
  \qbezier[70](-.5,.55)(-.3,-.5)(-.15,-.2)
  \qbezier[70](-.15,-.2)(0,.2)(.15,-.2)
  \qbezier[70](.5,.55)(.3,-.5)(.15,-.2)
 }
\put(10,-2){IV}
\put(12,-1){\circle{1.5}
  \path(-.7,0)(.7,0)
  \put(0,.5){\vector(0,1){.2}} \path(0,-.6)(0,.6)
  \qbezier[80](-.5,.55)(-.4,-.5)(-.25,0)
  \qbezier[50](-.25,0)(-.18,.3)(-.1,.1)
  \qbezier[50](-.1,.1)(0,-.1)(.1,.1)
  \qbezier[50](.25,0)(.18,.3)(.1,.1)
  \qbezier[80](.5,.55)(.4,-.5)(.25,0)
 }
\qbezier[35](5,0)(6.5,0)(8.5,-1)
\qbezier[65](8.5,-1)(10.6,-2.2)(14,-1.8) 
\qbezier[35](14,-1.8)(16,-1.25)(17,0)
\put(12.5,-2.95){\circle{1.5}
  \path(-.7,0)(.7,0)
  \put(0,.5){\vector(0,1){.2}} \path(0,-.6)(0,.6)
  \qbezier[50](-.5,.55)(-.4,0)(-.3,.2)
  \qbezier[50](-.3,.2)(-.2,.4)(-.12,.2)
  \qbezier[50](-.12,.2)(0,-.2)(.12,.2)
  \qbezier[50](.3,.2)(.2,.4)(.12,.2)
  \qbezier[50](.5,.55)(.4,0)(.3,.2)
 }
\put(16,1.5){III} \path(16.6,1.3)(17.25,0.65) 
\qbezier[30](17,0)(18,1.5)(18,3)
\qbezier[30](15,0)(16.5,0)(17.9,.8) 
\qbezier[12](17,0)(17.5,.5)(17.9,.8)
 %
\thicklines 
\qbezier[400](28,-14)(29.2,-5)(29,5)
\qbezier[600](27,5)(20,-13)(0,-12.95) \put(0,-13){\put(-.18,-.17){$\bullet$}}
\path(-2,-13)(30,-13) \put(15,-11.2){III}
\thinlines
\qbezier[100](20,-12.95)(26,-5)(28,5)
\qbezier[150](0,-12.95)(15,-13)(28.86,-5) 
\qbezier[60](20,-12.95)(25,-8)(28.86,-5)
\qbezier[7](20,-12.95)(19.4,-13.6)(19,-14)
\put(19,-11.2){\circle{1.5}
  \path(-.7,0)(.7,0)
  \put(0,.5){\vector(0,1){.2}} \path(0,-.6)(0,.6)
  \qbezier[80](-.5,.55)(-.42,-1)(-.3,0)
  \qbezier[50](-.3,0)(-.26,.3)(-.2,0)
  \qbezier[50](-.2,0)(-.14,-.3)(-.08,-.1)
  \qbezier[50](-.08,-.1)(0,.1)(.08,-.1)
  \qbezier[50](.2,0)(.14,-.3)(.08,-.1)
  \qbezier[50](.3,0)(.26,.3)(.2,0)
  \qbezier[80](.5,.55)(.42,-1)(.3,0)
 }
\put(22,-6.2){\circle{1.5}
  \path(-.7,0)(.7,0)
  \put(0,.5){\vector(0,1){.2}} \path(0,-.6)(0,.6)
  \qbezier[80](-.5,.55)(-.43,-.8)(-.3,-.33)
  \qbezier[50](-.3,-.33)(-.26,.1)(-.2,-.2)
  \qbezier[50](-.2,-.2)(-.14,-.4)(-.07,-.1)
  \qbezier[50](-.07,-.1)(0,.1)(.07,-.1)
  \qbezier[50](.2,-.2)(.14,-.4)(.07,-.1)
  \qbezier[50](.3,-.33)(.26,.1)(.2,-.2)
  \qbezier[80](.5,.55)(.43,-.8)(.3,-.33)
 }

\put(24.5,-11.2){\circle{1.5}
  \path(-.7,0)(.7,0)
  \put(0,.5){\vector(0,1){.2}} \path(0,-.6)(0,.6)
  \qbezier[50](-.5,.55)(-.4,0)(-.33,.2)
  \qbezier[50](-.33,.2)(-.29,.4)(-.22,0)
  \qbezier[50](-.22,0)(-.16,-.3)(-.1,-.1)
  \qbezier[50](-.1,-.1)(0,.1)(.1,-.1)
  \qbezier[50](.22,0)(.16,-.3)(.1,-.1)
  \qbezier[50](.33,.2)(.29,.4)(.22,0)
  \qbezier[50](.5,.55)(.4,0)(.33,.2)
 }

\path(25.6,-8.6)(24,-8.4)
\put(26.3,-8.8){\circle{1.5}
  \path(-.7,0)(.7,0)
  \put(0,.5){\vector(0,1){.2}} \path(0,-.6)(0,.6)
  \qbezier[80](-.5,.55)(-.42,-.4)(-.32,0)
  \qbezier[50](-.32,0)(-.26,.3)(-.2,0)
  \qbezier[50](-.2,0)(-.14,-.5)(-.07,-.1)
  \qbezier[50](-.07,-.1)(0,.1)(.07,-.1)
  \qbezier[50](.2,0)(.14,-.5)(.07,-.1)
  \qbezier[50](.32,0)(.26,.3)(.2,0)
  \qbezier[80](.5,.55)(.42,-.4)(.32,0)
 }

\put(27,-4.2){\circle{1.5}
  \path(-.7,0)(.7,0)
  \put(0,.5){\vector(0,1){.2}} \path(0,-.6)(0,.6)
  \qbezier[80](-.5,.55)(-.43,-.5)(-.32,-.25)
  \qbezier[50](-.32,-.25)(-.26,.1)(-.2,-.2)
  \qbezier[50](-.2,-.2)(-.14,-.8)(-.07,-.1)
  \qbezier[50](-.07,-.1)(0,.1)(.07,-.1)
  \qbezier[50](.2,-.2)(.14,-.8)(.07,-.1)
  \qbezier[50](.32,-.25)(.26,.1)(.2,-.2)
  \qbezier[80](.5,.55)(.43,-.5)(.32,-.25)
 }
\end{picture}
\end{center}

\vskip 0.5truecm
\centerline{{\bf Fig.~5:}
Bifurcations of functions at a singularity $A_{e,4}^+$
in typical families of even functions}
\end{figure}

Bifurcations of even functions in $\nu=1$ variable in a neighbourhood of a
caustic value corresponding to a singularity of type $A_{e,4}^+$ are shown in
Fig.~5 (compare Fig.~91~\cite{Arnold2}). This figure does not show the
whole caustic, but only its intersection with a cylinder
($\lam_1=\varepsilon\sin\varphi$, $\lam_2=\varepsilon\cos\varphi$,
$0\le\varphi\le2\pi$ for a small positive constant $\varepsilon$),
which intersects the plane $\Gamma_1\cup l_1\cup\Gamma_2$ transversally.
Let us describe the behaviour of critical points for parameter values $\lam$
close to $\overline{\lam}$:
\\In the domain $I$, there are no twin critical points.
Under a passage across the caustic $\Gamma_1$ to the domain $II$,
a twin point $(k_1,-k_1)$ separates from the basic critical point.
\\Under a passage across the caustic $\Gamma_1'$ from the domain $II$ to the
domain $III$, two twin critical points $(k_2,-k_2)$ and $(k_3,-k_3)$ arise.
Thus, for the domain $III$, there are 3 twin critical points in a neighbourhood
of the basic critical point.
\\Under a passage across the caustic $\Gamma_2'$ from the domain $III$ to the
domain $II$, the twin critical points $(k_1,-k_1)$ and $(k_2,-k_2)$ merge and
disappear.
\\Under a passage across the caustic $\Gamma_2$ from the domain $II$ to the
domain $IV$, the twin critical point $(k_2,-k_2)$ separates from the basic
critical point.
\\Under a passage across the caustic $\Gamma_1'$ from the domain $IV$ to the
domain $I$, the twin critical points $(k_2,-k_2)$ and $(k_3,-k_3)$ merge and
disappear.
\\ {\bf 9, 10.} ($X_{e,5}^{++}$, $X_{e,5}^{--}$, and $X_{e,5}^{+-}$)
As a ``model'' (more precisely, a one-parameter family of ``models'') of
a typical 3-parameter family of functions in $\nu=2$ variables, we consider
 $$
\pm x^4+ax^2y^2\pm y^4 + \lam_1 x^2 + \lam_2 xy + \lam_3 y^2.
 $$
Here $\lam_1,\lam_2,\lam_3$ are parameters of the family and $a$ is a parameter
of the model ($a\ne\pm2$ in the cases $++$ and $--$).
In a neighbourhood of a caustic value $\overline{\lam}$ corresponding
to a degeneracy of a {\it basic} critical point with singularity class
$X_{e,5}\subset X_9$, the set of caustic values is a surface which is
the union of two conic surfaces (more precisely, ``surfaces close to conic
ones'') $\Gamma$ and $\Gamma'$ with vertex at the point $\overline{\lam}$.
Moreover these surfaces can intersect each other on rays coming from
$\overline{\lam}$.
The surface $\Gamma$ (except the point $\overline{\lam}$)
corresponds to degeneracies of types $A_{e,2}$ and $A_{e,3}$ of a {\it basic}
critical point, while the surface $\Gamma'$ (except the point $\overline{\lam}$)
corresponds to degeneracies of types $A_2$ and $A_3$ of an {\it additional}
critical point. (Actually, for our model, as for any generic 1-parameter family
of ``models'' with a parameter $a$, degeneracies of types $A_4$ and $D_4^\pm$
as well as non-transversal intersections of different branches of the caustic
can happen in general for some exceptional values of $a$. However the number of
such exceptional values of $a$ is finite, thus the corresponding models are not
generic and, hence, one does not need to consider them.)
More precisely, the surfaces $\Gamma$ and $\Gamma'$ have the following form:

a) The surface $\Gamma$ is a circular cone having vertex at the point
$\overline{\lam}$ (in our model it is the cone $\lambda_2^2=4\lambda_1\lambda_3$).
Finitely many elements $l_1,l_2,\dots$ of the cone $\Gamma$ depending on $a$
are chosen:
\\In the case $++$, there are no chosen elements of the cone for $a>-2$, while
there are four chosen elements $l_1,l_2,l_3,l_4$ for $a<-2$ (for our model, they
are obtained by intersecting the cone with two planes
$\frac{\lam_3}{\lam_1}=\frac{-a+\sqrt{a^2-4}}{2}$ and
$\frac{\lam_3}{\lam_1}=\frac{-a-\sqrt{a^2-4}}{2}$).
\\ Similarly, in the case $--$, there are no chosen elements of the cone for
$a<2$, while there are four chosen elements $l_1,l_2,l_3,l_4$ for $a>2$ (for
our model, they are obtained by intersecting the cone with two planes
$\frac{\lam_1}{\lam_3}=\frac{a+\sqrt{a^2-4}}{2}$ and
$\frac{\lam_1}{\lam_3}=\frac{a-\sqrt{a^2-4}}{2}$).
\\ In the case $+-$ two elements $l_1$ and $l_2$ of the cone are chosen (for
our model, they are obtained by intersecting the cone by the plane
$\frac{\lam_1}{\lam_3}=\frac{-a+\sqrt{a^2+4}}{2}$).

All points of the cone $\Gamma$ except the point $\overline{\lam}$ and the
points lying on the chosen cone elements $l_1,l_2,\dots$ correspond to a
degeneracy of type $A_{e,2}$ of a {\it basic} critical point. All points of the
chosen cone elements apart from $\overline{\lam}$ correspond to a degeneracy of
type $A_{e,3}$ of a {\it basic} critical point.

Denote by $\Gamma^+$ and $\Gamma^-$ two half-cones of the cone $\Gamma$,
i.e.\ two connected components of $\Gamma\setminus\overline{\lam}$
(for our model, they are the parts of $\Gamma$ lying in the octants
$\lam_1>0$, $\lam_3>0$, and resp.\ $\lam_1<0$, $\lam_3<0$).
Denote by $\Omega^+$ and $\Omega^-$ the domains of $\RR^3$ bounded by
the surfaces $\Gamma^+$ and $\Gamma^-$, resp. Denote by $l_i^+$ and $l_i^-$
the rays of the cone element $l_i$ lying resp.\ in the half-cones $\Gamma^+$
and $\Gamma^-$.

b) The surface $\Gamma'$ is a conic surface which is formed by a number of
conic surfaces of type $A_2$. These surfaces adjoin the cone $\Gamma$ on its
chosen elements $l_1,l_2,\dots$ from above (more precisely, on their
rays $l_1^\pm,l_2^\pm,\dots$). Furthermore, these surfaces adjoin each
other on rays-cuspidal edges coming from $\overline{\lam}$ and having type
$A_3$. The shape of the surface $\Gamma'$ (in particular, the number
of cuspidal edges and the location of them with respect to the cone $\Gamma$)
depends in general on the signs $\pm\pm$ and the value of the module $a$
which appear in the normal form of an even germ of class $X_{e,5}$. In other
words, it depends on the degeneracy type of the basic critical point for
$\lam=\overline{\lam}$:
\\In the cases $++$ and $--$, if $a\ne0$ has a small enough absolute value
then the surface $\Gamma'$ consists of four surfaces $\Gamma_1'$, $\Gamma_2'$,
$\Gamma_3'$, and $\Gamma_4'$ of type $A_2$. These surfaces consecutively adjoin
each other on four rays $l_1'$, $l_2'$, $l_3'$, and $l_4'$ coming from
$\overline{\lam}$ and being cuspidal edges of type $A_3$.
More precisely, in the case $++$ with $0<a\ll1$, the surface $\Gamma'$ lies
inside the domain $\Omega^-$ (Fig.~6~a). But, in the case $++$ with $0<-a\ll1$,
each surface $\Gamma'_i$ intersects the cone $\Gamma$ transversally in a ray
corresponding a degeneracy of type $A_{e,2}+A_2$. In the latter case,
the cuspidal edges $l_2'$ and $l_4'$ lie inside the domain $\Omega^-$, while
the cuspidal edges $l_1'$ and $l_3'$ lie outside it (Fig.~6~b). In the case
$--$, the surface $\Gamma'$ has a similar form, see Remark~\ref {note:sym}.
\\In the case $+-$, if $a\ne0$ has a small enough absolute value, then
the surface $\Gamma'$ consists of four surfaces $\Gamma_1'$, $\Gamma_2'$,
$\Gamma_3'$, and $\Gamma_4'$ of type $A_2$. Here, two surfaces
$\Gamma_1'$ and $\Gamma_2'$ adjoin the cone $\Gamma$ on the rays $l_1^+$ and
resp.\ $l_2^+$ of the cone elements. Furthermore, these surfaces adjoin each
other on a cuspidal edge $l_1'$ (of type $A_3$) which is a ray coming from
$\overline{\lam}$. Other surfaces $\Gamma_3'$ and $\Gamma_4'$ adjoin the cone
$\Gamma$ on the rays $l_1^-$ and resp.\ $l_2^-$ of the cone elements.
Furthermore, these surfaces adjoin each other on a cuspidal edge $l_2'$
(of type $A_3$), which is a ray coming from $\overline{\lam}$.
Moreover, the surface $\Gamma_1'$ intersects each of the surfaces $\Gamma_4'$
and $\Gamma^-$ transversally in one ray coming from $\overline{\lam}$,
while the surface $\Gamma_2'$ intersects each of the surfaces $\Gamma_3'$ and
$\Gamma^-$ transversally in one ray coming from $\overline{\lam}$ (Fig.~6~c).

\begin{figure}
\unitlength = 5mm
\begin{center}
\begin{picture}(33,12)(-5,-7) 
\put(-10,-4) 
{
\qbezier[250](6.3,8.9)(10,9.8)(13.7,8.9)
\qbezier[150](6.3,8.9)(3.7,8)(6.3,7.1)
\qbezier[150](13.7,8.9)(16.3,8)(13.7,7.1)
\qbezier[250](6.3,7.1)(10,6.2)(13.7,7.1)
}
\path(-4.6,3.5)(4.6,-3.5) \path(-4.6,-3.5)(4.6,3.5) \put(-.18,-.17){$\bullet$}
\put(-10,-12) 
{
\qbezier[50](6.3,8.9)(10,9.8)(13.7,8.9)
\qbezier[150](6.3,8.9)(3.7,8)(6.3,7.1)
\qbezier[150](13.7,8.9)(16.3,8)(13.7,7.1)
\qbezier[250](6.3,7.1)(10,6.2)(13.7,7.1)
}
\qbezier[8](0,0)(.3,-.4)(.6,-.8) \put(.6,-.8){\vector(3,-4){.2}} \put(.8,-1.5) {\tiny$\lam_1$}
\put(0,0) {\vector(4,-1){.8}} \put(.8,-.3) {\tiny$\lam_2$}
\put(0,0) {\vector(-1,-2){.8}} \put(-1.3,-1.9) {\tiny$\lam_3$}

\qbezier[200](1.575,4.9)(0,4.1)(-1.3875,4.3375)
\qbezier[200](-1.3875,4.3375)(-.1,4)(-1.575,3.1) 
\qbezier[200](-1.575,3.1)(0,3.9)(1.3875,3.6625)
\qbezier[200](1.3875,3.6625)(.1,4)(1.575,4.9)

\qbezier[28](0,0)(0.5943396,1.8490566)(1.1886792,3.6981132)
\thicklines
 \path(1.1886792,3.6981132)(1.575,4.9)
\thinlines
\qbezier[27](0,0)(-0.54508925,1.70401785)(-1.0901785,3.4080357)
\thicklines
 \path(-1.0901785,3.4080357)(-1.3875,4.3375)
\thinlines
\qbezier[25](0,0)(-0.6909198,1.35990565)(-1.3818396,2.7198113)
\thicklines
 \path(-1.3818396,2.7198113)(-1.575,3.1)
\thinlines
\qbezier[25](0,0)(0.50296875,1.32765625)(1.0059375,2.6553125)
\thicklines
 \path(1.0059375,2.6553125)(1.3875,3.6625)

\thinlines
\qbezier[25](0,0)(-0.33804875,1.7121951)(-0.6760975,3.4243902)
 \path(-0.6760975,3.4243902)(-.77,3.9)
\qbezier[27](0,0)(0.3596899,1.84341085)(0.7193798,3.6868217)
 \path(0.7193798,3.6868217)(.8,4.1)

\put(-3.6,3.9){\tiny$\Omega^-$} \put(3,1.7){\tiny$\Gamma^-$}
\put(-3,-4.3){\tiny$\Omega^+$} \put(3,-2.1){\tiny$\Gamma^+$}
\put(1.7,4.75){\tiny$l'_1$} \put(.22,3.98){\tiny$\Gamma'_1$}
\put(1.4,3.2){\tiny$l'_2$} \put(-.15,3.15){\tiny$\Gamma'_2$}
\put(-2,3){\tiny$l'_3$} \put(-.73,3.82){\tiny$\Gamma'_3$}
\put(-1.8,4.2){\tiny$l'_4$} \put(-.3,4.6){\tiny$\Gamma'_4$}

\put(-1.3,-.05){\tiny$X_{e,5}$}
\put(-4.8,2.5){\tiny$A_{e,2}$}
\put(-4.8,-2.6){\tiny$A_{e,2}$}
\put(-.9,2.9){\tiny$A_2$}
\put(1.4,4){\tiny$A_3$}

\put(0,-7){\large a}
\put(11,0){
\put(-10,-4) 
{
\qbezier[250](6.3,8.9)(10,9.8)(13.7,8.9)
\qbezier[150](6.3,8.9)(3.7,8)(6.3,7.1)
\qbezier[150](13.7,8.9)(16.3,8)(13.7,7.1)
\qbezier[250](6.3,7.1)(10,6.2)(13.7,7.1)
}
\path(-4.6,3.5)(4.6,-3.5) \put(-.18,-.17){$\bullet$}
\path(-4.6,-3.5)(4.6,3.5)
\put(-10,-12) 
{
\qbezier[50](6.3,8.9)(10,9.8)(13.7,8.9)
\qbezier[150](6.3,8.9)(3.7,8)(6.3,7.1)
\qbezier[150](13.7,8.9)(16.3,8)(13.7,7.1)
\qbezier[250](6.3,7.1)(10,6.2)(13.7,7.1)
}
\qbezier[8](0,0)(.3,-.4)(.6,-.8) \put(.6,-.8){\vector(3,-4){.2}} \put(.8,-1.5) {\tiny$\lam_1$}
\put(0,0) {\vector(4,-1){.8}} \put(.8,-.3) {\tiny$\lam_2$}
\put(0,0) {\vector(-1,-2){.8}} \put(-1.3,-1.9) {\tiny$\lam_3$}

\qbezier[200](2.625,5.5)(0,4.1)(-2.3125,4.5625)
\qbezier[200](-2.3125,4.5625)(-.1,4)(-2.625,2.5) 
\qbezier[200](-2.625,2.5)(0,3.9)(2.3125,3.4375)
\qbezier[200](2.3125,3.4375)(.1,4)(2.625,5.5)

\qbezier[30](0,0)(1.22269735,2.5618421)(2.4453947,5.1236842)
\thicklines
 \path(2.4453947,5.1236842)(2.625,5.5)
\thinlines
\qbezier[27](0,0)(0.775,1.84642855)(1.55,3.6928571)
 \path(1.55,3.6928571)(2.17,5.17)
\qbezier[30](0,0)(1.0072254,2.55491325)(2.0144508,5.1098265)
 \path(2.0144508,5.1098265)(2.05,5.2)
\qbezier[30](0,0)(-0.8123022,1.6026503)(-1.6246044,3.2053006)
\thicklines
 \path(-1.6246044,3.2053006)(-2.3125,4.5625)
 \path(0,0)(-2.625,2.5)
\thinlines
 \path(0,0)(-2,2.76)
\qbezier[27](0,0)(-1.06184875,1.3764706)(-2.1236975,2.7529412)
 \path(-2.1236975,2.7529412)(-2.16,2.8)
\qbezier[25](0,0)(0.92138671875,1.36962890625)(1.8427734375,2.7392578125)
\thicklines
 \path(1.8427734375,2.7392578125)(2.3125,3.4375)

\thinlines
\qbezier[25](0,0)(-0.515,1.5892578125)(-1.03,3.178515625)
 \path(-1.03,3.178515625)(-1.28,3.95)
\qbezier[25](0,0)(0.54562045,1.7837591)(1.0912409,3.5675182)
 \path(1.0912409,3.5675182)(1.3,4.25)

\put(-3.6,3.9){\tiny$\Omega^-$} \put(3,1.7){\tiny$\Gamma^-$}
\put(-3,-4.3){\tiny$\Omega^+$} \put(3,-2.1){\tiny$\Gamma^+$}
\put(2.75,5.45){\tiny$l'_1$} \put(.6,3.98){\tiny$\Gamma'_1$}
\put(2.4,3.15){\tiny$l'_2$} \put(-.2,2.9){\tiny$\Gamma'_2$}
\put(-2.75,1.9){\tiny$l'_3$} \put(-1.2,3.82){\tiny$\Gamma'_3$}
\put(-2.7,4.3){\tiny$l'_4$} \put(-.3,4.8){\tiny$\Gamma'_4$}

\put(-1.3,-.05){\tiny$X_{e,5}$}
\put(-4.8,2.5){\tiny$A_{e,2}$}
\put(-4.8,-2.6){\tiny$A_{e,2}$}
\put(.8,3){\tiny$A_2$}
\put(-2.5,3.5){\tiny$A_3$}

\put(0,-7){\large b}
}
\put(22.5,0){
\put(-10,-4) 
{
\qbezier[250](6.3,8.9)(10,9.8)(13.7,8.9)
\qbezier[150](6.3,8.9)(3.7,8)(6.3,7.1)
\qbezier[150](13.7,8.9)(16.3,8)(13.7,7.1)
\qbezier[250](6.3,7.1)(10,6.2)(13.7,7.1)
}
\qbezier[30](0,0)(-1.63,1.234)(-3.26,2.468) \path(-3.26,2.468)(-4.6,3.5) 
\path(0,0)(4.6,-3.5)
\qbezier[30](0,0)(-2.3,-1.75)(-4.6,-3.5)
\path(0,0)(4.6,3.5)
\put(-.18,-.17){$\bullet$}
\put(-10,-12) 
{
\qbezier[50](6.3,8.9)(10,9.8)(13.7,8.9)
\qbezier[50](6.3,8.9)(3.7,8)(6.3,7.1)
\qbezier[150](13.7,8.9)(16.3,8)(13.7,7.1)
}
\qbezier[10](-3.7,-4.9)(-3.12,-5.05)(-2.48,-5.15)
\qbezier[200](-2.48,-5.15)(-.4,-5.495)(1.75,-5.25)
\qbezier[15](1.75,-5.25)(3,-5.09)(3.7,-4.9) 

\qbezier[8](0,0)(.3,-.4)(.6,-.8) \put(.6,-.8){\vector(3,-4){.2}} \put(.8,-1.5) {\tiny$\lam_1$}
\put(0,0) {\vector(4,-1){.8}} \put(.8,-.3) {\tiny$\lam_2$}
\qbezier[12](0,0)(-.3,-.8)(-.6,-1.6)
\put(-.54,-1.5) {\vector(-1,-3){.1}} \put(-.45,-1.9) {\tiny$\lam_3$}

\path(0,0)(3.15,2.95) 
\qbezier[25](0,0)(-1.575,-1.475)(-3.15,-2.95) 
\qbezier[25](0,0)(-1.285,1.435)(-2.57,2.87) \path(-2.57,2.87)(-4.2,4.7) 
\path(0,0)(4.2,-4.7)

\qbezier[200](-1.575,3.1)(-2.1,2.8)(-4.3,2.1)
\qbezier[200](-1.575,3.1)(-2.1,2.8)(-3.3,1.8) 
\qbezier[25](0,0)(-0.6909198,1.35990565)(-1.3818396,2.7198113) 
\thicklines \path(-1.3818396,2.7198113)(-1.575,3.1) \thinlines

\qbezier[6](-2.625,-5.5)(-2.8,-5.58)(-3.1,-5.65) 
\qbezier[50](-3.1,-5.65)(-3.5,-5.75)(-4.3,-5.9) 
\qbezier[200](-2.625,-5.5)(-2.9,-5.65)(-3.3,-6.2) 
\qbezier[15](0,0)(-.54,-1.11)(-1.08,-2.22) 
\thicklines \path(-1.08,-2.22)(-2.625,-5.5) \thinlines

\qbezier[200](3.15,2.95)(1,2.5)(.5,2) 
\qbezier[200](-4.2,4.7)(-5,4.4)(-5.7,3.7)

\qbezier[19](-3.15,-2.95)(-4,-3.15)(-4.86,-3.45) 
\qbezier[50](-4.86,-3.45)(-5.5,-3.65)(-5.9,-3.8) 
\qbezier[200](4.2,-4.7)(2.7,-5.5)(2.7,-6)

\qbezier[300](.5,2)(-1.9,-2.1)(-3.3,-6.2)
\qbezier[300](-3.3,1.8)(-1.7,-2.6)(2.7,-6)
\qbezier[300](-4.3,2.1)(-5.6,-.85)(-5.9,-3.8) 
\qbezier[300](-5.7,3.7)(-5.5,-1.1)(-4.3,-5.9)

 \path(0,0)(-2,2.76)  
\qbezier[27](0,0)(-1.06184875,1.3764706)(-2.1236975,2.7529412)
 \path(-2.1236975,2.7529412)(-2.16,2.8)

\path(0,0)(-1.45,-1.75) 
\qbezier[18](0,0)(-1.17,-.155)(-2.34,-.31) 
\path(-2.34,-.31)(-5.3,-.73)

\put(.5,4.3){\tiny$\Omega^-$} \put(3,1.7){\tiny$\Gamma^-$}
\put(-.5,-4.7){\tiny$\Omega^+$} \put(3,-2.1){\tiny$\Gamma^+$}
\put(-1.6,3.4){\tiny$l'_1$} \put(-4.3,1.4){\tiny$\Gamma'_1$}
\put(-2.6,-6){\tiny$l'_2$} \put(-3,1.4){\tiny$\Gamma'_2$}
\put(-2.2,-1.8){\tiny$l'_{23}$} \put(.5,1.5){\tiny$\Gamma'_3$}
\put(-6.05,-.75){\tiny$l'_{14}$} \put(-5.5,3.2){\tiny$\Gamma'_4$}
\put(-3.5,-3.5){\tiny$l_1^+$}
\put(3.1,3.3){\tiny$l_1^-$}
\put(4.25,-5.15){\tiny$l_2^+$}
\put(-4.4,5.1){\tiny$l_2^-$}

\put(.1,.3){\tiny$X_{e,5}$}
\put(2.9,-3.6){\tiny$A_{e,3}$}
\put(5,-4.5){\tiny$A_{e,2}$}
\put(5,4.5){\tiny$A_{e,2}$}
\put(-2.1,-4.5){\tiny$A_3$}
\put(2.85,-6){\tiny$A_2$}
\put(-5.05,-5.85){\tiny$A_2$}
\put(-5.7,-3.3){\tiny$A_2$}

\put(0,-7){\large c}
}
\end{picture}
\end{center}

\vskip 0.5truecm
\centerline{\begin{tabular}{rl} \bf Fig.~6:
& A typical singularity $X_{e,5}^{\pm\pm}$ of caustic surfaces,
with module $a\in\RR$
\\ \phantom{$I^{I^{I^I}}$}\!\!\!\!\!\!\!\!\!\! &
{\bf a: } $X_{e,5}^{++}$, \ \ $0<a\ll1$; \qquad
{\bf b: } $X_{e,5}^{++}$, \ \ $0<-a\ll1$; \qquad
{\bf c: } $X_{e,5}^{+-}$, \ \ $0<a\ll1$
\end{tabular}
}
\end{figure}

We obtained the above description of the form of the caustic $\Gamma\cup\Gamma'$
for $0<|a|\ll1$ from an explicit presentation of the surface $\Gamma'$ for $a=0$
in our model (taking into account the following Remarks~\ref {note:sym}
and~\ref {note:edges}). Namely, in the case $++$ with $a=0$,
the surface $\Gamma'$ is defined by the equation
 $$
(\varrho-1)^3=\frac{27}{16}(2-\sigma)\varrho \quad
\mbox{for $0<\varrho\le1$, $\lam_1<0$}.
 $$
Furthermore, in the case $+-$ with $a=0$,
the surface $\Gamma'$ is the closure of the surface defined by the equation
 $$
(\varrho-1)^3=\frac{27}{16}(2+\sigma)\varrho \quad
\begin{array}{l} \mbox{\phantom{and} for $\varrho<0$, $\lam_1<0$,} \\
                 \mbox{and for $\varrho\ge4$, $\lam_1(\lam_1^2-\lam_3^2)<0$}.
\end{array}
 $$
Here one denotes $\varrho=\frac{\lam_2^2}{\lam_1\lam_3}$,
$\sigma=\frac{\lam_1^2}{\lam_3^2}+\frac{\lam_3^2}{\lam_1^2}$.
Actually, in the case $+-$, the conic surface $\Gamma'$ is
obtained from the indicated surface by adding four rays, where two of these rays
lie on the plane $\lam_1=0$, and the other two rays lie on the plane $\lam_3=0$.

\begin{figure}
\unitlength = 6.5mm
\begin{center}
\begin{picture}(20,32)(-10,-16) 
\path(-9,-16)(-9,16) \path(-9,16)(9,16) \path(9,16)(9,-16) \path(9,-16)(-9,-16)


\put(0,8) 
{
\qbezier[300](-4.44,2.7)(0,5.4)(4.44,2.7)
\qbezier[200](-4.44,2.7)(-7.56,0)(-4.44,-2.7)
\qbezier[200](4.44,2.7)(7.56,0)(4.44,-2.7)
\qbezier[300](-4.44,-2.7)(0,-5.4)(4.44,-2.7)
}
\put(0,-8) 
{
\qbezier[300](-4.44,2.7)(0,5.4)(4.44,2.7)
\qbezier[200](-4.44,2.7)(-7.56,0)(-4.44,-2.7)
\qbezier[200](4.44,2.7)(7.56,0)(4.44,-2.7)
\qbezier[300](-4.44,-2.7)(0,-5.4)(4.44,-2.7)
}
\put(0,-12.06) {\vector(1,0){2}} \put(0,-12.06) {\put(1.8,-.5){$\lam_2$}}
\put(0,-12.06) {\vector(0,1){2}} \put(0,-12.06) {\put(-.5,1.8){$\varphi$}}
\qbezier[200](-3.5,8)(-.5,8.5)(0,11.5) \put(-3.5,8){\put(-.145,-.14){$\bullet$}}
\qbezier[200](3.5,8)(.5,8.5)(0,11.5) \put(0,11.5){\put(-.145,-.14){$\bullet$}}
\qbezier[200](-3.5,8)(-.5,7.5)(0,4.5) \put(0,4.5){\put(-.145,-.14){$\bullet$}}
\qbezier[200](3.5,8)(.5,7.5)(0,4.5) \put(3.5,8){\put(-.145,-.14){$\bullet$}}
\qbezier[80](-3.5,8)(0,8)(3.5,8)
\qbezier[80](0,11.5)(0,8)(0,4.5)
\put(0,8){
\qbezier[20](.25,0)(.15,.15)(0,.25) 
\qbezier[40](0,.25)(-.5,.5)(-.25,0)
\qbezier[20](-.25,0)(-.15,-.15)(0,-.25)
\qbezier[40](0,-.25)(.5,-.5)(.25,0)

\qbezier[40](.25,0)(.5,.5)(0,.25)
\qbezier[20](0,.25)(-.15,.15)(-.25,0)
\qbezier[40](-.25,0)(-.5,-.5)(0,-.25)
\qbezier[20](0,-.25)(.15,-.15)(.25,0)

\put(0,0){\put(-.066,-.144){$\cdot$}}
\put(-.21,.21){\put(-.066,-.144){$\cdot$}}
\put(.21,.21){\put(-.066,-.144){$\cdot$}}
\put(.21,-.21){\put(-.066,-.144){$\cdot$}}
\put(-.21,-.21){\put(-.066,-.144){$\cdot$}}
}

\put(-.75,8.735){\put(-.066,-.144){$\cdot$}
\put(.05,.215){\circle{.6}} 
\put(-.05,-.215){\circle{.6}}
\put(.02,.265){
\path(.065,-.06)(-.035,.04) 
\path(.065,.06)(-.035,-.04)
\qbezier[10](-.05,.05)(-.09,.08)(-.13,.05) 
\qbezier[10](-.13,.05)(-.18,0)(-.13,-.05)
\qbezier[10](-.05,-.05)(-.09,-.08)(-.13,-.05)
\qbezier[10](.06,.06)(.12,.12)(.18,.06) 
\qbezier[10](.18,.06)(.24,0)(.18,-.06)
\qbezier[10](.06,-.06)(.12,-.12)(.18,-.06)
\put(-.09,0){\put(-.066,-.144){$\cdot$}}
\put(.115,0){\put(-.066,-.144){$\cdot$}}
}
\put(-.02,-.265){
\path(-.045,.06)(.055,-.04) 
\path(-.045,-.06)(.055,.04)
\qbezier[10](.05,-.05)(.09,-.08)(.13,-.05) 
\qbezier[10](.13,-.05)(.18,0)(.13,.05)
\qbezier[10](.05,.05)(.09,.08)(.13,.05)
\qbezier[10](-.06,-.06)(-.12,-.12)(-.18,-.06) 
\qbezier[10](-.18,-.06)(-.24,0)(-.18,.06)
\qbezier[10](-.06,.06)(-.12,.12)(-.18,.06)
\put(.09,0){\put(-.066,-.144){$\cdot$}}
\put(-.115,0){\put(-.066,-.144){$\cdot$}}
}
}

 \put(0,9.775){\put(-.066,-.144){$\cdot$}
\put(0,.225){\circle{.6}} 
\put(0,-.225){\circle{.6}}
\put(0,.245){
\path(-.045,-.06)(.07,.06) 
\path(-.045,.06)(.07,-.06)
\qbezier[10](.06,-.06)(.12,-.12)(.18,-.06) 
\qbezier[10](.18,-.06)(.24,0)(.18,.06)
\qbezier[10](.06,.06)(.12,.12)(.18,.06)
\qbezier[10](-.06,.06)(-.12,.12)(-.18,.06) 
\qbezier[10](-.18,.06)(-.24,0)(-.18,-.06)
\qbezier[10](-.06,-.06)(-.12,-.12)(-.18,-.06)
\put(.12,0){\put(-.066,-.144){$\cdot$}}
\put(-.12,0){\put(-.066,-.144){$\cdot$}}
}
\put(0,-.245){
\path(-.04,-.05)(.06,.05) 
\path(-.04,.05)(.06,-.05)
\qbezier[10](.06,-.06)(.12,-.12)(.18,-.06) 
\qbezier[10](.18,-.06)(.24,0)(.18,.06)
\qbezier[10](.06,.06)(.12,.12)(.18,.06)
\qbezier[10](-.06,.06)(-.12,.12)(-.18,.06) 
\qbezier[10](-.18,.06)(-.24,0)(-.18,-.06)
\qbezier[10](-.06,-.06)(-.12,-.12)(-.18,-.06)
\put(.12,0){\put(-.066,-.144){$\cdot$}}
\put(-.12,0){\put(-.066,-.144){$\cdot$}}
}
}

\put(.7,8.95){\circle{.6}} 
\put(.8,8.52){\circle{.6}}
\put(.73,9){
\path(-.045,-.06)(.055,.04) 
\path(-.045,.06)(.055,-.04)
\qbezier[10](.05,.05)(.09,.08)(.13,.05) 
\qbezier[10](.13,.05)(.18,0)(.13,-.05)
\qbezier[10](.05,-.05)(.09,-.08)(.13,-.05)
\qbezier[10](-.06,.06)(-.12,.12)(-.18,.06) 
\qbezier[10](-.18,.06)(-.24,0)(-.18,-.06)
\qbezier[10](-.06,-.06)(-.12,-.12)(-.18,-.06)
\put(.09,0){\put(-.066,-.144){$\cdot$}}
\put(-.115,0){\put(-.066,-.144){$\cdot$}}
}
 \put(.75,8.735){\put(-.066,-.144){$\cdot$}}
\put(.77,8.47){
\path(.07,.06)(-.03,-.04) 
\path(.07,-.06)(-.03,.04)
\qbezier[10](-.05,-.05)(-.09,-.08)(-.13,-.05) 
\qbezier[10](-.13,-.05)(-.18,0)(-.13,.05)
\qbezier[10](-.05,.05)(-.09,.08)(-.13,.05)
\qbezier[10](.06,-.06)(.12,-.12)(.18,-.06) 
\qbezier[10](.18,-.06)(.24,0)(.18,.06)
\qbezier[10](.06,.06)(.12,.12)(.18,.06)
\put(-.09,0){\put(-.066,-.144){$\cdot$}}
\put(.115,0){\put(-.066,-.144){$\cdot$}}
}

\put(1.7,8){
\qbezier[10](.2,.05)(.15,.15)(.05,.2) 
\qbezier[40](.05,.2)(-.45,.45)(-.2,-.05)
\qbezier[10](-.2,-.05)(-.15,-.15)(-.05,-.2)
\qbezier[40](-.05,-.2)(.45,-.45)(.2,.05)

\qbezier[30](.2,.05)(.35,.35)(.05,.2)
\qbezier[30](.05,.2)(-.12,.12)(-.2,-.05)
\qbezier[30](-.2,-.05)(-.35,-.35)(-.05,-.2)
\qbezier[30](-.05,-.2)(.12,-.12)(.2,.05)

\put(0,0){\put(-.066,-.144){$\cdot$}}
\put(.18,.18){\put(-.066,-.144){$\cdot$}}
\put(-.18,.18){\put(-.066,-.144){$\cdot$}}
\put(-.18,-.18){\put(-.066,-.144){$\cdot$}}
\put(.18,-.18){\put(-.066,-.144){$\cdot$}}
}

 \put(.75,7.265){\put(-.066,-.144){$\cdot$}
\put(-.215,.05){\circle{.6}} 
\put(.215,-.05){\circle{.6}}
 \put(-.265,.02){
\path(-.05,.06)(.05,-.04) 
\path(.07,.06)(-.03,-.04)
\qbezier[10](.05,-.05)(.08,-.09)(.05,-.13) 
\qbezier[10](.05,-.13)(0,-.18)(-.05,-.13)
\qbezier[10](-.05,-.05)(-.08,-.09)(-.05,-.13)
\qbezier[10](.06,.06)(.12,.12)(.06,.18) 
\qbezier[10](.06,.18)(0,.24)(-.06,.18)
\qbezier[10](-.06,.06)(-.12,.12)(-.06,.18)
\put(0,-.09){\put(-.066,-.144){$\cdot$}}
\put(0,.115){\put(-.066,-.144){$\cdot$}}
}
\put(.265,-.02){
\path(.07,-.06)(-.03,.04) 
\path(-.05,-.06)(.05,.04)
\qbezier[10](-.05,.05)(-.08,.09)(-.05,.13) 
\qbezier[10](-.05,.13)(0,.18)(.05,.13)
\qbezier[10](.05,.05)(.08,.09)(.05,.13)
\qbezier[10](-.06,-.06)(-.12,-.12)(-.06,-.18) 
\qbezier[10](-.06,-.18)(0,-.24)(.06,-.18)
\qbezier[10](.06,-.06)(.12,-.12)(.06,-.18)
\put(0,.09){\put(-.066,-.144){$\cdot$}}
\put(0,-.115){\put(-.066,-.144){$\cdot$}}
}
}

 \put(0,6.4){\put(-.066,-.144){$\cdot$}
\put(.225,0){\circle{.6}} 
\put(-.225,0){\circle{.6}}
\put(.245,0){
\path(-.045,-.06)(.07,.06) 
\path(-.045,.06)(.07,-.06)
\qbezier[10](-.06,.06)(-.12,.12)(-.06,.18) 
\qbezier[10](-.06,.18)(0,.24)(.06,.18)
\qbezier[10](.06,.06)(.12,.12)(.06,.18)
\qbezier[10](.06,-.06)(.12,-.12)(.06,-.18) 
\qbezier[10](.06,-.18)(0,-.24)(-.06,-.18)
\qbezier[10](-.06,-.06)(-.12,-.12)(-.06,-.18)
\put(0,.12){\put(-.066,-.144){$\cdot$}}
\put(0,-.12){\put(-.066,-.144){$\cdot$}}
}
\put(-.245,0){
\path(-.04,-.05)(.06,.05) 
\path(-.04,.05)(.06,-.05)
\qbezier[10](-.06,.06)(-.12,.12)(-.06,.18) 
\qbezier[10](-.06,.18)(0,.24)(.06,.18)
\qbezier[10](.06,.06)(.12,.12)(.06,.18)
\qbezier[10](.06,-.06)(.12,-.12)(.06,-.18) 
\qbezier[10](.06,-.18)(0,-.24)(-.06,-.18)
\qbezier[10](-.06,-.06)(-.12,-.12)(-.06,-.18)
\put(0,.12){\put(-.066,-.144){$\cdot$}}
\put(0,-.12){\put(-.066,-.144){$\cdot$}}
}
}

 \put(-.75,7.265){\put(-.066,-.144){$\cdot$}
\put(-.215,-.05){\circle{.6}} 
\put(.215,.05){\circle{.6}}
 \put(-.265,-.02){
\path(-.05,-.06)(.05,.04) 
\path(.07,-.06)(-.03,.04)
\qbezier[10](.05,.05)(.08,.09)(.05,.13) 
\qbezier[10](.05,.13)(0,.18)(-.05,.13)
\qbezier[10](-.05,.05)(-.08,.09)(-.05,.13)
\qbezier[10](.06,-.06)(.12,-.12)(.06,-.18) 
\qbezier[10](.06,-.18)(0,-.24)(-.06,-.18)
\qbezier[10](-.06,-.06)(-.12,-.12)(-.06,-.18)
\put(0,.09){\put(-.066,-.144){$\cdot$}}
\put(0,-.115){\put(-.066,-.144){$\cdot$}}
}
\put(.265,.02){
\path(.07,.06)(-.03,-.04) 
\path(-.05,.06)(.05,-.04)
\qbezier[10](-.05,-.05)(-.08,-.09)(-.05,-.13) 
\qbezier[10](-.05,-.13)(0,-.18)(.05,-.13)
\qbezier[10](.05,-.05)(.08,-.09)(.05,-.13)
\qbezier[10](-.06,.06)(-.12,.12)(-.06,.18) 
\qbezier[10](-.06,.18)(0,.24)(.06,.18)
\qbezier[10](.06,.06)(.12,.12)(.06,.18)
\put(0,-.09){\put(-.066,-.144){$\cdot$}}
\put(0,.115){\put(-.066,-.144){$\cdot$}}
}
}

\put(-1.7,8){
\qbezier[10](-.2,.05)(-.15,.15)(-.05,.2) 
\qbezier[40](-.05,.2)(.45,.45)(.2,-.05)
\qbezier[10](.2,-.05)(.15,-.15)(.05,-.2)
\qbezier[40](.05,-.2)(-.45,-.45)(-.2,.05)

\qbezier[30](-.2,.05)(-.35,.35)(-.05,.2)
\qbezier[30](-.05,.2)(.12,.12)(.2,-.05)
\qbezier[30](.2,-.05)(.35,-.35)(.05,-.2)
\qbezier[30](.05,-.2)(-.12,-.12)(-.2,.05)

\put(0,0){\put(-.066,-.144){$\cdot$}}
\put(-.18,.18){\put(-.066,-.144){$\cdot$}}
\put(.18,.18){\put(-.066,-.144){$\cdot$}}
\put(.18,-.18){\put(-.066,-.144){$\cdot$}}
\put(-.18,-.18){\put(-.066,-.144){$\cdot$}}
}


\put(-2.5,10.5){\put(-.066,-.144){$\cdot$} 
\put(.05,.215){\circle{.6}}
\put(-.05,-.215){\circle{.6}}
\put(.06,.265){\put(-.066,-.144){$\cdot$}}
\put(-.06,-.265){\put(-.066,-.144){$\cdot$}}
}

\path(0,11.85)(.4,11.6)
\put(.7,11.3){\put(-.066,-.144){$\cdot$} 
\put(0,.225){\circle{.6}}
\put(0,-.225){\circle{.6}}
\put(0,.245){\put(-.066,-.144){$\cdot$}}
\put(0,-.245){\put(-.066,-.144){$\cdot$}}
}

\put(2.5,10.5){\put(-.066,-.144){$\cdot$} 
\put(-.05,.215){\circle{.6}}
\put(.05,-.215){\circle{.6}}
\put(-.06,.265){\put(-.066,-.144){$\cdot$}}
\put(.06,-.265){\put(-.066,-.144){$\cdot$}}
}

\put(4.8,8){
\put(-.15,.15){\circle{.6}} 
\put(.15,-.15){\circle{.6}}
\put(0,0){\put(-.066,-.144){$\cdot$}}
\put(-.2,.2){\put(-.066,-.144){$\cdot$}}
\put(.2,-.2){\put(-.066,-.144){$\cdot$}}
}

\put(2.5,5.5){\put(-.066,-.144){$\cdot$} 
\put(-.215,.05){\circle{.6}}
\put(.215,-.05){\circle{.6}}
\put(-.265,.06){\put(-.066,-.144){$\cdot$}}
\put(.265,-.06){\put(-.066,-.144){$\cdot$}}
}

\path(0,4.15)(-.38,4.38)
\put(-.85,4.55){\put(-.066,-.144){$\cdot$} 
\put(.225,0){\circle{.6}}
\put(-.225,0){\circle{.6}}
\put(.245,0){\put(-.066,-.144){$\cdot$}}
\put(-.245,0){\put(-.066,-.144){$\cdot$}}
}

\put(-2.5,5.5){\put(-.066,-.144){$\cdot$} 
\put(-.215,-.05){\circle{.6}}
\put(.215,.05){\circle{.6}}
\put(-.265,-.06){\put(-.066,-.144){$\cdot$}}
\put(.265,.06){\put(-.066,-.144){$\cdot$}}
}

\put(-4.8,8){
\put(.15,.15){\circle{.6}} 
\put(-.15,-.15){\circle{.6}}
\put(0,0){\put(-.066,-.144){$\cdot$}}
\put(.2,.2){\put(-.066,-.144){$\cdot$}}
\put(-.2,-.2){\put(-.066,-.144){$\cdot$}}
}


\put(-7.5,15){ 
\qbezier[40](0,0)(-.08,.35)(.1,.5) 
\qbezier[40](.1,.5)(.6,.6)(.5,.1)
\qbezier[40](0,0)(.35,-.08)(.5,.1)
\qbezier[40](0,0)(.08,-.35)(-.1,-.5) 
\qbezier[40](-.1,-.5)(-.6,-.6)(-.5,-.1)
\qbezier[40](0,0)(-.35,.08)(-.5,-.1)
\put(.25,.25){\put(-.066,-.144){$\cdot$}}
\put(-.25,-.25){\put(-.066,-.144){$\cdot$}}
}

\put(-3.75,15){ 
\qbezier[40](0,0)(-.3,.25)(-.12,.52) 
\qbezier[40](-.12,.52)(.23,.76)(.38,.35)
\qbezier[40](0,0)(.37,.05)(.38,.35)
\qbezier[40](0,0)(.3,-.25)(.12,-.52) 
\qbezier[40](.12,-.52)(-.23,-.76)(-.38,-.35)
\qbezier[40](0,0)(-.37,-.05)(-.38,-.35)
\put(.1,.3){\put(-.066,-.144){$\cdot$}}
\put(-.1,-.3){\put(-.066,-.144){$\cdot$}}
}

\put(0,15){ 
\qbezier[40](0,0)(.4,.25)(.25,.5) 
\qbezier[30](.25,.5)(0,.75)(-.25,.5)
\qbezier[40](0,0)(-.4,.25)(-.25,.5)
\qbezier[40](0,0)(.4,-.25)(.25,-.5) 
\qbezier[30](.25,-.5)(0,-.75)(-.25,-.5)
\qbezier[40](0,0)(-.4,-.25)(-.25,-.5)
\put(0,.35){\put(-.066,-.144){$\cdot$}}
\put(0,-.35){\put(-.066,-.144){$\cdot$}}
}

\put(3.75,15){ 
\qbezier[40](0,0)(.3,.25)(.12,.52) 
\qbezier[40](.12,.52)(-.23,.76)(-.38,.35)
\qbezier[40](0,0)(-.37,.05)(-.38,.35)
\qbezier[40](0,0)(-.3,-.25)(-.12,-.52) 
\qbezier[40](-.12,-.52)(.23,-.76)(.38,-.35)
\qbezier[40](0,0)(.37,-.05)(.38,-.35)
\put(-.1,.3){\put(-.066,-.144){$\cdot$}}
\put(.1,-.3){\put(-.066,-.144){$\cdot$}}
}

\put(7.5,15){ 
\qbezier[40](0,0)(.08,.35)(-.1,.5) 
\qbezier[40](-.1,.5)(-.6,.6)(-.5,.1)
\qbezier[40](0,0)(-.35,-.08)(-.5,.1)
\qbezier[40](0,0)(-.08,-.35)(.1,-.5) 
\qbezier[40](.1,-.5)(.6,-.6)(.5,-.1)
\qbezier[40](0,0)(.35,.08)(.5,-.1)
\put(-.25,.25){\put(-.066,-.144){$\cdot$}}
\put(.25,-.25){\put(-.066,-.144){$\cdot$}}
}

\put(7.5,8){ 
\qbezier[40](0,0)(.08,.35)(-.1,.5) 
\qbezier[40](-.1,.5)(-.6,.6)(-.5,.1)
\qbezier[40](0,0)(-.35,-.08)(-.5,.1)
\qbezier[40](0,0)(-.08,-.35)(.1,-.5) 
\qbezier[40](.1,-.5)(.6,-.6)(.5,-.1)
\qbezier[40](0,0)(.35,.08)(.5,-.1)
\put(-.25,.25){\put(-.066,-.144){$\cdot$}}
\put(.25,-.25){\put(-.066,-.144){$\cdot$}}
}

\put(3.75,0){ 
\qbezier[40](0,0)(-.25,-.3)(-.52,-.12) 
\qbezier[40](-.52,-.12)(-.76,.23)(-.35,.38)
\qbezier[40](0,0)(-.05,.37)(-.35,.38)
\qbezier[40](0,0)(.25,.3)(.52,.12) 
\qbezier[40](.52,.12)(.76,-.23)(.35,-.38)
\qbezier[40](0,0)(.05,-.37)(.35,-.38)
\put(-.3,.1){\put(-.066,-.144){$\cdot$}}
\put(.3,-.1){\put(-.066,-.144){$\cdot$}}
}

\put(7.5,0){ 
\qbezier[40](0,0)(.08,.35)(-.1,.5) 
\qbezier[40](-.1,.5)(-.6,.6)(-.5,.1)
\qbezier[40](0,0)(-.35,-.08)(-.5,.1)
\qbezier[40](0,0)(-.08,-.35)(.1,-.5) 
\qbezier[40](.1,-.5)(.6,-.6)(.5,-.1)
\qbezier[40](0,0)(.35,.08)(.5,-.1)
\put(-.25,.25){\put(-.066,-.144){$\cdot$}}
\put(.25,-.25){\put(-.066,-.144){$\cdot$}}
}

\put(0,0){ 
\qbezier[40](0,0)(.25,.4)(.5,.25) 
\qbezier[30](.5,.25)(.75,0)(.5,-.25)
\qbezier[40](0,0)(.25,-.4)(.5,-.25)
\qbezier[40](0,0)(-.25,.4)(-.5,.25) 
\qbezier[30](-.5,.25)(-.75,0)(-.5,-.25)
\qbezier[40](0,0)(-.25,-.4)(-.5,-.25)
\put(.35,0){\put(-.066,-.144){$\cdot$}}
\put(-.35,0){\put(-.066,-.144){$\cdot$}}
}

\put(-3.75,0){ 
\qbezier[40](0,0)(.25,-.3)(.52,-.12) 
\qbezier[40](.52,-.12)(.76,.23)(.35,.38)
\qbezier[40](0,0)(.05,.37)(.35,.38)
\qbezier[40](0,0)(-.25,.3)(-.52,.12) 
\qbezier[40](-.52,.12)(-.76,-.23)(-.35,-.38)
\qbezier[40](0,0)(-.05,-.37)(-.35,-.38)
\put(.3,.1){\put(-.066,-.144){$\cdot$}}
\put(-.3,-.1){\put(-.066,-.144){$\cdot$}}
}

\put(-7.5,0){ 
\qbezier[40](0,0)(-.08,.35)(.1,.5) 
\qbezier[40](.1,.5)(.6,.6)(.5,.1)
\qbezier[40](0,0)(.35,-.08)(.5,.1)
\qbezier[40](0,0)(.08,-.35)(-.1,-.5) 
\qbezier[40](-.1,-.5)(-.6,-.6)(-.5,-.1)
\qbezier[40](0,0)(-.35,.08)(-.5,-.1)
\put(.25,.25){\put(-.066,-.144){$\cdot$}}
\put(-.25,-.25){\put(-.066,-.144){$\cdot$}}
}

\put(-7.5,8){ 
\qbezier[40](0,0)(-.08,.35)(.1,.5) 
\qbezier[40](.1,.5)(.6,.6)(.5,.1)
\qbezier[40](0,0)(.35,-.08)(.5,.1)
\qbezier[40](0,0)(.08,-.35)(-.1,-.5) 
\qbezier[40](-.1,-.5)(-.6,-.6)(-.5,-.1)
\qbezier[40](0,0)(-.35,.08)(-.5,-.1)
\put(.25,.25){\put(-.066,-.144){$\cdot$}}
\put(-.25,-.25){\put(-.066,-.144){$\cdot$}}
}


\put(-7.5,-8){ 
\qbezier[40](0,0)(-.08,.35)(.1,.5) 
\qbezier[40](.1,.5)(.6,.6)(.5,.1)
\qbezier[40](0,0)(.35,-.08)(.5,.1)
\qbezier[40](0,0)(.08,-.35)(-.1,-.5) 
\qbezier[40](-.1,-.5)(-.6,-.6)(-.5,-.1)
\qbezier[40](0,0)(-.35,.08)(-.5,-.1)
\put(.25,.25){\put(-.066,-.144){$\cdot$}}
\put(-.25,-.25){\put(-.066,-.144){$\cdot$}}
}

\put(-7.5,-15){ 
\qbezier[40](0,0)(-.08,.35)(.1,.5) 
\qbezier[40](.1,.5)(.6,.6)(.5,.1)
\qbezier[40](0,0)(.35,-.08)(.5,.1)
\qbezier[40](0,0)(.08,-.35)(-.1,-.5) 
\qbezier[40](-.1,-.5)(-.6,-.6)(-.5,-.1)
\qbezier[40](0,0)(-.35,.08)(-.5,-.1)
\put(.25,.25){\put(-.066,-.144){$\cdot$}}
\put(-.25,-.25){\put(-.066,-.144){$\cdot$}}
}

\put(-3.75,-15){ 
\qbezier[40](0,0)(-.3,.25)(-.12,.52) 
\qbezier[40](-.12,.52)(.23,.76)(.38,.35)
\qbezier[40](0,0)(.37,.05)(.38,.35)
\qbezier[40](0,0)(.3,-.25)(.12,-.52) 
\qbezier[40](.12,-.52)(-.23,-.76)(-.38,-.35)
\qbezier[40](0,0)(-.37,-.05)(-.38,-.35)
\put(.1,.3){\put(-.066,-.144){$\cdot$}}
\put(-.1,-.3){\put(-.066,-.144){$\cdot$}}
}

\put(0,-15){ 
\qbezier[40](0,0)(.4,.25)(.25,.5) 
\qbezier[30](.25,.5)(0,.75)(-.25,.5)
\qbezier[40](0,0)(-.4,.25)(-.25,.5)
\qbezier[40](0,0)(.4,-.25)(.25,-.5) 
\qbezier[30](.25,-.5)(0,-.75)(-.25,-.5)
\qbezier[40](0,0)(-.4,-.25)(-.25,-.5)
\put(0,.35){\put(-.066,-.144){$\cdot$}}
\put(0,-.35){\put(-.066,-.144){$\cdot$}}
}

\put(3.75,-15){ 
\qbezier[40](0,0)(.3,.25)(.12,.52) 
\qbezier[40](.12,.52)(-.23,.76)(-.38,.35)
\qbezier[40](0,0)(-.37,.05)(-.38,.35)
\qbezier[40](0,0)(-.3,-.25)(-.12,-.52) 
\qbezier[40](-.12,-.52)(.23,-.76)(.38,-.35)
\qbezier[40](0,0)(.37,-.05)(.38,-.35)
\put(-.1,.3){\put(-.066,-.144){$\cdot$}}
\put(.1,-.3){\put(-.066,-.144){$\cdot$}}
}

\put(7.5,-15){ 
\qbezier[40](0,0)(.08,.35)(-.1,.5) 
\qbezier[40](-.1,.5)(-.6,.6)(-.5,.1)
\qbezier[40](0,0)(-.35,-.08)(-.5,.1)
\qbezier[40](0,0)(-.08,-.35)(.1,-.5) 
\qbezier[40](.1,-.5)(.6,-.6)(.5,-.1)
\qbezier[40](0,0)(.35,.08)(.5,-.1)
\put(-.25,.25){\put(-.066,-.144){$\cdot$}}
\put(.25,-.25){\put(-.066,-.144){$\cdot$}}
}

\put(7.5,-8){ 
\qbezier[40](0,0)(.08,.35)(-.1,.5) 
\qbezier[40](-.1,.5)(-.6,.6)(-.5,.1)
\qbezier[40](0,0)(-.35,-.08)(-.5,.1)
\qbezier[40](0,0)(-.08,-.35)(.1,-.5) 
\qbezier[40](.1,-.5)(.6,-.6)(.5,-.1)
\qbezier[40](0,0)(.35,.08)(.5,-.1)
\put(-.25,.25){\put(-.066,-.144){$\cdot$}}
\put(.25,-.25){\put(-.066,-.144){$\cdot$}}
}


\put(1.5,-7.5) {\vector(1,0){.78}} \put(1.5,-7.5) {\put(.5,-.35){\tiny$x$}}
\put(1.5,-7.5) {\vector(0,1){.8}} \put(1.5,-7.5) {\put(-.35,.55){\tiny$y$}}
\put(1.5,-7.5) {\put(-.066,-.144){$\cdot$}}
\put(1.5,-7.5) {\circle{.7}}
\put(1.5,-7.5) {
\path(-.84,-.65)(-.84,.65) \qbezier[15](-.85,.65)(-.85,.85)(-.65,.85)
\path(-.65,.85)(.65,.85) \qbezier[15](.65,.85)(.85,.85)(.85,.65)
\path(.86,.65)(.86,-.65) \qbezier[15](.85,-.65)(.85,-.85)(.65,-.85)
\path(.65,-.85)(-.65,-.85) \qbezier[15](-.65,-.85)(-.85,-.85)(-.85,-.65)
}

\put(-.6,11.2){$l_1'$}
\put(3.45,7.35){$l_2'$} \put(3.4,8.3){\tiny$A_3$}
\put(.3,4.5){$l_3'$}
\put(-3.45,7.35){$l_4'$} \put(-3.5,8.3){\tiny$A_3$}
\put(0,8){
\put(1.6,1){\tiny$A_2$}
\put(-2,1){\tiny$A_2$}
\put(1.6,-1.1){\tiny$A_2$}
\put(-2,-1.1){\tiny$A_2$}
}
\put(-5,10.85){$\Gamma^-$} \put(4.7,10.85){\tiny$A_{e,2}$}
\put(-5,-5.25){$\Gamma^+$} \put(4.7,-5.25){\tiny$A_{e,2}$}

\end{picture}
\end{center}

\vskip 0.5truecm
\centerline{\begin{tabular}{ll}
\bf Fig.~7: & Bifurcations of functions at a singularity $X_{e,5}^{++}$ with
              module $a\in\RR\setminus\{2,-2\}$ \\
            & in typical families of even functions, $0<a\ll1$
\end{tabular}
}
\end{figure}

\begin{figure}
\unitlength = 6.5mm
\begin{center}
\begin{picture}(20,32)(-10,-16) 
\path(-9,-16)(-9,16) \path(-9,16)(9,16) \path(9,16)(9,-16) \path(9,-16)(-9,-16)


\put(0,8) 
{
\qbezier[300](-4.44,2.7)(0,5.4)(4.44,2.7)
\qbezier[200](-4.44,2.7)(-7.56,0)(-4.44,-2.7)
\qbezier[200](4.44,2.7)(7.56,0)(4.44,-2.7)
\qbezier[300](-4.44,-2.7)(0,-5.4)(4.44,-2.7)
}
\put(0,-8) 
{
\qbezier[300](-4.44,2.7)(0,5.4)(4.44,2.7)
\qbezier[200](-4.44,2.7)(-7.56,0)(-4.44,-2.7)
\qbezier[200](4.44,2.7)(7.56,0)(4.44,-2.7)
\qbezier[300](-4.44,-2.7)(0,-5.4)(4.44,-2.7)
}
\put(0,-12.06) {\vector(1,0){2}} \put(0,-12.06) {\put(1.8,-.5){$\lam_2$}}
\put(0,-12.06) {\vector(0,1){2}} \put(0,-12.06) {\put(-.5,1.8){$\varphi$}}
\qbezier[400](-3.5,.6)(-8,-2.5)(-5.95,-8.5) \put(-5.95,-8.5){\put(-.145,-.14){$\bullet$}}
\qbezier[200](-3.5,.6)(-.5,2)(0,4.5) \put(0,4.5){\put(-.145,-.14){$\bullet$}}
\qbezier[200](3.5,.6)(.5,2)(0,4.5)
\qbezier[400](3.5,.6)(8,-2.5)(5.95,-8.5) \put(5.95,-8.5){\put(-.145,-.14){$\bullet$}}
\qbezier[400](-3.5,.7)(-5.9,2)(-5.95,7.5) \put(-5.95,7.5){\put(-.145,-.14){$\bullet$}}
\qbezier[200](-3.5,.7)(-.5,-1)(0,-3.4) \put(0,-3.4){\put(-.145,-.14){$\bullet$}}
\qbezier[200](3.5,.7)(.5,-1)(0,-3.4)
\qbezier[400](3.5,.7)(5.9,2)(5.95,7.5) \put(5.95,7.5){\put(-.145,-.14){$\bullet$}}
\qbezier[80](-5.95,-8.5)(-7,-3.5)(0,.6)
\qbezier[80](0,.6)(5.6,4)(5.95,7.5)
\qbezier[80](5.95,-8.5)(7,-3.5)(0,.6)
\qbezier[80](0,.6)(-5.6,4)(-5.95,7.5)
\qbezier[80](0,-3.4)(0,0)(0,4.5)


\put(-.002,.58){\put(-.066,-.144){$\cdot$} 
\circle{.6}
\path(-.4,-.4)(.4,.4)
\path(-.4,.4)(.4,-.4)
\put(.18,0){\put(-.066,-.144){$\cdot$}}
\put(-.18,0){\put(-.066,-.144){$\cdot$}}
\put(0,.18){\put(-.066,-.144){$\cdot$}}
\put(0,-.18){\put(-.066,-.144){$\cdot$}}
}

\put(0,2.85){\put(-.066,-.144){$\cdot$} 
\circle{.6}
\qbezier[80](.35,-.3)(-.05,0)(.35,.3) \put(.23,0){\put(-.066,-.144){$\cdot$}}
\qbezier[80](-.35,-.3)(.02,0)(-.35,.3) \put(-.23,0){\put(-.066,-.144){$\cdot$}}
\path(-.04,-.045)(.08,.07) 
\path(.06,-.045)(-.06,.07)
\qbezier[10](-.06,.06)(-.12,.12)(-.06,.18) 
\qbezier[10](-.06,.18)(0,.24)(.06,.18)
\qbezier[10](.06,.06)(.12,.12)(.06,.18)
\qbezier[10](.06,-.06)(.12,-.12)(.06,-.18) 
\qbezier[10](.06,-.18)(0,-.24)(-.06,-.18)
\qbezier[10](-.06,-.06)(-.12,-.12)(-.06,-.18)
\put(0,.12){\put(-.066,-.144){$\cdot$}}
\put(0,-.12){\put(-.066,-.144){$\cdot$}}
}

\put(1.1,1.3){ 
\qbezier[200](-.55,-.45)(-.2,-.1)(0,0)
\qbezier[200](.55,.45)(.2,.1)(0,0) 
\qbezier[200](.55,-.3)(.45,-.1)(.18,.15)
\qbezier[200](.18,.15)(-.15,.35)(0,0)
\qbezier[200](-.55,.3)(-.45,.1)(-.18,-.15)
\qbezier[200](-.18,-.15)(.15,-.35)(0,0)
\put(.05,.11){\put(-.066,-.144){$\cdot$}}
\put(-.05,-.11){\put(-.066,-.144){$\cdot$}}
\path(.1,.45)(-.22,.19) 
\path(-.3,.38)(-.12,.16)
\put(-.16,.22){
\qbezier[10](.04,-.08)(.05,-.18)(0,-.19) 
\qbezier[10](0,-.19)(-.05,-.22)(-.1,-.18)
\qbezier[10](-.08,-.04)(-.16,-.11)(-.1,-.18)
\put(-.04,-.11){\put(-.066,-.144){$\cdot$}}
}
\path(-.1,-.45)(.26,-.18) 
\path(.3,-.38)(.14,-.15)
\put(.16,-.22){
\qbezier[10](-.04,.08)(-.05,.18)(0,.19) 
\qbezier[10](0,.19)(.05,.22)(.1,.18)
\qbezier[10](.08,.04)(.16,.11)(.1,.18)
\put(.04,.11){\put(-.066,-.144){$\cdot$}}
}
}

\put(1.1,-.05){ 
\qbezier[200](.45,-.55)(.1,-.2)(0,0)
\qbezier[200](-.45,.55)(-.1,.2)(0,0) 
\qbezier[200](.3,.55)(.1,.45)(-.15,.18)
\qbezier[200](-.15,.18)(-.35,-.15)(0,0)
\qbezier[200](-.3,-.55)(-.1,-.45)(.15,-.18)
\qbezier[200](.15,-.18)(.35,.15)(0,0)
\put(-.11,.05){\put(-.066,-.144){$\cdot$}}
\put(.11,-.05){\put(-.066,-.144){$\cdot$}}
\path(.45,-.07)(.2,.23) 
\path(.38,.28)(.16,.12)
\put(.22,.16){
\qbezier[10](-.08,-.04)(-.18,-.05)(-.19,0) 
\qbezier[10](-.19,0)(-.22,.05)(-.18,.1)
\qbezier[10](-.04,.08)(-.11,.16)(-.18,.1)
\put(-.11,.04){\put(-.066,-.144){$\cdot$}}
}
\path(-.45,.1)(-.17,-.25) 
\path(-.38,-.3)(-.15,-.13)
\put(-.22,-.16){
\qbezier[10](.08,.04)(.18,.05)(.19,0) 
\qbezier[10](.19,0)(.22,-.05)(.18,-.1)
\qbezier[10](.04,-.08)(.11,-.16)(.18,-.1)
\put(.11,-.04){\put(-.066,-.144){$\cdot$}}
}
}

\put(0,-1.95){\put(-.066,-.144){$\cdot$} 
\circle{.6}
\qbezier[80](-.3,.35)(0,-.04)(.3,.35) \put(0,.23){\put(-.066,-.144){$\cdot$}}
\qbezier[80](-.3,-.35)(0,.04)(.3,-.35) \put(0,-.23){\put(-.066,-.144){$\cdot$}}
\path(-.045,-.06)(.07,.06) 
\path(-.045,.06)(.07,-.06)
\qbezier[10](.06,-.06)(.12,-.12)(.18,-.06) 
\qbezier[10](.18,-.06)(.24,0)(.18,.06)
\qbezier[10](.06,.06)(.12,.12)(.18,.06)
\qbezier[10](-.06,.06)(-.12,.12)(-.18,.06) 
\qbezier[10](-.18,.06)(-.24,0)(-.18,-.06)
\qbezier[10](-.06,-.06)(-.12,-.12)(-.18,-.06)
\put(.12,0){\put(-.066,-.144){$\cdot$}}
\put(-.12,0){\put(-.066,-.144){$\cdot$}}
}

\put(-.55,-.6){\put(-.066,-.144){$\cdot$}  
\qbezier[200](-.5,-.55)(1.05,.15)(-.5,.55) 
\qbezier[200](.5,.55)(-1.05,-.15)(.5,-.55)
\path(-.03,-.08)(.05,.07) 
\path(-.06,.04)(.09,-.04)
\qbezier[10](.08,-.04)(.18,-.05)(.19,0) 
\qbezier[10](.19,0)(.22,.05)(.18,.1)
\qbezier[10](.04,.08)(.11,.16)(.18,.1)
\qbezier[10](-.08,.04)(-.18,.05)(-.19,0) 
\qbezier[10](-.19,0)(-.22,-.05)(-.18,-.1)
\qbezier[10](-.04,-.08)(-.11,-.16)(-.18,-.1)
\put(.11,.04){\put(-.066,-.144){$\cdot$}}
\put(-.11,-.04){\put(-.066,-.144){$\cdot$}}
\path(-.5,-.1)(-.28,.34) 
\path(-.5,.4)(-.24,.24)
\put(-.32,.28){
\qbezier[10](.08,-.04)(.18,-.05)(.19,0) 
\qbezier[10](.19,0)(.22,.05)(.18,.1)
\qbezier[10](.04,.08)(.11,.16)(.18,.1)
\put(.11,.04){\put(-.066,-.144){$\cdot$}}
}
\path(.5,.1)(.28,-.38) 
\path(.5,-.3)(.25,-.24)
\put(.32,-.28){
\qbezier[10](-.08,.04)(-.18,.05)(-.19,0) 
\qbezier[10](-.19,0)(-.22,-.05)(-.18,-.1)
\qbezier[10](-.04,-.08)(-.11,-.16)(-.18,-.1)
\put(-.11,-.04){\put(-.066,-.144){$\cdot$}}
}
}

\put(-1.8,.58){\put(-.066,-.144){$\cdot$} 
\qbezier[100](.55,-.55)(0,-.45)(0,0) 
\qbezier[100](.55,.55)(.45,0)(0,0)
\qbezier[100](-.55,.55)(0,.45)(0,0)
\qbezier[100](-.55,-.55)(-.5,0)(0,0)
\path(-.1,.55)(.21,.26) 
\path(.4,.55)(.11,.22)
\put(.13,.3){
\qbezier[10](-.04,-.08)(-.05,-.18)(0,-.19) 
\qbezier[10](0,-.19)(.05,-.22)(.1,-.18)
\qbezier[10](.08,-.04)(.16,-.11)(.1,-.18)
\put(.04,-.11){\put(-.066,-.144){$\cdot$}}
}
\path(.1,-.55)(-.21,-.26) 
\path(-.4,-.55)(-.09,-.23)
\put(-.13,-.3){
\qbezier[10](.04,.08)(.05,.18)(0,.19) 
\qbezier[10](0,.19)(-.05,.22)(-.1,.18)
\qbezier[10](-.08,.04)(-.16,.11)(-.1,.18)
\put(-.04,.11){\put(-.066,-.144){$\cdot$}}
}
\path(-.55,-.1)(-.26,.19) 
\path(-.55,.4)(-.22,.09)
\put(-.3,.13){
\qbezier[10](.08,-.04)(.18,-.05)(.19,0) 
\qbezier[10](.19,0)(.22,.05)(.18,.1)
\qbezier[10](.04,.08)(.11,.16)(.18,.1)
\put(.11,.04){\put(-.066,-.144){$\cdot$}}
}
\path(.55,.1)(.26,-.23) 
\path(.55,-.4)(.23,-.09)
\put(.3,-.13){
\qbezier[10](-.08,.04)(-.18,.05)(-.19,0) 
\qbezier[10](-.19,0)(-.22,-.05)(-.18,-.1)
\qbezier[10](-.04,-.08)(-.11,-.16)(-.18,-.1)
\put(-.11,-.04){\put(-.066,-.144){$\cdot$}}
}
}

\put(-.55,1.8){\put(-.066,-.144){$\cdot$}  
\qbezier[200](.55,-.5)(-.15,1.05)(-.55,-.5) 
\qbezier[200](-.55,.5)(.15,-1.05)(.55,.5)
\path(.09,-.04)(-.06,.04) 
\path(-.03,-.08)(.05,.08)
\qbezier[10](.04,.08)(.05,.18)(0,.19) 
\qbezier[10](0,.19)(-.05,.22)(-.1,.18)
\qbezier[10](-.08,.04)(-.16,.11)(-.1,.18)
\qbezier[10](-.04,-.08)(-.05,-.18)(0,-.19) 
\qbezier[10](0,-.19)(.05,-.22)(.1,-.18)
\qbezier[10](.08,-.04)(.16,-.11)(.1,-.18)
\put(-.04,.11){\put(-.066,-.144){$\cdot$}}
\put(.04,-.11){\put(-.066,-.144){$\cdot$}}
\path(.1,-.55)(-.36,-.28) 
\path(-.4,-.55)(-.24,-.26)
\put(-.28,-.32){
\qbezier[10](.04,.08)(.05,.18)(0,.19) 
\qbezier[10](0,.19)(-.05,.22)(-.1,.18)
\qbezier[10](-.08,.04)(-.16,.11)(-.1,.18)
\put(-.04,.11){\put(-.066,-.144){$\cdot$}}
}
\path(-.1,.55)(.36,.28) 
\path(.4,.55)(.26,.25)
\put(.28,.32){
\qbezier[10](-.04,-.08)(-.05,-.18)(0,-.19) 
\qbezier[10](0,-.19)(.05,-.22)(.1,-.18)
\qbezier[10](.08,-.04)(.16,-.11)(.1,-.18)
\put(.04,-.11){\put(-.066,-.144){$\cdot$}}
}
}


\put(-3.75,8){\put(-.066,-.144){$\cdot$} 
\qbezier[200](-.75,.7)(.3,-1.1)(.75,.7) 
\qbezier[200](.75,-.7)(-.3,1.1)(-.75,-.7)
\qbezier[40](-.75,-.2)(-.55,.2)(-.75,.4) 
\qbezier[40](.75,.2)(.55,-.2)(.75,-.4)
\qbezier[60](-.55,.7)(.2,.05)(.55,.7)
\qbezier[60](.55,-.7)(-.2,-.05)(-.55,-.7)
}

\put(0,8){\put(-.066,-.144){$\cdot$} 
\qbezier[200](-.75,-.7)(0,1)(.75,-.7) 
\qbezier[200](-.75,.7)(0,-1)(.75,.7)
\qbezier[60](-.75,-.3)(-.5,0)(-.75,.3) 
\qbezier[60](.75,-.3)(.5,0)(.75,.3)
\qbezier[80](-.55,-.7)(0,-.05)(.55,-.7)
\qbezier[80](-.55,.7)(0,.05)(.55,.7)
}

\put(3.75,8){\put(-.066,-.144){$\cdot$} 
\qbezier[200](-.75,-.7)(.3,1.1)(.75,-.7) 
\qbezier[200](.75,.7)(-.3,-1.1)(-.75,.7)
\qbezier[40](-.75,.2)(-.55,-.2)(-.75,-.4) 
\qbezier[40](.75,-.2)(.55,.2)(.75,.4)
\qbezier[60](-.55,-.7)(.2,-.05)(.55,-.7)
\qbezier[60](.55,.7)(-.2,.05)(-.55,.7)
}


\put(-3.75,-8){\put(-.066,-.144){$\cdot$} 
\qbezier[200](-.7,-.75)(1.1,.3)(-.7,.75) 
\qbezier[200](.7,.75)(-1.1,-.3)(.7,-.75)
\qbezier[40](.2,-.75)(-.2,-.55)(-.4,-.75) 
\qbezier[40](-.2,.75)(.2,.55)(.4,.75)
\qbezier[60](-.7,-.55)(-.05,.2)(-.7,.55)
\qbezier[60](.7,.55)(.05,-.2)(.7,-.55)
}

\put(0,-8){\put(-.066,-.144){$\cdot$} 
\qbezier[200](-.7,-.75)(1,0)(-.7,.75) 
\qbezier[200](.7,-.75)(-1,0)(.7,.75)
\qbezier[60](-.3,-.75)(0,-.5)(.3,-.75) 
\qbezier[60](-.3,.75)(0,.5)(.3,.75)
\qbezier[80](-.7,-.55)(-.05,0)(-.7,.55)
\qbezier[80](.7,-.55)(.05,0)(.7,.55)
}

\put(3.75,-8){\put(-.066,-.144){$\cdot$} 
\qbezier[200](.7,-.75)(-1.1,.3)(.7,.75) 
\qbezier[200](-.7,.75)(1.1,-.3)(-.7,-.75)
\qbezier[40](-.2,-.75)(.2,-.55)(.4,-.75) 
\qbezier[40](.2,.75)(-.2,.55)(-.4,.75)
\qbezier[60](.7,-.55)(.05,.2)(.7,.55)
\qbezier[60](-.7,.55)(-.05,-.2)(-.7,-.55)
}


\put(-4.1,2.3){\put(-.066,-.144){$\cdot$} 
\qbezier[100](.55,.45)(.45,0)(0,0) 
\qbezier[100](-.55,.45)(0,.45)(0,0)
\qbezier[100](-.55,-.45)(-.45,0)(0,0)
\qbezier[100](.55,-.45)(0,-.5)(0,0)
\path(-.55,-.1)(-.26,.19) 
\path(-.55,.3)(-.22,.09)
\put(-.3,.13){
\qbezier[10](.08,-.04)(.18,-.05)(.19,0) 
\qbezier[10](.19,0)(.22,.05)(.18,.1)
\qbezier[10](.04,.08)(.11,.16)(.18,.1)
\put(.11,.04){\put(-.066,-.144){$\cdot$}}
}
\path(.55,.1)(.26,-.23) 
\path(.55,-.3)(.23,-.09)
\put(.3,-.13){
\qbezier[10](-.08,.04)(-.18,.05)(-.19,0) 
\qbezier[10](-.19,0)(-.22,-.05)(-.18,-.1)
\qbezier[10](-.04,-.08)(-.11,-.16)(-.18,-.1)
\put(-.11,-.04){\put(-.066,-.144){$\cdot$}}
}
}

\put(-2.8,2.65){ 
\qbezier[200](.55,-.45)(.2,-.1)(0,0)
\qbezier[200](-.55,.45)(-.2,.1)(0,0) 
\qbezier[200](-.55,-.3)(-.45,-.1)(-.18,.15)
\qbezier[200](-.18,.15)(.15,.35)(0,0)
\qbezier[200](.55,.3)(.45,.1)(.18,-.15)
\qbezier[200](.18,-.15)(-.15,-.35)(0,0)
\put(-.05,.11){\put(-.066,-.144){$\cdot$}}
\put(.05,-.11){\put(-.066,-.144){$\cdot$}}
}

\put(-2,3.55){\put(-.066,-.144){$\cdot$}  
\qbezier[200](.55,-.5)(-.15,1.15)(-.55,-.5) 
\qbezier[200](-.55,.5)(.15,-1.15)(.55,.5)
\path(.09,-.04)(-.06,.04) 
\path(-.03,-.08)(.05,.08)
\qbezier[10](.04,.08)(.05,.18)(0,.19) 
\qbezier[10](0,.19)(-.05,.22)(-.1,.18)
\qbezier[10](-.08,.04)(-.16,.11)(-.1,.18)
\qbezier[10](-.04,-.08)(-.05,-.18)(0,-.19) 
\qbezier[10](0,-.19)(.05,-.22)(.1,-.18)
\qbezier[10](.08,-.04)(.16,-.11)(.1,-.18)
\put(-.04,.11){\put(-.066,-.144){$\cdot$}}
\put(.04,-.11){\put(-.066,-.144){$\cdot$}}
}

\put(2,3.55){\put(-.066,-.144){$\cdot$}  
\qbezier[200](-.55,-.5)(.15,1.15)(.55,-.5) 
\qbezier[200](.55,.5)(-.15,-1.15)(-.55,.5)
\path(-.07,-.04)(.08,.04) 
\path(.05,-.07)(-.03,.08)
\qbezier[10](-.04,.08)(-.05,.18)(0,.19) 
\qbezier[10](0,.19)(.05,.22)(.1,.18)
\qbezier[10](.08,.04)(.16,.11)(.1,.18)
\qbezier[10](.04,-.08)(.05,-.18)(0,-.19) 
\qbezier[10](0,-.19)(-.05,-.22)(-.1,-.18)
\qbezier[10](-.08,-.04)(-.16,-.11)(-.1,-.18)
\put(.04,.11){\put(-.066,-.144){$\cdot$}}
\put(-.04,-.11){\put(-.066,-.144){$\cdot$}}
}

\put(2.8,2.65){ 
\qbezier[200](-.55,-.45)(-.2,-.1)(0,0)
\qbezier[200](.55,.45)(.2,.1)(0,0) 
\qbezier[200](.55,-.3)(.45,-.1)(.18,.15)
\qbezier[200](.18,.15)(-.15,.35)(0,0)
\qbezier[200](-.55,.3)(-.45,.1)(-.18,-.15)
\qbezier[200](-.18,-.15)(.15,-.35)(0,0)
\put(.05,.11){\put(-.066,-.144){$\cdot$}}
\put(-.05,-.11){\put(-.066,-.144){$\cdot$}}
}

\put(4.1,2.3){\put(-.066,-.144){$\cdot$} 
\qbezier[100](-.55,.45)(-.45,0)(0,0) 
\qbezier[100](.55,.45)(0,.45)(0,0)
\qbezier[100](.55,-.45)(.45,0)(0,0)
\qbezier[100](-.55,-.45)(0,-.5)(0,0)
\path(.55,-.1)(.28,.2) 
\path(.55,.3)(.24,.09)
\put(.3,.13){
\qbezier[10](-.08,-.04)(-.18,-.05)(-.19,0) 
\qbezier[10](-.19,0)(-.22,.05)(-.18,.1)
\qbezier[10](-.04,.08)(-.11,.16)(-.18,.1)
\put(-.11,.04){\put(-.066,-.144){$\cdot$}}
}
\path(-.55,.1)(-.25,-.22) 
\path(-.55,-.3)(-.23,-.1)
\put(-.3,-.13){
\qbezier[10](.08,.04)(.18,.05)(.19,0) 
\qbezier[10](.19,0)(.22,-.05)(.18,-.1)
\qbezier[10](.04,-.08)(.11,-.16)(.18,-.1)
\put(.11,-.04){\put(-.066,-.144){$\cdot$}}
}
}


\put(-4.5,-1.3){\put(-.066,-.144){$\cdot$} 
\qbezier[100](.45,-.55)(0,-.45)(0,0) 
\qbezier[100](.45,.55)(.45,0)(0,0)
\qbezier[100](-.45,.55)(0,.45)(0,0)
\qbezier[100](-.45,-.55)(-.5,0)(0,0)
\path(-.1,.55)(.21,.26) 
\path(.3,.55)(.11,.22)
\put(.13,.3){
\qbezier[10](-.04,-.08)(-.05,-.18)(0,-.19) 
\qbezier[10](0,-.19)(.05,-.22)(.1,-.18)
\qbezier[10](.08,-.04)(.16,-.11)(.1,-.18)
\put(.04,-.11){\put(-.066,-.144){$\cdot$}}
}
\path(.1,-.55)(-.21,-.26) 
\path(-.3,-.55)(-.09,-.23)
\put(-.13,-.3){
\qbezier[10](.04,.08)(.05,.18)(0,.19) 
\qbezier[10](0,.19)(-.05,.22)(-.1,.18)
\qbezier[10](-.08,.04)(-.16,.11)(-.1,.18)
\put(-.04,.11){\put(-.066,-.144){$\cdot$}}
}
}

\put(-3.5,-2){ 
\qbezier[200](-.45,-.55)(-.1,-.2)(0,0)
\qbezier[200](.45,.55)(.1,.2)(0,0) 
\qbezier[200](-.3,.55)(-.1,.45)(.15,.18)
\qbezier[200](.15,.18)(.35,-.15)(0,0)
\qbezier[200](.3,-.55)(.1,-.45)(-.15,-.18)
\qbezier[200](-.15,-.18)(-.35,.15)(0,0)
\put(.11,.05){\put(-.066,-.144){$\cdot$}}
\put(-.11,-.05){\put(-.066,-.144){$\cdot$}}
}

\put(-2.5,-3.2){\put(-.066,-.144){$\cdot$}  
\qbezier[200](-.5,-.55)(1.15,.15)(-.5,.55) 
\qbezier[200](.5,.55)(-1.15,-.15)(.5,-.55)
\path(-.03,-.08)(.05,.07) 
\path(-.06,.04)(.09,-.04)
\qbezier[10](.08,-.04)(.18,-.05)(.19,0) 
\qbezier[10](.19,0)(.22,.05)(.18,.1)
\qbezier[10](.04,.08)(.11,.16)(.18,.1)
\qbezier[10](-.08,.04)(-.18,.05)(-.19,0) 
\qbezier[10](-.19,0)(-.22,-.05)(-.18,-.1)
\qbezier[10](-.04,-.08)(-.11,-.16)(-.18,-.1)
\put(.11,.04){\put(-.066,-.144){$\cdot$}}
\put(-.11,-.04){\put(-.066,-.144){$\cdot$}}
}

\path(0,-3.76)(.5,-3.41)
\put(.8,-3.3){\put(-.066,-.144){$\cdot$}  
\qbezier[200](-.4,-.5)(1.1,0)(-.4,.5) 
\qbezier[200](.4,.5)(-1.1,0)(.4,-.5)
\path(-.045,-.06)(.07,.06) 
\path(-.045,.06)(.07,-.06)
\qbezier[10](.06,-.06)(.12,-.12)(.18,-.06) 
\qbezier[10](.18,-.06)(.24,0)(.18,.06)
\qbezier[10](.06,.06)(.12,.12)(.18,.06)
\qbezier[10](-.06,.06)(-.12,.12)(-.18,.06) 
\qbezier[10](-.18,.06)(-.24,0)(-.18,-.06)
\qbezier[10](-.06,-.06)(-.12,-.12)(-.18,-.06)
\put(.12,0){\put(-.066,-.144){$\cdot$}}
\put(-.12,0){\put(-.066,-.144){$\cdot$}}
}

\put(2.5,-3.2){\put(-.066,-.144){$\cdot$}  
\qbezier[200](.5,-.55)(-1.15,.15)(.5,.55) 
\qbezier[200](-.5,.55)(1.15,-.15)(-.5,-.55)
\path(.05,-.08)(-.03,.07) 
\path(.08,.04)(-.07,-.04)
\qbezier[10](-.08,-.04)(-.18,-.05)(-.19,0) 
\qbezier[10](-.19,0)(-.22,.05)(-.18,.1)
\qbezier[10](-.04,.08)(-.11,.16)(-.18,.1)
\qbezier[10](.08,.04)(.18,.05)(.19,0) 
\qbezier[10](.19,0)(.22,-.05)(.18,-.1)
\qbezier[10](.04,-.08)(.11,-.16)(.18,-.1)
\put(-.11,.04){\put(-.066,-.144){$\cdot$}}
\put(.11,-.04){\put(-.066,-.144){$\cdot$}}
}

\put(3.5,-2){ 
\qbezier[200](.45,-.55)(.1,-.2)(0,0)
\qbezier[200](-.45,.55)(-.1,.2)(0,0) 
\qbezier[200](.3,.55)(.1,.45)(-.15,.18)
\qbezier[200](-.15,.18)(-.35,-.15)(0,0)
\qbezier[200](-.3,-.55)(-.1,-.45)(.15,-.18)
\qbezier[200](.15,-.18)(.35,.15)(0,0)
\put(-.11,.05){\put(-.066,-.144){$\cdot$}}
\put(.11,-.05){\put(-.066,-.144){$\cdot$}}
}

\put(4.5,-1.3){\put(-.066,-.144){$\cdot$} 
\qbezier[100](-.45,-.55)(0,-.45)(0,0) 
\qbezier[100](-.45,.55)(-.45,0)(0,0)
\qbezier[100](.45,.55)(0,.45)(0,0)
\qbezier[100](.45,-.55)(.5,0)(0,0)
\path(.1,.55)(-.19,.26) 
\path(-.3,.55)(-.09,.24)
\put(-.13,.3){
\qbezier[10](.04,-.08)(.05,-.18)(0,-.19) 
\qbezier[10](0,-.19)(-.05,-.22)(-.1,-.18)
\qbezier[10](-.08,-.04)(-.16,-.11)(-.1,-.18)
\put(-.04,-.11){\put(-.066,-.144){$\cdot$}}
}
\path(-.1,-.55)(.23,-.26) 
\path(.3,-.55)(.11,-.23)
\put(.13,-.3){
\qbezier[10](-.04,.08)(-.05,.18)(0,.19) 
\qbezier[10](0,.19)(.05,.22)(.1,.18)
\qbezier[10](.08,.04)(.16,.11)(.1,.18)
\put(.04,.11){\put(-.066,-.144){$\cdot$}}
}
}


\put(-7.7,14.5){\put(-.066,-.144){$\cdot$} 
\qbezier[100](.75,-.75)(0,-.6)(0,0) 
\qbezier[100](.75,.75)(.6,0)(0,0)
\qbezier[100](-.75,.75)(0,.6)(0,0)
\qbezier[100](-.75,-.75)(-.6,0)(0,0)
\qbezier[80](.2,-.75)(-.4,.1)(-.6,-.75) 
\qbezier[80](.6,.75)(.4,-.1)(-.2,.75)
\qbezier[80](.75,-.6)(-.1,-.4)(.75,.2)
\qbezier[80](-.75,-.2)(.1,.4)(-.75,.6)
}

\put(-3.85,14.5){\put(-.066,-.144){$\cdot$} 
\qbezier[100](.75,-.75)(.2,-.55)(0,0) 
\qbezier[100](.75,.75)(.55,.2)(0,0)
\qbezier[100](-.75,.75)(-.2,.55)(0,0)
\qbezier[100](-.75,-.75)(-.55,-.2)(0,0)
\qbezier[80](.4,-.75)(-.275,0)(-.6,-.75) 
\qbezier[80](.6,.75)(.275,0)(-.4,.75)
\qbezier[80](.75,-.6)(0,-.275)(.75,.4)
\qbezier[80](-.75,-.4)(0,.275)(-.75,.6)
}

\put(0,14.5) {\put(-.066,-.144){$\cdot$} 
\path(-.75,-.75)(.75,.75) \path(-.75,.75)(.75,-.75)
\qbezier[80](-.65,-.75)(0,-.1)(.65,-.75) 
\qbezier[80](-.65,.75)(0,.1)(.65,.75)
\qbezier[80](-.75,-.65)(-.1,0)(-.75,.65)
\qbezier[80](.75,-.65)(.1,0)(.75,.65)
}

\put(3.85,14.5){\put(-.066,-.144){$\cdot$} 
\qbezier[100](-.75,-.75)(-.2,-.55)(0,0) 
\qbezier[100](-.75,.75)(-.55,.2)(0,0)
\qbezier[100](.75,.75)(.2,.55)(0,0)
\qbezier[100](.75,-.75)(.55,-.2)(0,0)
\qbezier[80](-.4,-.75)(.275,0)(.6,-.75) 
\qbezier[80](-.6,.75)(-.275,0)(.4,.75)
\qbezier[80](-.75,-.6)(0,-.275)(-.75,.4)
\qbezier[80](.75,-.4)(0,.275)(.75,.6)
}

\put(7.7,14.5){\put(-.066,-.144){$\cdot$} 
\qbezier[100](-.75,-.75)(0,-.6)(0,0) 
\qbezier[100](-.75,.75)(-.6,0)(0,0)
\qbezier[100](.75,.75)(0,.6)(0,0)
\qbezier[100](.75,-.75)(.6,0)(0,0)
\qbezier[80](-.2,-.75)(.4,.1)(.6,-.75) 
\qbezier[80](-.6,.75)(-.4,-.1)(.2,.75)
\qbezier[80](-.75,-.6)(.1,-.4)(-.75,.2)
\qbezier[80](.75,-.2)(-.1,.4)(.75,.6)
}


\put(-7.7,-14.5){\put(-.066,-.144){$\cdot$} 
\qbezier[100](.75,-.75)(0,-.6)(0,0) 
\qbezier[100](.75,.75)(.6,0)(0,0)
\qbezier[100](-.75,.75)(0,.6)(0,0)
\qbezier[100](-.75,-.75)(-.6,0)(0,0)
\qbezier[80](.2,-.75)(-.4,.1)(-.6,-.75) 
\qbezier[80](.6,.75)(.4,-.1)(-.2,.75)
\qbezier[80](.75,-.6)(-.1,-.4)(.75,.2)
\qbezier[80](-.75,-.2)(.1,.4)(-.75,.6)
}

\put(-3.85,-14.5){\put(-.066,-.144){$\cdot$} 
\qbezier[100](.75,-.75)(.2,-.55)(0,0) 
\qbezier[100](.75,.75)(.55,.2)(0,0)
\qbezier[100](-.75,.75)(-.2,.55)(0,0)
\qbezier[100](-.75,-.75)(-.55,-.2)(0,0)
\qbezier[80](.4,-.75)(-.275,0)(-.6,-.75) 
\qbezier[80](.6,.75)(.275,0)(-.4,.75)
\qbezier[80](.75,-.6)(0,-.275)(.75,.4)
\qbezier[80](-.75,-.4)(0,.275)(-.75,.6)
}

\put(0,-14.5) {\put(-.066,-.144){$\cdot$} 
\path(-.75,-.75)(.75,.75) \path(-.75,.75)(.75,-.75)
\qbezier[80](-.65,-.75)(0,-.1)(.65,-.75) 
\qbezier[80](-.65,.75)(0,.1)(.65,.75)
\qbezier[80](-.75,-.65)(-.1,0)(-.75,.65)
\qbezier[80](.75,-.65)(.1,0)(.75,.65)
}
\put(0,-14.5) {\vector(1,0){.78}} \put(0,-14.5) {\put(.6,-.35){\tiny$x$}}
\put(0,-14.5) {\vector(0,1){.8}} \put(0,-14.5) {\put(-.35,.65){\tiny$y$}}
\put(0,-14.5) {
\path(-.84,-.65)(-.84,.65) \qbezier[15](-.85,.65)(-.85,.85)(-.65,.85) 
\path(-.65,.85)(.65,.85) \qbezier[15](.65,.85)(.85,.85)(.85,.65)
\path(.86,.65)(.86,-.65) \qbezier[15](.85,-.65)(.85,-.85)(.65,-.85)
\path(.65,-.85)(-.65,-.85) \qbezier[15](-.65,-.85)(-.85,-.85)(-.85,-.65)
}

\put(3.85,-14.5){\put(-.066,-.144){$\cdot$} 
\qbezier[100](-.75,-.75)(-.2,-.55)(0,0) 
\qbezier[100](-.75,.75)(-.55,.2)(0,0)
\qbezier[100](.75,.75)(.2,.55)(0,0)
\qbezier[100](.75,-.75)(.55,-.2)(0,0)
\qbezier[80](-.4,-.75)(.275,0)(.6,-.75) 
\qbezier[80](-.6,.75)(-.275,0)(.4,.75)
\qbezier[80](-.75,-.6)(0,-.275)(-.75,.4)
\qbezier[80](.75,-.4)(0,.275)(.75,.6)
}

\put(7.7,-14.5){\put(-.066,-.144){$\cdot$} 
\qbezier[100](-.75,-.75)(0,-.6)(0,0) 
\qbezier[100](-.75,.75)(-.6,0)(0,0)
\qbezier[100](.75,.75)(0,.6)(0,0)
\qbezier[100](.75,-.75)(.6,0)(0,0)
\qbezier[80](-.2,-.75)(.4,.1)(.6,-.75) 
\qbezier[80](-.6,.75)(-.4,-.1)(.2,.75)
\qbezier[80](-.75,-.6)(.1,-.4)(-.75,.2)
\qbezier[80](.75,-.2)(-.1,.4)(.75,.6)
}


\multiput(-7.7,8)(0,-8){3}{\put(-.066,-.144){$\cdot$} 
\qbezier[100](.75,-.75)(0,-.6)(0,0) 
\qbezier[100](.75,.75)(.6,0)(0,0)
\qbezier[100](-.75,.75)(0,.6)(0,0)
\qbezier[100](-.75,-.75)(-.6,0)(0,0)
\qbezier[80](.2,-.75)(-.4,.1)(-.6,-.75) 
\qbezier[80](.6,.75)(.4,-.1)(-.2,.75)
\qbezier[80](.75,-.6)(-.1,-.4)(.75,.2)
\qbezier[80](-.75,-.2)(.1,.4)(-.75,.6)
}

\multiput(7.7,8)(0,-8){3}{\put(-.066,-.144){$\cdot$} 
\qbezier[100](-.75,-.75)(0,-.6)(0,0) 
\qbezier[100](-.75,.75)(-.6,0)(0,0)
\qbezier[100](.75,.75)(0,.6)(0,0)
\qbezier[100](.75,-.75)(.6,0)(0,0)
\qbezier[80](-.2,-.75)(.4,.1)(.6,-.75) 
\qbezier[80](-.6,.75)(-.4,-.1)(.2,.75)
\qbezier[80](-.75,-.6)(.1,-.4)(-.75,.2)
\qbezier[80](.75,-.2)(-.1,.4)(.75,.6)
}

\put(.3,4.5){$l_1'$} \put(-.7,4.7){\tiny$A_3$}
\put(-.6,-3.6){$l_2'$}
\put(4.2,.5){$l_{23}'$}
\put(-4.6,.5){$l_{14}'$}
\put(-6.6,-8.6){$l_1^+$} \put(-7.3,-4.5){$\Gamma_1'$} \put(-5.7,-8.6){\tiny$A_{e,3}$}
\put(6.3,7.5){$l_1^-$} \put(6.1,4){$\Gamma_3'$} \put(4.9,7.5){\tiny$A_{e,3}$}
\put(6.3,-8.6){$l_2^+$} \put(6.9,-4.5){$\Gamma_2'$} \put(4.9,-8.6){\tiny$A_{e,3}$}
\put(-6.6,7.5){$l_2^-$} \put(-6.4,4){$\Gamma_4'$} \put(-5.7,7.5){\tiny$A_{e,3}$}
\put(0,-4){
\put(1.3,2.5){\tiny$A_2$}
\put(-1.7,2.5){\tiny$A_2$}
}
\put(2.2,-5.1){\tiny$A_{e,2}$}
\put(-5,10.85){$\Gamma^-$} \put(4.7,10.85){\tiny$A_{e,2}$}
\put(-5,-11.5){$\Gamma^+$} \put(4.7,-11.3){\tiny$A_{e,2}$}
\end{picture}
\end{center}

\vskip 0.5truecm
\centerline{\begin{tabular}{ll}
\bf Fig.~8: & Bifurcations of functions at a singularity $X_{e,5}^{+-}$ with
              module $a\in\RR$ \\
            & in typical families of even functions, $0<a\ll1$
\end{tabular}
}
\end{figure}

In the cases $++$ and $+-$ with $0<a\ll1$, bifurcations of level lines of
functions in $\nu=2$ variables are shown on Fig.~7 and~8, respectively.
These figures do not show the whole caustic, but only its intersection with the
cylinder $\lam_1=\varepsilon\cos\varphi$, $\lam_3=\varepsilon\sin\varphi$,
$-\frac{\pi}{4}\le\varphi\le\frac{7\pi}{4}$ for a small positive constant
$\varepsilon$. In the both cases, the bifurcation diagram is invariant under
the transformation $(\lam_2,\varphi,x,y)\mapsto(-\lam_2,\varphi,-x,y)$.
One can see from these figures that, if the parameter $\lam$ goes around a loop
which envelops the cone $\Gamma$, the following transformation happens.
In the case $++$, the twin critical point $(k,-k)$ (corresponding to a pair of
minima of an even function) transforms into $(-k,k)$. In the case $+-$, each of
two separatrices of the saddle point $0$ turns by $\pi$.

\begin{note} \label {note:sym} \rm
The conic surface $\Gamma'$ is symmetric. Namely, it is invariant under
reflections with respect to one or two planes which leave the cone $\Gamma$
invariant. In our model, it is the plane $\lambda_2=0$ in the case $+-$
and two planes $\lambda_2=0$ and $\lambda_1=\lambda_3$ in each of the cases
$++$ and $--$. Besides, in the case $+-$, the caustic surfaces corresponding
to opposite values of the parameter $a$ are obtained from each other by
reflection with respect to the plane $\lambda_1=-\lambda_3$. Furthermore,
the caustic surfaces in the case $++$ and the case $--$, corresponding to
opposite values of the parameter $a$, are obtained from each other
by the reflection with respect to this plane.
\end{note}

\begin{note} \label {note:edges} \rm
If the parameter $a$ continuously varies in the domain of admissible values
($a\ne\pm2$ in the cases $++$ and $--$), the conic surface $\Gamma'$ is
deformed. In particular, locations of the cuspidal edges $l_1',l_2',l_3',l_4'$
(in the cases $++$ and $--$) and resp.\ $l_1',l_2'$ (in the case $+-$) are
continuously deformed with respect to $a$. Moreover, each of these cuspidal
edges is preserved for all admissible values of $a$ by just deforming
continuously on $a$, with the edge remaining in the corresponding plane of
symmetry:
\\In the cases $++$ and $--$, the cuspidal edges $l_1'$ and $l_3'$ lie in
the plane $\lambda_2=0$, while the cuspidal edges $l_2'$ and $l_4'$ lie in
the plane $\lambda_1=\lambda_3$. More precisely: in the case $++$, the cuspidal
edges $l_1'$ and $l_3'$ correspond to the twin critical points
$\left( (\sqrt{-\frac{\lam_1}{2}},0),(-\sqrt{-\frac{\lam_1}{2}},0)\right)$ and
resp.\
$\left( (0,\sqrt{-\frac{\lam_3}{2}}),(0,-\sqrt{-\frac{\lam_3}{2}})\right)$
of type $A_3$. They are defined by
$\lam_3=\frac{a}{2}\lam_1$, $\lam_1<\lam_2=0$ and resp.\
$\lam_1=\frac{a}{2}\lam_3$, $\lam_3<\lam_2=0$.
Furthermore, the cuspidal edges $l_2'$ and $l_4'$ correspond to the twin
critical points
$\left( (\sqrt{-\frac{\lam_1}{2}},\sqrt{-\frac{\lam_1}{2}}),
(-\sqrt{-\frac{\lam_1}{2}},-\sqrt{-\frac{\lam_1}{2}})\right)$ and resp.\
$\left( (\sqrt{-\frac{\lam_1}{2}},-\sqrt{-\frac{\lam_1}{2}}),
(-\sqrt{-\frac{\lam_1}{2}},\sqrt{-\frac{\lam_1}{2}})\right)$
of type $A_3$. They are defined by
$\lam_2=\frac{a-2}{2}\lam_1$, $\lam_1=\lam_3<0$ and resp.\
$\lam_2=\frac{2-a}{2}\lam_1$, $\lam_1=\lam_3<0$.
\\In the case $+-$, the cuspidal edges $l_1'$ and $l_2'$ lie in the plane
$\lambda_2=0$ and correspond to the twin critical points
$\left( (\sqrt{-\frac{\lam_1}{2}},0),(-\sqrt{-\frac{\lam_1}{2}},0)\right)$ and
resp.\
$\left( (0,\sqrt{\frac{\lam_3}{2}}),(0,-\sqrt{\frac{\lam_3}{2}})\right)$
of type $A_3$. They are defined by
$\lam_3=\frac{a}{2}\lam_1$, $\lam_1<\lam_2=0$ and resp.\
$\lam_1=-\frac{a}{2}\lam_3$, $\lam_3>\lam_2=0$.

>From the indicated explicit presentation of the cuspidal edges $l_1',l_2'$
(in the case $+-$) and $l_1',l_2',l_3',l_4'$ (in the cases $++$ and $--$),
we find their location with respect to the cone $\Gamma$ for any $a$:
\\In the case $++$ (the case $--$ is analogous due to~\ref {note:sym}),
all of the cuspidal edges $l_1',l_2',l_3',l_4'$ remain inside the domain
$\Omega^-$ for $0<a<2$ and for $2<a<6$, while they are outside the domains
$\Omega^+$ and $\Omega^-$ for $a<-2$ (moreover, all these cuspidal edges
``tend'' to the axis of the cone $\Gamma$ as $a\to2$).
\\In the case $++$ with $-2<a<0$, the cuspidal edges $l_2'$ and $l_4'$ lie
inside the domain $\Omega^-$, while the cuspidal edges $l_1'$ and $l_3'$ lie
outside the domains $\Omega^+$ and $\Omega^-$.
\\In the case $++$ with $a>6$, the cuspidal edges $l_1'$ and $l_3'$ lie inside
the domain $\Omega^-$, while the cuspidal edges $l_2'$ and $l_4'$ lie outside
the domains $\Omega^+$ and $\Omega^-$.
\\In the case $+-$ with $a>0$, the cuspidal edge $l_1'$ lies inside the domain
$\Omega^-$, while the cuspidal edge $l_2'$ lies outside the domains $\Omega^+$
and $\Omega^-$.
\\In the case $+-$ with $a<0$, the cuspidal edge $l_2'$ lies inside the domain
$\Omega^+$, while the cuspidal edge $l_1'$ lies outside the domains $\Omega^+$
and $\Omega^-$.

Unfortunately, the authors do not know whether the form of the conic surface
$\Gamma'$ (in particular, the number of its cuspidal edges) and its location
with respect to the cone $\Gamma$ change under a continuous change of the
value of $a$. In particular, we do not know whether {\it bifurcations of
cuspidal edges} happen for some exceptional values of $a$ lying outside the
interval $0<|a|\ll1$ considered above.
Possible bifurcations of cuspidal edges should be similar to bifurcations of
cusps occuring in generic families of sections by parallel planes of
a swallow-tail (a rise/destruction of a pair of cuspidal edges),
a purse (a bifurcation of one cuspidal edge), and
a pyramid (a bifurcation of three cuspidal edges).
\end{note}

Of course, except for the singularities listed above, transversal intersections of
different branches of caustics are also possible:
\\ {\bf 11, 12.} ($A_2+A_3^\pm$, $A_{e,2}^\pm+A_3^\pm$)
In a neighbourhood of a caustic value $\overline{\lam}$ corresponding
to a degeneracy of an {\it additional} critical point of degeneracy type
$A_3^\pm$ and another {\it basic or additional} critical point of degeneracy
type $A_2$ or resp.\ $A_{e,2}^\pm$, the set of caustic values is the union of
a surface from case~3 ($A_3^\pm$) and a plane (from case~2 or resp.~1)
which intersects the exceptional curve of type $A_3^\pm$ of the first surface
transversally.
\\ {\bf 13, 14.} ($A_2+A_{e,3}^\pm$, $A_{e,2}^\pm+A_{e,3}^\pm$)
In a neighbourhood of a caustic value $\overline{\lam}$ corresponding to
a degeneracy of a {\it basic} critical point of degeneracy type $A_{e,3}^\pm$
and another {\it basic or additional} critical point of degeneracy type $A_2$
or resp.\ $A_{e,2}^\pm$, the set of caustic values is the union of a surface
from the case~4 ($A_{e,3}^\pm$) and a plane (from the case~2 or resp.~1) which
intersects the exceptional curve of type $A_{e,3}^\pm$ of the first surface
transversally.
\\ {\bf 15, 16, 17, 18.} ($A_2+A_2+A_2$, $A_2+A_2+A_{e,2}^\pm$,
$A_2+A_{e,2}^\pm+A_{e,2}^\pm$, $A_{e,2}^\pm+A_{e,2}^\pm+A_{e,2}^\pm$)
In a neighbourhood of a caustic value $\overline{\lam}$ corresponding
to a degeneracy of three different {\it basic or additional} critical points
of degeneracy types $A_2$ or $A_{e,2}^\pm$, the set of caustic values is the
union of three planes which intersect each other transversally at the point
$\overline{\lam}$.

\section{A determinant of singularity types} \label {sec:def}

In practice, the following questions often arise:
whether the germ of a given function $f$ has a certain singularity type, and
which normal form the function has near the singular point.
For singularities of codimension 0 (i.e.\ Morse functions and even Morse
functions), answers are given by the Morse lemma for arbitrary functions and
even functions, see Theorem~\ref {Mors}.

Below, we answer these questions for the singularities which appear in typical
$l-$parameter families with $l=1,2,3,4$ parameters for arbitrary functions and
even functions, see Statements~\ref {st:def},~\ref {st:defD4},~\ref {st:defD},
and resp.~\ref {st:defe},~\ref {st:defXe5},~\ref {st:defXe}.

\begin{statement} \label {st:def}
{\bf (A singularity $A_\mu$)}
Let $0$ be a critical point of a smooth function $f:\RR^\nu\to\RR$.
Suppose that the kernel of the second differential of $f$ at the origin
coincides with the coordinate axis $Ok_1$ (in particular, $0$ is a critical
point of co-rank~$1$). Then there exists a sequence of real numbers
$a_\mu$, $\mu=2,3,\dots$, depending on partial derivatives of $f$ at the
origin, and possessing the following properties:

1. $a_\mu=\frac{\partial^{\mu+1}f(0)}{\partial k_1^{\mu+1}}+\dots$, where
the additional terms form a polynomial in the values of the partial derivatives
$\frac{\partial^{i_1+\dots+i_\nu}f(0)}{\partial k_1^{i_1}\dots k_\nu^{i_\nu}}$,
$3\le i_1+\dots+i_\nu\le\mu$ (except the derivatives with respect to the
variable $k_1$ only) and the components of the inverse of the matrix
$\|\frac{\partial^2f(0)}{\partial k_i\partial k_j}\|$, $i,j=2,\dots,\nu$.
This polynomial has integer coefficients and a vanishing free term.
The numbers $a_\mu$ are given in Table~$(\ref {tab:def})$ for $\mu=2,3,4,5$.

2. The germ of $f$ at the origin has a singularity of type $A_\mu$
if and only if
 $$
 \refstepcounter{theorem} \label{critAmu}
a_2=\dots=a_{\mu-1}=0, \quad a_\mu\ne0.
 \eqno (\thetheorem)
 $$

Moreover, if~$(\ref {critAmu})$ is fulfilled then $f$ reduces to the form
$f=a+a_\mu\tilde k_1^{\mu+1}+Q(\tilde k_2,\dots,\tilde k_\nu)$
in some neighbourhood of the origin by means of a regular change of variables
$k\to\tilde k$ leaving the origin fixed.
Here $a$ is a constant, $Q$ is a nondegenerate quadratic form in $\nu-1$
variables.

In particular, $f\in A^+_\mu$ if $a_\mu>0$; $f\in A^-_\mu$ if $a_\mu<0$.
\end{statement}

\begin{statement} \label {st:defe}
{\bf (An even singularity $A_{e,\mu}$)}
Let $f:\RR^\nu\to\RR$ be a smooth even function defined in a neighbourhood
of the point $0$ in $\RR^\nu$.
Suppose that the kernel of the second differential of $f$ at the origin
coincides with the coordinate axis $Ok_1$ (in particular, $0$ is a critical
point of co-rank~$1$). Then there exists a sequence of real numbers
$a_{e,\mu}$, $\mu=2,3,\dots$, depending on partial derivatives of $f$ at the
origin, and possessing the following properties:

1. $a_{e,\mu}=a_{2\mu-1}=\frac{\partial^{2\mu}f(0)}{\partial k_1^{2\mu}}+\dots$
where $a_{2\mu-1}$ is assigned according to Statement~$\ref {st:def}$ to the
even function $f$, which is considered as a usual function.
The numbers $a_{e,\mu}$ are given in Table~$(\ref {tab:def})$ for $\mu=2,3,4,5$.

2. The even germ of $f$ at the origin has an even singularity of type
$A_{e,\mu}$ if and only if
 $$
 \refstepcounter{theorem} \label{critAemu}
a_{e,2}=\dots=a_{e,\mu-1}=0, \quad a_{e,\mu}\ne0.
 \eqno (\thetheorem)
 $$

Moreover, if~$(\ref {critAemu})$ is fulfilled then $f$ reduces to the form
$f=a+a_{e,\mu}\tilde k_1^{2\mu}+Q(\tilde k_2,\dots,\tilde k_\nu)$
in some neighbourhood of the origin by means of a regular odd change of variables
$k\to\tilde k$.
Here $a$ and $Q$ are the same as in Statement~$\ref {st:def}$.

In particular, $f\in A^+_{e,\mu}$ if $a_{e,\mu}>0$;
$f\in A^-_{e,\mu}$ if $a_\mu<0$.
\end{statement}

 $$
\refstepcounter{theorem}\label{tab:def}
\begin{tabular}{|c|l|l|}
\hline
$\mu$ & $a_\mu$ & $a_{e,\mu}$ \\
\hline
$2$ & $f_{y^3}$
    & $f_{y^4}$
\\
$3$ & $f_{y^4}-3 f_{z^2}^{-1} [f_{y^2z}]^2$
    & $f_{y^6}-10 f_{z^2}^{-1} [f_{y^3z}]^2$
\\
$4$ & $f_{y^5}- f_{y^2z} f_{z^2}^{-1}
      \left(10f_{y^3z}-15 f_{yz^2} f_{z^2}^{-1} f_{y^2z}\right)$
    & $f_{y^8}- 56 f_{y^3z} f_{z^2}^{-1}
      \left(f_{y^5z}- 5 f_{y^2z^2} f_{z^2}^{-1} f_{y^3z}\right)$
\\
$5$ & $f_{y^6} -15 f_{y^2z} f_{z^2}^{-1}
       (f_{y^4z} - 4 f_{yz^2} f_{z^2}^{-1} f_{y^3z}$
    & $f_{y^{10}} -120 f_{y^3z} f_{z^2}^{-1}
       (f_{y^7z} - 21 f_{y^2z^2} f_{z^2}^{-1} f_{y^5z}$
\\
& $\phantom{f_{y^6}}
        + 6 [f_{yz^2} f_{z^2}^{-1}]^2 f_{y^2z}  )
     -10 f_{z^2}^{-1} [f_{y^3z}]^2$
    & $\phantom{f_{y^10}} +105 [f_{y^2z^2} f_{z^2}^{-1}]^2 f_{y^3z} )
                          -126 f_{z^2}^{-1} [f_{y^5z}]^2$
\\
& $\phantom{f_{y^6}}
     +45 f_{y^2z^2} [f_{z^2}^{-1} f_{y^2z}]^2
     -15 f_{z^3} [f_{z^2}^{-1} f_{y^2z}]^3$
& $\phantom{f_{y^10}} +2100 f_{y^4z^2} [f_{z^2}^{-1} f_{y^3z}]^2
                      -2620 f_{yz^3} [f_{z^2}^{-1} f_{y^3z}]^3$
\\
 \hline
\end{tabular}
 \eqno (\thetheorem)
 $$
{\bf Comment to Table~(\ref {tab:def}).}
In this table, the variables are denoted by
$y=k_1$ and $z=(k_2,\dots,k_\nu)$, moreover
$f_{y^az^b}:=\frac{\partial^{a+b}f(0)}{\partial y^a\partial z^b}$.
The formulae in the table are written for $\nu=2$, but they are easily
transformed to the corresponding formulae for any $\nu\ge1$. For example:
 $$
\refstepcounter{theorem}\label{note:tab}
a_\mu = \frac{\partial^{\mu+1}f(0)}{\partial k_1^{\mu+1}}
\quad \mbox{for $\nu=1$}, \qquad
a_3 = \frac{\partial^4 f(0)}{\partial k_1^4} - 3 \sum_{i=2}^\nu\sum_{j=2}^\nu
\frac{\partial^3 f(0)}{\partial k_1^2 k_i} A_{ij}
\frac{\partial^3 f(0)}{\partial k_1^2 k_j} \quad \mbox{for $\nu\ge2$}
 \eqno (\thetheorem)
 $$
where $A_{ij}$, $i,j=2,\dots,\nu$ denote the components of the inverse of the
matrix
$\|\frac{\partial^2f(0)}{\partial k_i\partial k_j}\|$, $i,j=2,\dots,\nu$.

We will prove Statement~\ref {st:def} using the technique of the proof of
the lemma about classifying finite-multiple singularities of
co-rank~$1$~\cite [v.~1,~9.6]{Arnold}.
Due to the parametric Morse lemma for arbitrary smooth functions (analogous to
Theorem~\ref {th:paramMors}), in a small neighbourhood of the origin,
there exists a smooth change of variables
$(k_1,k_2,\dots,k_\nu)\to (k_1,\tilde k_2,\dots,\tilde k_\nu)$
leaving the origin fixed and reducing $f$ to the form
$f=a+\psi(k_1)+Q(\tilde k_2,\dots,\tilde k_\nu)$, where $a$ is a constant,
$Q$ is a nondegenerate quadratic form, and $\psi(k_1)=O(|k_1|^3)$ as $k_1\to0$.
It is obvious that the germ at the origin of a function $\psi$ in one variable
has a singularity of type $A_\mu$ if and only if the order of the function
$\psi$ at the origin is $\mu+1$, i.e.\
$\psi'(0)=\psi''(0)=\dots=\psi^{(\mu)}(0)=0$ and $\psi^{(\mu+1)}(0)\ne0$.

Therefore, it remains to express the numbers $a_\mu:=\psi^{(\mu+1)}(0)$ in terms
of the values of the partial derivatives of $f$ at the origin.
Observe that $\psi(k_1)=f(k_1,k_2(k_1),\dots,k_\nu(k_1))$ where the functions
$k_i=k_i(k_1)$, $i=2,\dots,\nu$ define the coordinate axis
$Ok_1$ of the coordinate system $k_1,\tilde k_2,\dots,\tilde k_\nu$ and,
hence, are implicitly defined by the system of equations
$\frac{\partial f}{\partial k_i}(k_1,k_2(k_1),\dots,k_\nu(k_1))=0$,
$i=2,\dots,\nu$.
One sets the values at the origin, of the derivatives of order $\le\mu-1$
of the left-hand sides of these equations, to zero.
This will give expressions for $k_i'(0),\dots,k_i^{(\mu-1)}(0)$ in terms of
the partial derivatives of order $\le\mu$ of $f$ at the origin.
By substituting these expressions into the derivative of order $\mu+1$ at the
origin of the composite function $\psi$, we obtain
$a_\mu=\psi^{(\mu+1)}(0)=\frac{\partial^{\mu+1}f(0)}{\partial k_1^{\mu+1}}+\dots$
where <<$\dots$>> is expressed by means of the partial derivatives
of order $\le\mu$ of $f$ at the origin.

Statement~\ref {st:defe} is similarly proved based on
Theorem~\ref {th:paramMors}.

\begin {statement} \label {st:defD4}
{\bf (A singularity $D_4$)}
Let $0$ be a critical point of a smooth function $f:\RR^\nu\to\RR$.
Suppose that the kernel of the second differential of $f$ at the origin
coincides with the coordinate plane $Ok_1k_2$ (in particular, $0$ is a critical
point of co-rank~$2$). The germ of $f$ at the origin has a singularity of type
$D_4$ if and only if $\Delta_3\ne0$, where $\Delta_3$ is the discriminant
 $$
\Delta_3 = 4 (B^3D+AC^3) + 27A^2D^2 - B^2C^2 - 18ABCD,
 $$
 $$
A=\frac{1}{3!}\frac{\partial^3f(0)}{\partial k_1^3}, \
B=\frac{1}{2!}\frac{\partial^3f(0)}{\partial k_1^2\partial k_2}, \
C=\frac{1}{2!}\frac{\partial^3f(0)}{\partial k_1\partial k_2^2}, \
D=\frac{1}{3!}\frac{\partial^3f(0)}{\partial k_2^3}.
 $$

Moreover, if $\Delta_3\ne0$ then $f$ reduces to the form
$f=d+\tilde k_1^2\tilde k_2+\Delta_3\tilde k_2^3+Q(\tilde k_3,\dots,\tilde k_\nu)$
in some neighbourhood of the origin by means of a regular change of variables
$k\to\tilde k$ leaving the point $0$ fixed. Here $d$ is a constant,
$Q$ is a nondegenerate quadratic form in $\nu-2$ variables.

In particular, $f\in D^+_4$ if $\Delta_3>0$; $f\in D^-_4$ if $\Delta_3<0$.
 \end{statement}

\begin{statement} \label {st:defD}
{\bf (A singularity $D_\mu$)}
Let $0$ be a critical point of a smooth function $f:\RR^\nu\to\RR$.
Suppose that the kernel of the second differential of $f$ at the origin
coincides with the coordinate plane $Ok_1k_2$ (in particular, $0$ is a critical
point of co-rank~$2$). Suppose that
$\frac{\partial^3f(0)}{\partial k_1^3}=\frac{\partial^3f(0)}{\partial k_1k_2^2}=0$
and $\frac{\partial^3f(0)}{\partial k_1^2\partial k_2}=1$ (this can be
achieved by means of a suitable linear change of the variables $k_1,k_2$
provided that the cubic part of the Taylor series at the origin of the function
$f(k_1,k_2,0,\dots,0)$ does not vanish and is not a perfect cube).
Then there exists a sequence of real numbers $d_\mu$, $\mu=4,5,\dots$,
depending on partial derivatives of $f$ at the origin, and possessing
the following properties:

1. $d_\mu=
\frac{\partial^{\mu-1}f(0)}{\partial k_2^{\mu-1}}+\dots$, where
the additional terms form a polynomial in the values of the partial derivatives
$\frac{\partial^{i_1+\dots+i_\nu}f(0)}{\partial k_1^{i_1}\dots k_\nu^{i_\nu}}$,
$3\le i_1+\dots+i_\nu\le\mu-2$
(except the derivatives with respect to the variable $k_2$ only)
and the components of the inverse of the matrix
$\|\frac{\partial^2f(0)}{\partial k_i\partial k_j}\|$, $i,j=3,\dots,\nu$.
This polynomial has rational coefficients and a vanishing free term.
The numbers $d_\mu$ are given in Table~$(\ref {tab:defD})$ for $\mu=4,5,6,7,8$.

2. For $\mu\ge4$, the germ of $f$ at the origin has a singularity of type
$D_\mu^\pm$ if and only if
 $$
 \refstepcounter{theorem} \label{critDmu}
d_4=\dots=d_{\mu-1}=0, \quad d_\mu\ne0.
 \eqno (\thetheorem)
 $$

Moreover, if~$(\ref {critDmu})$ holds then $f$ reduces to the form
$f=d+\tilde k_1^2\tilde k_2+d_\mu\tilde k_2^{\mu-1}+Q(\tilde k_3,\dots,\tilde k_\nu)$
in some neighbourhood of the origin by means of a regular change of variables
$k\to\tilde k$ leaving the origin fixed.
Here $d$ and $Q$ are similar to those in Statement~$\ref {st:defD4}$.

In particular, $f\in D^+_\mu$ if $d_\mu>0$; $f\in D^-_\mu$ if $d_\mu<0$.
 \end{statement}

\begin {statement} \label {st:defXe5}
{\bf (An even singularity $X_{e,5}$)}
Suppose that the kernel of the second differential at the origin of a smooth
even function $f:\RR^\nu\to\RR$ coincides with the coordinate plane $Ok_1k_2$
(in particular, $0$ is a critical point of co-rank~$2$).
Then the even germ of $f$ at the origin has an even singularity of class
$X_{e,5}^{\pm\pm}$ if and only if $\Delta_4\ne0$, where $\Delta_4$ is the
discriminant
 $$
\Delta_4 =
-4(AC^3D^2+B^3D^3+B^2C^3E)-27(A^2D^4+B^4E^2)+B^2C^2D^2+18(ABCD^3+B^3CDE)
 $$
 $$
+144(A^2CD^2E+AB^2CE^2)-6AB^2D^2E-80ABC^2DE-192A^2BDE^2+16AE(C^2-4AE)^2,
 $$
 $$
A=\frac{1}{4!}\frac{\partial^4f(0)}{\partial k_1^4}, \
B=\frac{1}{3!}\frac{\partial^4f(0)}{\partial k_1^3\partial k_2}, \
C=\frac{1}{4}\frac{\partial^3f(0)}{\partial k_1^2\partial k_2^2}, \
D=\frac{1}{3!}\frac{\partial^3f(0)}{\partial k_1\partial k_2^3}, \
E=\frac{1}{4!}\frac{\partial^4f(0)}{\partial k_2^4}.
 $$

Moreover, $f\in X^{++}_{e,5} \cup X^{--}_{e,5}$ if $\Delta_4>0$;
$f\in X^{+-}_{e,5}$ if $\Delta_4<0$.
 \end{statement}

\begin{note} \label{note:discrim} \rm
Suppose that, under hypothesis of Statement~\ref {st:defXe5},
$\Delta_4\ne0$. This is equivalent to the fact that the polynomial
$Ax^4+Bx^3+Cx^2+Dx+E$ does not have multiple roots~\cite{vanderWaerden}.
Therefore, by means of a linear change of variables, the corresponding
homogeneous polynomial of the fourth degree in two variables reduces to the
form $Ax^4+Cx^2y^2+Ey^4$. In these variables, we have $B=D=0$ and
$\Delta_4 = 16 AE(C^2-4AE)^2$. Consequently, $A\ne0$, $E\ne0$, and
$C^2\ne 4AE$. Therefore one can assume that $B=D=0$, $|A|=|E|=1$, and
that $|C|\ne2$ as soon as $AE=1$.
\end{note}

\begin{statement} \label {st:defXe}
{\bf (An even singularity $X_{e,\mu}$)}
Suppose that the kernel of the second differential at the origin of a smooth
even function $f:\RR^\nu\to\RR$ coincides with the coordinate plane $Ok_1k_2$
(in particular, $0$ is a critical point of co-rank~$2$).
Suppose that
 $$
\frac{\partial^4f(0)}{\partial k_1^3\partial k_2}
=\frac{\partial^4f(0)}{\partial k_1\partial k_2^3}=0, \quad
\frac{\partial^4f(0)}{\partial k_1^4}=24\varepsilon, \quad
\frac{\partial^4f(0)}{\partial k_1^2\partial k_2^2}=4\eta, \quad
\frac{\partial^4f(0)}{\partial k_2^4}\ne 6\varepsilon  \quad
\mbox{where} \quad \varepsilon,\eta=\pm1
 $$
(this can be achieved by means of a suitable linear change of the variables
$k_1,k_2$ provided that the homogeneous polynomial of the $4$-th degree
$Ak_1^4+Bk_1^3k_2+Ck_1^2k_2^2+Dk_1k_2^3+Ek_2^4$ of the Taylor series at the
origin of the function $f(k_1,k_2,0,\dots,0)$ does not vanish, is not a
perfect square, and reduces
either to the form~$1$ from~$\ref{subsec:4forms}$ with $a\ne0$,
or to the form~$2$ from~$\ref{subsec:4forms}$).
Then there exists a sequence of real numbers
$x_{e,\mu}=x_{e,\mu}^{\varepsilon,\eta}$, $\mu=5,6,\dots$
depending on partial derivatives of $f$ at the
origin and possessing the following properties:

1. $x_{e,\mu}^{\varepsilon,\eta}=
\frac{\partial^{2\mu-6}f(0)}{\partial k_2^{2\mu-6}}+\dots$,
where
the additional terms form a polynomial in the values of the partial derivatives
$\frac{\partial^{i_1+\dots+i_\nu}f(0)}{\partial k_1^{i_1}\dots k_\nu^{i_\nu}}$,
$4\le i_1+\dots+i_\nu\le2\mu-8$ (except the derivatives with respect to the
variable $k_2$ only) and the components of the inverse of the matrix
$\|\frac{\partial^2f(0)}{\partial k_i\partial k_j}\|$, $i,j=3,\dots,\nu$.
This polynomial has rational coefficients and a vanishing free term.
The numbers $x_{e,\mu}=x_{e,\mu}^{\varepsilon,\eta}$ are given in
Table~$(\ref {tab:defD})$ for $\mu=5,6,7,8$.

2. For $\mu\ge5$, the even germ of $f$ at the origin has a singularity of class
$X_{e,\mu}^{\pm\pm}$ if and only if
 $$
 \refstepcounter{theorem} \label{critXemu}
x_{e,5}^{\varepsilon,\eta}=\dots=x_{e,\mu-1}^{\varepsilon,\eta}=0, \quad
x_{e,\mu}^{\varepsilon,\eta}\ne0 .
 \eqno (\thetheorem)
 $$

Moreover, if~$(\ref {critXemu})$ holds then $f$ reduces to the form
$f=d+\varepsilon\tilde k_1^4+\eta\tilde k_1^2\tilde k_2^2
+ \frac{x_{e,\mu}^{\varepsilon,\eta}}{(2\mu-6)!} \tilde k_2^{2\mu-6}
+ Q(\tilde k_3,\dots,\tilde k_\nu)$
in some neighbourhood of the origin by means of a regular odd change of
variables $k\to\tilde k$.
Here $d$ and $Q$ are similar to those in Statement~$\ref {st:defD4}$.

In particular,
$\frac{x_{e,\mu}^{\varepsilon,\eta}}{(2\mu-6)!}$ is the even module
of the singularity $X_{e,\mu}^{\varepsilon,\eta}$ if $\mu\ge6$.
 \end{statement}

 $$
\refstepcounter{theorem}
\label{tab:defD}
\begin{tabular}{|c|l|l|l|}
\hline
$\mu$ & $d_\mu$ & $x_{e,\mu}=x_{e,\mu}^{\varepsilon,\eta}$ \\
\hline
$4$ & $a_{0,3}=f_{y^3}$ & $-$
\\
$5$ &
 $a_{0,4}=f_{y^4}-3 f_{z^2}^{-1} [f_{y^2z}]^2$ & $a_{0,4}=f_{y^4}$
\\
$6$ &
 $a_{0,5}-\frac{5}{3}a_{1,3}^2$ & $a_{0,6}=f_{y^6}-10 f_{z^2}^{-1} [f_{y^3z}]^2$
\\
$7$ & $a_{0,6}-5a_{1,3}a_{1,4}+5a_{1,3}^2a_{2,2}$
    & $a_{0,8}-\frac{7\eta}{10}a_{1,5}^2$
\phantom{$I^{I^{I^I}}$}\!\!\!\!\!\!\!\!\!\!
\\
$8$ & $a_{0,7}-7a_{1,3}a_{1,5}+\frac{35}{3}a_{1,3}^2a_{2,3}-
       \frac{35}{9}a_{1,3}^3a_{3,1} $
    & $a_{0,10}-3\eta\, a_{1,5}a_{1,7}+\frac{21}{16}a_{1,5}^2a_{2,4}$
\phantom{$I^{I^{I^I}}$}\!\!\!\!\!\!\!\!\!\!
\\
    & $\phantom{a_{0,7}} -\frac{35}{8} (a_{1,4}-2a_{1,3}a_{2,2})^2$
    &
\\
 \hline
\end{tabular}
 \eqno (\thetheorem)
 $$
{\bf Comment to Table~(\ref {tab:defD}).}
In this table, the variables are denoted by
$x=k_1$, $y=k_2$, and $z=(k_3,\dots,k_\nu)$, moreover
$f_{x^ay^bz^c}:=\frac{\partial^{a+b+c}f(0)}{\partial x^a\partial y^b\partial z^c}$.
Furthermore, one denotes $a_{0,j}=a_{j-1}$ where $a_\mu$ are the real numbers
assigned according to Statement~$\ref {st:def}$ to the function
$f(0,k_2,\dots,k_\nu)$, see Table~(\ref {tab:def}). Other real numbers
$a_{i,j}=\frac{\partial f^{i+j}(0)}{\partial x^i\partial y^j}+\dots$
in Table~(\ref {tab:defD}) are similar to the numbers $a_{0,j}$ and
are defined by the formulae
 $$
a_{1,3}=f_{xy^3}-3 f_{xyz} f_{z^2}^{-1} f_{y^2z},\quad
a_{2,2}=f_{x^2y^2}-2 f_{z^2}^{-1} [f_{xyz}]^2 - f_{x^2z} f_{z^2}^{-1} f_{y^2z},
 $$
 $$
a_{1,4}=f_{xy^4}-6 f_{xy^2z} f_{z^2}^{-1} f_{y^2z}
-4 f_{y^3z} f_{z^2}^{-1} f_{xyz}
+3 f_{xz^2} [f_{z^2}^{-1} f_{y^2z}]^2
+12 f_{y^2z} f_{z^2}^{-1} f_{yz^2} f_{z^2}^{-1} f_{xyz},
 $$
 $$
a_{2,3}=f_{x^2y^3}-f_{x^2z} f_{z^2}^{-1} f_{y^3z}
-3 (f_{x^2yz} - f_{x^2z} f_{z^2}^{-1} f_{yz^2} -2 f_{xyz} f_{z^2}^{-1} f_{xz^2})
f_{z^2}^{-1} f_{y^2z}
 $$
 $$
-6 (f_{xyz} f_{z^2}^{-1} f_{xy^2z} - f_{yz^2} [f_{z^2}^{-1} f_{xyz}]^2) ,
 $$
 $$
a_{1,5}=f_{xy^5} -10 f_{y^3z} f_{z^2}^{-1} f_{xy^2z}
-30 f_{xz^2} f_{z^2}^{-1} f_{yz^2} [f_{z^2}^{-1}  f_{y^2z}]^2
-60 f_{xyz} f_{z^2}^{-1} [f_{yz^2} f_{z^2}^{-1}]^2 f_{y^2z}
-5 f_{y^4z} f_{z^2}^{-1} f_{xyz}
 $$
 $$
-10 f_{xy^3z} f_{z^2}^{-1} f_{y^2z}
-15 f_{z^3} [f_{z^2}^{-1} f_{xyz}]^3
+15 f_{xyz^2} [f_{z^2}^{-1} f_{y^2z}]^2
+30 f_{xyz} f_{z^2}^{-1} f_{y^2z^2} f_{z^2}^{-1} f_{y^2z}
 $$
 $$
+10 f_{y^2z} f_{z^2}^{-1} f_{xz^2} f_{z^2}^{-1} f_{y^3z}
+20 f_{xyz} f_{z^2}^{-1} f_{yz^2} f_{z^2}^{-1} f_{y^3z}
+30 f_{y^2z} f_{z^2}^{-1} f_{yz^2} f_{z^2}^{-1} f_{xy^2z}.
 $$
Moreover, if $f$ is an even function then
 $$
a_{2,4}=f_{x^2y^4} - 4 f_{x^2yz} f_{z^2}^{-1} f_{y^3z}
                   - 6 f_{z^2}^{-1} [f_{xy^2z}]^2 ,
 $$
 $$
a_{1,7}=f_{xy^7} - 35 f_{xy^4z} f_{z^2}^{-1} f_{y^3z}
                 - 21 f_{xy^2z} f_{z^2}^{-1} ( f_{y^5z}
                      - 10 [f_{y^2z^2} f_{z^2}^{-1}]^2 f_{y^3z})
                 + 70 f_{xyz^2} [f_{z^2}^{-1}  f_{y^3z}]^2 .
 $$
One can obtain $a_{j,i}$ from $a_{i,j}$ by replacing
all partial derivatives with respect to $x$
by partial derivatives with respect to $y$, and vica-versa.
The above formulae are written for $\nu=3$, but they are easily
transformed to the corresponding formulae for any $\nu\ge2$,
see~(\ref {note:tab}).
In the right column of Table~(\ref {tab:defD}), the function $f$ is supposed to
be even. This leads to a simplification of the formulae for the numbers
$a_{i,j}$.

One proves Statements~\ref {st:defD4} and~\ref {st:defD} similarly to the
proof of Statement~\ref {st:def}, using the Newton ruler
method~\cite[v.~1,~12.6]{Arnold}, see also~\ref{subsec:4forms}.
Namely, due to the parametric Morse lemma for arbitrary smooth functions
(analogous to Theorem~\ref {th:paramMors}), in a small neighbourhood of the
origin, there exists a smooth change of variables
$(k_1,k_2,k_3,\dots,k_\nu)\to (k_1,k_2,\tilde k_3,\dots,\tilde k_\nu)$
leaving the origin fixed and reducing $f$ to the form
$f=a+\psi(k_1,k_2)+Q(\tilde k_3,\dots,\tilde k_\nu)$, where $a$ is a constant,
$Q$ is a nondegenerate quadratic form, and $\psi(k_1,k_2)=o(k_1^2+k_2^2)$ as
$(k_1,k_2)\to0$.
Let $P(k_1,k_2)=Ak_1^3+Bk_1^2k_2+Ck_1k_2^2+Dk_2^3$ be the Taylor polynomial of
degree 3 of the function $\psi$ at zero.
It follows from the Newton ruler method that the germ at the origin of the
function $\psi$ in two variables has a singularity of type $D_4$ if and only if
three lines, which are defined by the linear factors of the decomposition of
$P$, are pairwise different, compare~\ref {subsec:4forms}.
This is equivalent to the fact that the discriminant $\Delta_3$ of the
polynomial $Ax^3+Bx^2+Cx+D$ does not vanish~\cite{vanderWaerden}.
To see the validity of the formulae for the coefficients $A,B,C,D$ in
Statement~\ref {st:defD4}, one uses the fact that the definition of the
<<cubical>> part (terms of degree $\le3$) of a function at a critical
point is well defined~\cite{Arnold}, and the fact
that the coordinate plane $Ok_1k_2$ coincides with the kernel of the second
differential of $f$ at zero, both for the initial variables and for the new
variables. Statement~\ref {st:defD4} is thus proved.

In order to prove Statement~\ref {st:defD}, we observe that $A=C=0$ and
$B=\frac{1}{2}$. If $D\ne0$ then $\Delta_3\ne0$ and, by the above arguments,
$f$ has a singularity at the origin of type $D_4$. Now we assume that $D=0$ and
$\mu\ge5$.
Following the method of ruler turning, we perform the change of variables
$(k_1,k_2)\to(\tilde k_1,\tilde k_2)$ given by
$\tilde k_1=k_1+\lam_2k_2^2+\dots+\lam_{[\mu/2]-1}k_2^{[\mu/2]-1}$,
$\tilde k_2=k_2$, where the $\lam_j$ are determined implicitely (also uniquely
and independently of $\mu$) by the following condition: after this change of
variables the coefficients at the terms
$\tilde k_1\tilde k_2^j$, $3\le j\le[\mu/2]$, of the Taylor series of
$\psi$ in $\tilde k_1,\tilde k_2$, centred at the origin, vanish.
Due to the ruler method, the function $\psi$ (and, hence, the function $f$)
has a singularity at the origin of type $D_\mu$ if and only if the coefficient
at the term $\tilde k_2^j$ in the Taylor series of $\psi$ at zero
vanishes for any $j\le\mu-2$ and does not vanish for $j=\mu-1$.

Therefore, it remains to express the numbers
$d_\mu:=\frac{\partial^{\mu-1}\psi(0)}{\partial\tilde k_2^{\mu-1}}$ in terms
of the values of the partial derivatives of $f$ at the origin.
First we will express the numbers
$a_{q,p}:=\frac{\partial^{q+p}\psi(0)}{\partial k_1^q\partial k_2^p}=
\frac{\partial^{q+p}f(0)}{\partial k_1^q\partial k_2^p}+\dots$.
Computations, analogous to those for the numbers
$a_\mu=\psi^{(\mu+1)}(0)=\frac{\partial^{\mu+1}f(0)}{\partial k_1^{\mu+1}}+\dots$
(see the proof of Statement~\ref {st:def}), lead to the above formulas for
$a_{q,p}$, see the comment to Table~(\ref{tab:defD}).
By means of the method of ruler turning, one easily obtains an expression
for the coefficient $\lam_j$ of the change of variables
$(k_1,k_2)\to(\tilde k_1,\tilde k_2)$ in terms of the numbers $a_{q,p}$,
$q\ge1$, $2q+p\le j+3$. After performing this change of variables,
one obtains the required expression $d_\mu=a_{0,\mu-1}+\dots$ in terms
of the numbers $a_{q,p}$, $2q+p\le\mu-1$. This proves Statement~\ref {st:defD}.

Statements~\ref {st:defXe5} and~\ref {st:defXe} are similarly proved based on
Theorem~\ref {th:paramMors} and the results of Section~\ref {subsec:4forms},
see the cases~1 and~2.

Statements~\ref {st:def},~\ref {st:defe},~\ref {st:defD4},~\ref {st:defD},~\ref
{st:defXe5},~\ref {st:defXe}, and Theorem~\ref {th:cork}
demonstrate that singularities of type $A_\mu$, $D_\mu$, and even singularities
of type $A_{e,\mu}$ have codimension $c=\mu-1$, while even singularities of
class $X_{e,\mu}^{\pm\pm}$ have even codimension $c=\mu-2$. That is, these
singularities are non-removable for typical families with $l\ge c$ parameters
on the entire family.
Besides, it follows from these assertions that the only singularities which
appear in typical families with $l=1,2,3,4$ parameters are $A_1$, $A_2$, $A_3$,
$A_4$, $A_5$, $D_4$, and $D_5$ (the seven of Thom).
In additional, in the case of typical even families, the only even singularities
which appear in such families are
$A_{e,1}$, $A_{e,2}$, $A_{e,3}$, $A_{e,4}$, $A_{e,5}$, $X_{e,5}$, and $X_{e,6}$.

Here, by {\it typical} (even) families of functions, we understand (even)
families which form an open dense subspace in the space of all (even) families
of functions (in the convergence topology with a finite number
of derivatives on each compact set).

Actually, due to Tables~(\ref {obKl}) and~(\ref {evKl}), in the case of five
parameters, three more singularities $A_6$, $D_6$, and $E_6$, as well as
five more even singularities $A_{e,6}$, $X_{e,7}$, $Y_{e,3,3}$,
$\tilde Y_{e,3}$, and $Z_{e,7}$ appear.

\section{Conditions for versal deformations} \label {sec:cond:versal}

In this section, we describe conditions on a deformation of a germ
to be versal or, equivalently, stable, see~\ref {note:vozmTh}.

It is known that the versal property is fulfilled for {\it typical}
deformations with $l\ge\mu-1$ parameters, i.e.\ deformations which form an open
dense subspace in the space of all $l-$parameter deformations of a given germ
(in the convergence topology with a finite number of derivatives on each
compact set). However,
the deformations which appear in practice are often not arbitrary,
but have a special form. Therefore, in addition to the question about
the determination of the type of a singular point of a given function $f(k)$,
the following questions are also of interest for applications:
whether an $l-$parameter deformation $F(k,\lam)$ of the germ of this function
is versal, and which normal form this deformation has near the singular point.

In the case of singularities of codimension $0$ (Morse functions and
even Morse functions), answers are given by the parametric Morse lemma for
arbitrary functions and even functions. The following assertions give
conditions on $l-$parameter deformations with $l=1,2,3,4$ parameters to be
versal for <<usual>> singularities and even singularities.
As we have mentioned above, these assertions can be considered as analogues
of the parametric Morse lemma for degenerate singularities.

In the following assertions, one considers a family of functions $F=F(k,\lam)$
in a variable $k=(k_1,\dots,k_\nu)$ with a parameter
$\lam=(\lam_1,\dots,\lam_l)$. One denotes by
$\frac{\partial^{i_1+\dots+i_\nu+1}F(0,0)}
{\partial k_1^{i_1}\dots\partial k_\nu^{i_\nu}\partial\lam}$
the $l-$dimensional vector (i.e.\ the element of $\RR^l$) with components
$\frac{\partial^{i_1+\dots+i_\nu+1}F(0,0)}
{\partial k_1^{i_1}\dots\partial k_\nu^{i_\nu}\partial\lam_j}$, $1\le j\le l$.

\begin{statement} \label {st:versal}
{\bf (A deformation of a singularity $A_\mu$)}
Let, under the hypothesis of Statement~$\ref {st:def}$,
$F:\RR^\nu\times\RR^l\to \RR$ be a smooth $l-$parameter
deformation of the germ of $f$ at the origin, i.e.\ $f=F(\cdot,0)$.
Then there exists a sequence of
vectors $\v_\mu\in\RR^l$, $\mu=2,3,\dots$, depending on partial derivatives
of $f$ at the origin, and possessing the following properties:

1. $\v_\mu=\frac{\partial^{\mu}F(0,0)}{\partial k_1^{\mu-1}\partial\lam}+\dots$,
where
the additional terms form a linear combination of vectors
$\frac{\partial^{i_1+\dots+i_\nu+1}F(0,0)}
{\partial k_1^{i_1}\dots\partial k_\nu^{i_\nu}\partial\lam}$,
$i_1+\dots+i_\nu+1<\mu$, $(i_2,\dots,i_\nu)\ne(0,\dots,0)$, the coefficients of
which are polynomials in the values of the partial derivatives of order $\le\mu$
of $f$ at the origin (except the derivatives with respect to the
variable $k_1$ only) and the components of the inverse of the matrix
$\|\frac{\partial^2f(0)}{\partial k_i\partial k_j}\|$, $i,j=2,\dots,\nu$.
These polynomials have integer coefficients and vanishing free terms. The
vectors $\v_\mu$ are given in Table~$(\ref {tab:versal})$ for $\mu=2,3,4,5$.

2. Suppose that the conditions~$(\ref {critAmu})$ are fulfilled, i.e.\
the germ of $f$ at the origin has a singularity of type $A_\mu$.
The deformation $F$ is $R^+-$versal if and only if
the vectors
 $$
\v_2,\ \dots,\ \v_\mu
 $$
form a linearly independent system in $\RR^l$ (in particular, $l\ge\mu-1$).

Moreover, if~$(\ref {critAmu})$ holds and the above vectors form a linearly
independent system in $\RR^l$ then $F$ reduces to the form
$F=a(\tilde\lam)+a_\mu\tilde k_1^{\mu+1}+Q(\tilde k_2,\dots,\tilde k_\nu)+
\tilde\lam_1\tilde k_1+\dots+\tilde\lam_{\mu-1}\tilde k_1^{\mu-1}$
in some neighbourhood of the origin by means of a regular change of variables
$(k,\lam)\to(\tilde k,\tilde\lam)$ leaving the origin fixed and having the form
$\tilde k=\tilde k(k,\lam)$, $\tilde\lam=\tilde\lam(\lam)$.
Here $a$ is smooth function, and $Q$ is a nondegenerate quadratic form in
$\nu-1$ variables.
\end{statement}

\begin{statement} \label {st:versale}
{\bf (An even deformation of a singularity $A_{e,\mu}$)}
Let, under the hypothesis of Statement~$\ref {st:defe}$, $f(0)=0$ and let
$F:\RR^\nu\times\RR^l\to \RR$ be a smooth $l-$parameter even deformation
of the germ of $f$ at the origin. That is, $f=F(\cdot,0)$ and, for each $\lam$,
the function $F(\cdot,\lam)$ is even and has a vanishing value at the origin.
Then there exists a sequence of
vectors $\v_{e,\mu}\in\RR^l$, $\mu=2,3,\dots$, depending on partial derivatives
of $F$ at the origin, and possessing the following properties:

1. $\v_{e,\mu}=\v_{2\mu-1}
      =\frac{\partial^{2\mu-1}F(0,0)}{\partial k_1^{2\mu-2}\partial\lam}+\dots$,
where the vectors $\v_{2\mu-1}$ are assigned according to
Statement~$\ref {st:versal}$ to the even deformation $F$, which is considered
as a usual deformation.
The vectors $\v_{e,\mu}$ are given in Table~$(\ref {tab:versal})$ for
$\mu=2,3,4$.

2. Suppose that the conditions~$(\ref {critAemu})$ are fulfilled, i.e.\
the even germ of $f$ at the origin has a singularity of type $A_{e,\mu}$.
Then the even deformation $F$ is $R_O-$versal if and only if the vectors
 $$
\v_{e,2},\ \dots,\ \v_{e,\mu}
 $$
form a linearly independent system in $\RR^l$ (in particular, $l\ge\mu-1$).

Moreover, if~$(\ref {critAemu})$ holds and the above vectors form a linearly
independent system in $\RR^l$ then $F$ reduces to the form
$F=a_{e,\mu}\tilde k_1^{2\mu}+Q(\tilde k_2,\dots,\tilde k_\nu)+
\tilde\lam_1\tilde k_1^2+\dots+\tilde\lam_{\mu-1}\tilde k_1^{2(\mu-1)}$
in some neighbourhood of the origin by means of a regular change of variables
$(k,\lam)\to(\tilde k,\tilde\lam)$ leaving the origin fixed and having the form
$\tilde k=\tilde k(k,\lam)$, $\tilde\lam=\tilde\lam(\lam)$ with
$\tilde k(-k,\lam)=-\tilde k(k,\lam)$.
Here $Q$ is as in Statement~$\ref {st:versal}$.
\end{statement}

 $$
\refstepcounter{theorem}
\label{tab:versal}
\begin{tabular}{|c|l|l|}
\hline
$\mu$ & $\v_\mu$ & $\v_{e,\mu}$ \\
\hline
$2$ & $F_{y\lam}$
    & $F_{y^2\lam}$
\\
$3$ & $F_{y^2\lam}-f_{y^2z}f_{z^2}^{-1} F_{z\lam}$
    & $F_{y^4\lam}-4f_{y^3z}f_{z^2}^{-1} F_{yz\lam}$
\\
$4$ & $F_{y^3\lam}
    - \left( f_{y^3z}
        -3f_{y^2z} f_{z^2}^{-1} f_{yz^2}
      \right) f_{y^2}^{-1} F_{z\lam}$
    & $F_{y^6\lam}
      -6 \left( f_{y^5z} - 10 f_{y^3z} f_{z^2}^{-1} f_{y^2z^2}
         \right) f_{z^2}^{-1} F_{yz\lam}
$ \\ & $\phantom{F_{y^3\lam}}
    -3f_{y^2z} f_{z^2}^{-1} F_{yz\lam}$
     & $\phantom{F_{y^6\lam}}
      -20 f_{y^3z} f_{z^2}^{-1} F_{y^3z\lam}
      +4 F_{z^2\lam} [f_{z^2}^{-1} f_{y^3z}]^2$
\\
$5$ & $F_{y^4\lam}
      -4 \left( f_{y^3z} -3f_{y^2z} f_{z^2}^{-1} f_{yz^2}
       \right) f_{z^2}^{-1} F_{yz\lam}
$ & $F_{y^8\lam} + \dots$ \\ & $\phantom{F_{y^4\lam}}
      -6f_{y^2z}f_{z^2}^{-1} F_{y^2z\lam}
      +F_{z^2\lam} [f_{z^2}^{-1} f_{y^2z}]^2
$ & \\ & $
   -\left( f_{y^4z}
   -4f_{y^3z} f_{z^2}^{-1} f_{yz^2}
   +9f_{y^2z} [f_{z^2}^{-1} f_{yz^2}]^2
    \right.
$ & \\ & $\left.
   -6f_{y^2z} f_{z^2}^{-1} f_{y^2z^2}
   - f_{z^3} [f_{z^2}^{-1} f_{y^2z}]^2
   \right) f_{z^2}^{-1} F_{z\lam}
$ &
\\
 \hline
\end{tabular}
\eqno (\thetheorem)
 $$
{\bf Comment to Table~(\ref {tab:versal}).}
In this table, the variables and the parameters are denoted by
$y=k_1$, $z=(k_2,\dots,k_\nu)$, and $\lam=(\lam_1,\dots,\lam_l)$, resp.,
moreover $f_{y^az^b}:=\frac{\partial^{a+b}f(0)}{\partial y^a\partial z^b}$ and
$F_{y^az^b\lam}:=
 \frac{\partial^{a+b+1}F(0,0)}{\partial y^a\partial z^b\partial\lam}$.
The formulae in the table are written for $\nu=2$, $l=1$, but they are easily
transformed to the corresponding formulae for any $\nu,l\ge1$,
see~(\ref {note:tab}).

For a proof of Statement~\ref {st:versal}, let us perform, for any
sufficiently small value of the parameter $\lam$, a change of variables
similar to that in the proof of Statement~\ref {st:def}. In the new
variables $k_1,\tilde K_2,\dots,\tilde K_\nu,\lam$, we have
$F=A(\lam)+\Psi(k_1,\lam)+Q(\tilde K_2,\dots,\tilde K_\nu)$, where
$A$ and $\Psi$ are smooth functions, $A(0)=a$, $\Psi(\cdot,0)=\psi$,
$\tilde K=\tilde K(k,\lam)$. Since the function $\psi(k_1)$ has order $\mu+1$
at the origin, the germ of the function $k_1^\mu$ belongs to the Jacobian
ideal $I_{\nabla\psi}$. It follows that the germs of the functions
$k_1,\dots,k_1^{\mu-1}$ are generators of the local algebra $Q_{\nabla\psi}$
of the gradient map of
the function $\psi$ at the origin (considered as a vector space), and that
the differential operators
$\frac{\partial}{\partial k_1}|_0,\frac{\partial^2}{\partial k_1^2}|_0,\dots,
\frac{\partial^{\mu-1}}{\partial k_1^{\mu-1}}|_0$ are generators of
its dual space $Q_{\nabla\psi}^*$.
It follows from this, from the conditions on a deformation to be
infinitesimally $R^+-$versal, and from the versality theorem
(see the theorems from~\cite [v.~1,~8.2 and~8.3]{Arnold}), that
a deformation $\Psi$ of a function $\psi$ is $R^+-$versal if and only if the
vectors
$\v_i:=\frac{\partial^i\Psi(0,0)}{\partial k_1^{i-1}\partial\lam}\in\RR^l$,
$i=2,\dots,\mu$ are linearly independent.

In order to compute the vectors $\v_\mu$, observe that
$\Psi(k_1,\lam)=F(k_1,K_2(k_1,\lam),\dots,K_\nu(k_1,\lam);\lam)$ where
$K_i=K_i(k_1,\lam)$, $i=2,\dots,\nu$ are the functions which define the
coordinate plane $Ok_1\lam$ of the coordinate system
$k_1,\tilde K_2,\dots,\tilde K_\nu,\lam$ and, therefore, are implicitly
defined by the system of equations
 $$
\frac{\partial F}{\partial k_i}(k_1,K_2(k_1,\lam),\dots,K_\nu(k_1,\lam);\lam)=0,
\quad i=2,\dots,\nu.
 $$
In particular, $K_i(k_1,0)=k_i(k_1)$, $i=2,\dots,\nu$, see the proof of
Statement~\ref {st:def}.
One now sets to zero the value of each composite partial derivative of bi-order
$(m,1)$, $m\le\mu-2$, at the origin with respect to $(k_1,\lam)$ on the
left-hand sides of these equations.
This will give expressions for composite partial derivatives at the origin
of bi-order $(m,1)$ with respect to $(k_1,\lam)$ of the functions
$K_i(k_1,\lam)$.
By substituting these expressions and the expressions for
$k_i'(0),\dots,k_i^{(\mu)}(0)$ into the composite partial derivative of order
$(\mu-1,1)$ of the composite function $\Psi$ at the origin with respect to
$(k_1,\lam)$, we obtain the required expression for $\v_\mu$.

Statement~\ref {st:versale} is similarly proved on the basis of
Theorem~\ref {verMO}.

\begin{statement} \label {st:versalD}
{\bf (A deformation of a singularity $D_\mu$)}
Let, under the hypothesis of Statement~$\ref {st:defD}$,
$F:\RR^\nu\times\RR^l\to \RR$ be a smooth $l-$parameter
deformation of the germ of $f$ at the origin, i.e.\ $f=F(\cdot,0)$.
Then there exists a sequence of
vectors $\w_\mu\in\RR^l$, $\mu=1,2,3,\dots$, depending on partial derivatives
of $F$ at the origin, and possessing the following properties:

1.
$\w_2=\frac{\partial^2F(0,0)}{\partial k_1\partial\lam}$ and
$\w_1=\frac{\partial^3F(0,0)}{\partial k_1^2\partial\lam}+\dots$
are equal to the
vectors $\v_2$ and $\v_3$, resp., which are assigned according to
Statement~$\ref {st:versal}$ to the deformation
$F(k_1,0,k_3,\dots,k_\nu;\lam)$ of the germ of $f(k_1,0,k_3,\dots,k_\nu)$
at the origin. The vectors
$\w_3=\frac{\partial^2F(0,0)}{\partial k_2\partial\lam}$ and
$\w_4=\frac{\partial^3F(0,0)}{\partial k_2^2\partial\lam}+\dots$
are similarly defined by means of the deformation
$F(0,k_2,k_3,\dots,k_\nu;\lam)$ of the germ of $f(0,k_2,k_3,\dots,k_\nu)$
at the origin. For $\mu\ge5$, one has
$\w_\mu=\frac{\partial^{\mu-1}F(0,0)}{\partial k_2^{\mu-2}\partial\lam}+\dots$,
where
the additional terms form a linear combination of vectors
$\frac{\partial^{i_1+\dots+i_\nu+1}F(0,0)}
{\partial k_1^{i_1}\dots\partial k_\nu^{i_\nu}\partial\lam}$,
$i_1+\dots+i_\nu+1<\mu-1$, $(i_1,i_3,\dots,i_\nu)\ne(0,\dots,0)$, the
coefficients of which are polynomials in the values of the partial derivatives
of order $\le\mu-1$ of $f$ at the origin (except the derivatives with respect to
the variable $k_2$ only) and the components of the inverse of the matrix
$\|\frac{\partial^2f(0)}{\partial k_i\partial k_j}\|$, $i,j=3,\dots,\nu$.
These polynomials have rational coefficients and vanishing free terms. The
vectors $\w_\mu$ are given in Table~$(\ref {tab:versalD})$ for
$\mu=1,2,3,4,5,6,7$.

2. Suppose that the conditions~$(\ref {critDmu})$ are fulfilled, i.e.\
the germ of $f$ at the origin has a singularity of type $D_\mu^\pm$, $\mu\ge4$.
The deformation $F$ is $R^+-$versal if and only if the vectors
 $$
\w_2,\ \dots,\ \w_{\mu-1},\ \w_\mu - d_\mu \w_1
 $$
form a linearly independent system in $\RR^l$ (in particular, $l\ge\mu-1$).

Moreover, if~$(\ref {critDmu})$ holds and the above vectors form a linearly
independent system in $\RR^l$ then $F$ reduces to the form
$F=d(\tilde\lam)+ \tilde k_1^2\tilde k_2
+ d_\mu\tilde k_2^{\mu-1}+Q(\tilde k_3,\dots,\tilde k_\nu)+
\tilde\lam_1\tilde k_1+\tilde\lam_2\tilde k_1^2+
\tilde\lam_3\tilde k_2+\tilde\lam_4\tilde k_2^2+\dots
 +\tilde\lam_{\mu-1}\tilde k_2^{\mu-3}$
in some neighbourhood of the origin by means of a regular change of variables
$(k,\lam)\to(\tilde k,\tilde\lam)$ leaving the origin fixed and having the form
$\tilde k=\tilde k(k,\lam)$, $\tilde\lam=\tilde\lam(\lam)$.
Here $d$ is a smooth function, and $Q$ is a nondegenerate quadratic form in
$\nu-2$ variables.
\end{statement}

\begin{statement} \label {st:versalXe}
{\bf (An even deformation of a singularity $X_{e,\mu}$)}
Let, under the hypothesis of Statement~$\ref {st:defXe}$, $f(0)=0$ and let
$F:\RR^\nu\times\RR^l\to \RR$ be a smooth $l-$parameter even deformation
of the germ of $f$ at the origin. That is, $f=F(\cdot,0)$ and, for each $\lam$,
the function $F(\cdot,\lam)$ is even and has a vanishing value at the origin.
Then there exists a sequence of
vectors $\x_{e,\mu}=\x_{e,\mu}^{\varepsilon,\eta}\in\RR^l$, $\mu=1,2,3,\dots$
depending on partial derivatives of $F$ at the origin and possessing the
following properties:

1.
$\x_{e,2}^{\varepsilon,\eta}
         =\frac{\partial^3F(0,0)}{\partial k_1^2\partial\lam}$,
$\x_{e,3}^{\varepsilon,\eta}
         =\frac{\partial^3F(0,0)}{\partial k_1\partial k_2\partial\lam}$,
$\x_{e,1}^{\varepsilon,\eta}
        = \varepsilon(\frac{\partial^5F(0,0)}{\partial k_1^4\partial\lam}+\dots)
-12\eta(\frac{\partial^5F(0,0)}{\partial k_1^2\partial k_2^2\partial\lam}+\dots)$,
see Table~$(\ref {tab:versalD})$.
For $\mu\ge4$, one has $\x_{e,\mu}^{\varepsilon,\eta}
 =\frac{\partial^{2\mu-5}F(0,0)}{\partial k_2^{2\mu-6}\partial\lam}+\dots$,
where
the additional terms form a linear combination of vectors
$\frac{\partial^{i_1+\dots+i_\nu+1}F(0,0)}
{\partial k_1^{i_1}\dots\partial k_\nu^{i_\nu}\partial\lam}$,
$i_1+\dots+i_\nu+1<2\mu-5$, $(i_1,i_3,\dots,i_\nu)\ne(0,\dots,0)$, the
coefficients of which are polynomials in the values of the partial derivatives
of order $\le2\mu-4$ of $f$ at the origin (except the derivatives with respect
to the variable $k_2$ only) and the components of the inverse of the matrix
$\|\frac{\partial^2f(0)}{\partial k_i\partial k_j}\|$, $i,j=3,\dots,\nu$.
These polynomials have rational coefficients and vanishing free terms. The
vectors $\x_{e,\mu}^{\varepsilon,\eta}$ are given in
Table~$(\ref {tab:versalD})$ for $\mu=1,2,3,4,5,6$.

2. Suppose that the conditions~$(\ref {critXemu})$ are fulfilled, i.e.\
the even germ of $f$ at the origin has a singularity of class
$X_{e,\mu}^{\varepsilon,\eta}$, $\mu\ge5$. The even deformation $F$ is
$R_O-$versal if and only if the vectors
 $$
\x_{e,2}^{\varepsilon,\eta},\ \dots,\ \x_{e,\mu-1}^{\varepsilon,\eta},\
\x_{e,\mu}^{\varepsilon,\eta} + \frac{\mu-3}{48}\,x_{e,\mu}^{\varepsilon,\eta} \,
\x_{e,1}^{\varepsilon,\eta}
 $$
form a linearly independent system in $\RR^l$ (in particular, $l\ge\mu-1$).

Moreover, if~$(\ref {critXemu})$ holds and the above vectors form a linearly
independent system in $\RR^l$, then $F$ reduces to the form
$F =
\varepsilon\tilde k_1^4+\eta\tilde k_1^2\tilde k_2^2
+ \frac{x_{e,\mu}^{\varepsilon,\eta}}{(2\mu-6)!} \tilde k_2^{2\mu-6}
+ Q(\tilde k_3,\dots,\tilde k_\nu)
+ \tilde\lam_1\tilde k_1^2 + \tilde\lam_2\tilde k_1\tilde k_2
+ \tilde\lam_3\tilde k_2^2 + \tilde\lam_4\tilde k_1^4
+ \tilde\lam_5\tilde k_2^4 + \dots + \tilde\lam_{\mu-1}\tilde k_2^{2\mu-8}$
in some neighbourhood of the origin by means of a regular change of variables
$(k,\lam)\to(\tilde k,\tilde\lam)$ leaving the origin fixed and having the form
$\tilde k=\tilde k(k,\lam)$, $\tilde\lam=\tilde\lam(\lam)$ with
$\tilde k(-k,\lam)=-\tilde k(k,\lam)$.
Here $Q$ is as in Statement~$\ref {st:versalD}$.
\end{statement}

 $$
\refstepcounter{theorem}
\label{tab:versalD}
\begin{tabular}{|c|l|l|}
\hline
$\mu$ &
$\w_\mu\phantom{I^{I^I}}$ & $\x_{e,\mu}=\x_{e,\mu}^{\varepsilon,\eta}$\\
\hline
$1$ & $\v_{2,0} = F_{x^2\lam}-f_{x^2z}f_{z^2}^{-1} F_{z\lam}$
    & $\varepsilon\v_{4,0}-12\eta\v_{2,2}$
\\
$2$ & $\v_{1,0} = F_{x\lam}$
    & $\v_{2,0} = F_{x^2\lam}$
\\
$3$ & $\v_{0,1} = F_{y\lam}$
    & $\v_{1,1} = F_{xy\lam}$
\\
$4$ & $\v_{0,2} = F_{y^2\lam}-f_{y^2z}f_{z^2}^{-1} F_{z\lam}$
    & $\v_{0,2} = F_{y^2\lam}$
\\
$5$ & $\v_{0,3} - 3a_{1,3} \v_{1,1}$
    & $\v_{0,4} - \frac{\eta}{10}a_{1,5}\v_{1,1}$
\\
$6$ & $\v_{0,4}-2a_{1,3}\v_{1,2}-(a_{1,4}-2a_{2,2}a_{1,3})\v_{1,1}
+\frac{1}{3}a_{1,3}^2\v_{2,0}$
  & $\v_{0,6} - \frac{\eta}{2}a_{1,5}\v_{1,3} + \frac{1}{160}a_{1,5}^2\v_{2,0}$
\\ &
   & \ \ $+(\frac{1}{16}a_{1,5}a_{2,4}-\frac{1}{14}\eta \, a_{1,7})\v_{1,1}$
\\
$7$ & $\v_{0,5} -\frac{10}{3}a_{1,3} \v_{1,3}+\frac{5}{3} a_{1,3}^2 \v_{2,1}$
    & $F_{y^8\lam} + ...$
\\ & \ \ $-(a_{1,5}-\frac{10}{3}a_{2,3}a_{1,3}+\frac{5}{3}a_{3,1}a_{1,3}^2)\v_{1,1}$
    &
\\ & \ \ $+\frac{5}{6}(a_{1,4}
           -2a_{2,2}a_{1,3})(a_{1,3}\v_{2,0}-3\v_{1,2}-3a_{2,2}\v_{1,1})$
    &
\\
 \hline
\end{tabular}
\eqno (\thetheorem)
 $$
{\bf Comment to Table~(\ref {tab:versalD}).}
In this table, the variables and the parameters are denoted by
$x=k_1$, $y=k_2$, $z=(k_3,\dots,k_\nu)$, and $\lam=(\lam_1,\dots,\lam_l)$,
resp., moreover
$f_{x^ay^bz^c}:=\frac{\partial^{a+b+c}f(0)}{\partial x^a\partial y^b\partial z^c}$
and
$F_{x^ay^bz^c\lam}:=
 \frac{\partial^{a+b+c+1}F(0,0)}{\partial x^a\partial y^b\partial z^c\partial\lam}$.
The real numbers
$a_{i,j}=\frac{\partial f^{i+j}(0)}{\partial x^i\partial y^j}+\dots$
are the same as those in Table~(\ref {tab:defD}). Furthermore, one denotes
$\v_{0,j}=\v_{j+1}$ where $\v_\mu\in\RR^l$ are the
vectors assigned according to Statement~$\ref {st:versal}$ to the deformation
$F(0,k_2,\dots,k_\nu,\lam)$ of $f(0,k_2,\dots,k_\nu)$, see
Table~(\ref {tab:versal}). Other
vectors $\v_{i,j}
=\frac{\partial F^{i+j+1}(0,0)}{\partial x^i\partial y^j\partial\lam}
+\dots \in\RR^l$
in Table~(\ref {tab:versalD}) are similar to the vectors $\v_{0,j}$ and are
defined by the formulae
 $$
\v_{1,1} = F_{xy\lam} - f_{xyz} f_{z^2}^{-1} F_{z\lam},
 $$
 $$
\v_{1,2} = F_{xy^2\lam} - f_{y^2z} f_{z^2}^{-1} F_{xz\lam}
                        -2 f_{xyz} f_{z^2}^{-1} F_{yz\lam}
          - (f_{xy^2z} - f_{y^2z} f_{z^2}^{-1} f_{xz^2}
                       -2 f_{xyz} f_{z^2}^{-1} f_{yz^2}) f_{z^2}^{-1} F_{z\lam},
 $$
 $$
\v_{1,3} = F_{xy^3\lam}
- (f_{y^3z} -3 f_{y^2z} f_{z^2}^{-1} f_{yz^2}) f_{z^2}^{-1} F_{xz\lam}
-3 (f_{xy^2z} - f_{y^2z} f_{z^2}^{-1} f_{xz^2} -2 f_{xyz} f_{z^2}^{-1} f_{yz^2})
   f_{z^2}^{-1} F_{yz\lam}
 $$
 $$
-3 f_{xyz} f_{z^2}^{-1} F_{y^2z\lam}
-3 f_{y^2z} f_{z^2}^{-1} (F_{xyz\lam} - F_{z^2\lam} f_{z^2}^{-1} f_{xyz})
- (f_{xy^3z} -  f_{y^3z} f_{z^2}^{-1} f_{xz^2}
-3 f_{xy^2z} f_{z^2}^{-1} f_{yz^2}
 $$
 $$
-3 (f_{xyz} (f_{z^2}^{-1} f_{y^2z^2} - 2 [f_{z^2}^{-1} f_{yz^2}]^2 )
+ f_{y^2z} f_{z^2}^{-1} ( f_{xyz^2} - 2f_{yz^2} f_{z^2}^{-1} f_{xz^2}
 - f_{z^3} f_{z^2}^{-1} f_{xyz}))) f_{z^2}^{-1} F_{z\lam} .
 $$
Moreover, if $F(\cdot,\lam)$ is an even function for any $\lam$ then
 $$
\v_{2,2} = F_{x^2y^2\lam}
         - 2 f_{xy^2z} f_{z^2}^{-1}F_{xz\lam}
         - 2 f_{x^2yz} f_{z^2}^{-1}F_{yz\lam} .
 $$
One can obtain $\v_{j,i}$ from $\v_{i,j}$ by replacing all partial derivatives
with respect to $x$ by partial derivatives with respect to $y$, and vica-versa.
The above formulae are written for $\nu=3$, $l=1$, but they are easily
transformed to the corresponding formulae for any $\nu\ge2$, $l\ge1$,
see~(\ref {note:tab}).
In the right column of Table~(\ref {tab:versalD}), the function $F(\cdot,\lam)$
is supposed to be even for each $\lam$. This leads to a simplification of the
formulae for the vectors $\v_{i,j}$.

For a proof of Statement~\ref {st:versalD}, let us perform, for any
sufficiently small value of the parameter $\lam$, a change of variables
similar to that in the proof of Statement~\ref {st:defD}. In the new
variables $k_1,k_2,\tilde K_3,\dots,\tilde K_\nu,\lam$, we have
$F=A(\lam)+\Psi(k_1,k_2,\lam)+Q(\tilde K_3,\dots,\tilde K_\nu)$, where
$A$ and $\Psi$ are smooth functions, $A(0)=a$, $\Psi(\cdot,0)=\psi$,
$\tilde K=\tilde K(k,\lam)$.
Now let us perform the change of variables
$(k_1,k_2)\to(\tilde k_1,\tilde k_2)$ given by
$\tilde k_1=k_1+\lam_2k_2^2+\dots+\lam_{\mu-3}k_2^{\mu-3}$,
$\tilde k_2=k_2$, where the $\lam_j$ are determined implicitely (also uniquely
and independently of $\mu$) by the following condition: after this change of
variables the coefficients at the terms
$\tilde k_1\tilde k_2^j$, $3\le j\le\mu-2$, of the Taylor series of
$\psi$ in $\tilde k_1,\tilde k_2$, centred at the origin, vanish,
see the proof of Statement~\ref {st:defD}.
Denote by $\tilde\psi(\tilde k_1,\tilde k_2)$ and
$\tilde\Psi(\tilde k_1,\tilde k_2,\lam)$ the functions which are obtained
from the functions $\psi(k_1,k_2)$ and $\Psi(k_1,k_2,\lam)$, resp., after
this change of variables. Thus, the coefficients at the terms
$\tilde k_1\tilde k_2^j$, $j\le\mu-2$, of the Taylor series of
$\tilde\psi$, centred at zero, vanish.
It follows from the proof of Statement~\ref {st:defD} that, for each
$j\le2\mu-4$, the coefficient at the term $\tilde k_2^j$ equals
$\frac{d_{j+1}}{j!}$,
therefore the coefficients at the terms $\tilde k_2^j$, $j\le\mu-2$,
also vanish, and the coefficient at the term $\tilde k_2^{\mu-1}$ equals
$\frac{d_\mu}{(\mu-1)!}\ne0$.
Using the above properties of the coefficients of the Taylor series of the
function $\tilde\psi$ at zero, one easily proves that the germ of the
function $\tilde k_1\tilde k_2$ at zero (and, thus, the germs of the functions
$\tilde k_1^3$, $\tilde k_1^2+\frac{2d_\mu}{(\mu-2)!}\tilde k_2^{\mu-2}$, and
$\tilde k_2^{\mu-1}$)
belongs to the Jacobian ideal $I_{\nabla\tilde\psi}$.
It follows that the germs of the functions
$\tilde k_1,\tilde k_1^2,\tilde k_2,\tilde k_2^2,\dots,\tilde k_2^{\mu-3}$
are generators of the local algebra $Q_{\nabla\tilde\psi}$
of the gradient map of the function $\tilde\psi$ at the origin (considered as
a vector space), and that the differential operators
$\frac{\partial}{\partial\tilde k_1}|_0$,
$\frac{\partial}{\partial\tilde k_2}|_0$,
$\frac{\partial^2}{\partial\tilde k_2^2}|_0,\dots,
\frac{\partial^{\mu-3}}{\partial\tilde k_2^{\mu-3}}|_0$,
$\frac{\partial^{\mu-2}}{\partial\tilde k_2^{\mu-2}}|_0 - d_\mu
\frac{\partial^2}{\partial\tilde k_1^2}|_0$ are generators of
its dual space $Q_{\nabla\tilde\psi}^*$.
It follows from this, from the conditions on a deformation to be
infinitesimally $R^+-$versal, and from the versality theorem
(see the theorems from~\cite [v.~1,~8.2 and~8.3]{Arnold}), that
the deformation $\Psi$ of the function $\psi$ is $R^+-$versal if and only if
the vectors $\w_2,\dots,\w_{\mu-1},\w_\mu-d_\mu\w_1\in\RR^l$
are linearly independent. Here
$\w_i:=\frac{\partial^{4-i}\tilde\Psi(0,0)}{\partial\tilde k_1^{3-i}\partial\lam}$
for $i=1,2$, and
$\w_i:=\frac{\partial^{i-1}\tilde\Psi(0,0)}{\partial\tilde k_2^{i-2}\partial\lam}$
for $i\ge3$.

In order to compute the vectors $\w_\mu$, let us express the vectors $\v_{q,p}
:=\frac{\partial^{q+p+1}\Psi(0)}{\partial k_1^q\partial k_2^p\partial\lam}
=\frac{\partial^{q+p+1}F(0)}{\partial k_1^q\partial k_2^p\partial\lam}+\dots$
in terms of the partial derivatives of the function $F$ at the origin.
Computations, similar to those for the vectors $\v_\mu
=\frac{\partial^\mu F(0,0)}{\partial k_1^{\mu-1}\partial\lam}+\dots\in\RR^l$,
see the proof of Statement~\ref {st:versal}, lead to the above formulae for
$\v_{q,p}$, see the comment to Table~(\ref{tab:versalD}).
Performing the change of variables $(k_1,k_2)\to(\tilde k_1,\tilde k_2)$
and using the values of the coefficients $\lam_j$ from the proof of
Statement~\ref {st:defD}, one obtains the required expressions for the
vectors $\w_\mu$ in terms of the numbers $a_{q,p}$ and the vectors $\v_{q,p}$.
This proves Statement~\ref {st:versalD}.

Statement~\ref {st:versalXe} is similarly proved on the basis of
Theorem~\ref {verMO}. Here one proves that the germ of the function
$\tilde k_1^3\tilde k_2$ at the origin (and, thus, the germs of the even
functions
$\tilde k_1^6$, $\varepsilon\tilde k_1^4+2\eta\tilde k_1^2\tilde k_2^2$,
$\eta\tilde k_1^2\tilde k_2^2
+\frac{\mu-3}{(2\mu-6)!} x^{\varepsilon,\eta}_{e,\mu} \tilde k_2^{2\mu-6}$,
$\tilde k_1\tilde k_2^3$, and $\tilde k_2^{2\mu-4}$)
belongs to the even Jacobian ideal $I^e_{\nabla\tilde\psi}$.
It follows that the germs of the even functions
$\tilde k_1^2,\tilde k_1\tilde k_2,\tilde k_2^2,
\tilde k_1^4,\tilde k_2^4,\dots,\tilde k_2^{2\mu-8}$ are generators of the
even local algebra $Q^e_{\nabla\tilde\psi}$
of the gradient map of the function $\tilde\psi$ at the origin (considered as
a vector space), and that the differential operators
$\frac{\partial^2}{\partial\tilde k_1^2}|_0$,
$\frac{\partial^2}{\partial\tilde k_1\partial\tilde k_2}|_0$,
$\frac{\partial^2}{\partial\tilde k_2^2}|_0$,
$\frac{\partial^4}{\partial\tilde k_2^4}|_0,\dots,
\frac{\partial^{2\mu-8}}{\partial\tilde k_2^{2\mu-8}}|_0$,
$\frac{\partial^{2\mu-6}}{\partial\tilde k_2^{2\mu-6}}|_0
+ \frac{\mu-3}{48} x^{\varepsilon,\eta}_{e,\mu}
(\varepsilon \frac{\partial^4}{\partial\tilde k_1^4}|_0
-12 \eta \frac{\partial^4}{\partial\tilde k_1^2\partial\tilde k_2^2}|_0)$
are generators of its dual space $(Q^e_{\nabla\tilde\psi})^*$.
It follows from this and Theorem~\ref {verMO} that
the deformation $\Psi$ of the function $\psi$ is $R_O-$versal if and only if
the vectors $\x_{e,2}^{\varepsilon,\eta} ,\dots,\x_{e,\mu-1}^{\varepsilon,\eta}$,
$\x_{e,\mu}^{\varepsilon,\eta} + \frac{\mu-3}{48} x_{e,\mu}^{\varepsilon,\eta}
\x^{\varepsilon,\eta}_{e,1}\in\RR^l$ are linearly independent.

Statements~\ref {st:versal},~\ref {st:versale},~\ref {st:versalD},
and~\ref {st:versalXe} demonstrate that versal deformations are typical if
the number $l$ of parameters is large enough. More precisely: if $l\ge\mu-1$
then the conditions from~\ref {st:versal},~\ref {st:versalD} on a deformation
of a germ $f$ to be $R^+-$versal
(or the conditions from~\ref {st:versale},~\ref {st:versalXe} on an even
deformation of an even germ $f$ to be $R_O-$versal) are fulfilled for
``typical'' (resp.\ ``typical even'') $l-$parameter deformations of $f$.

Here, by typical $l-$parameter (even) deformations of $f$ we mean (even)
deformations of $f$ which form an open dense subspace in the space of all
(even) $l-$parameter deformations of $f$ in the convergence topology with a
finite number of derivatives on each compact set.


\begin{thebibliography}{02}

\bibitem{Arnold}  V.I.~Arnold, S.M.~Gusein-Zade, A.N.~Varchenko.
{\it Singularities of Differentiable Maps}.
Monographs in Mathematics, Birkhauser, Boston, 1985 and 1988.

\bibitem{Arnold2} V.I.~Arnold,
{\it Singularities of Caustics and Wave Fronts}.
Kluwer Academic Publishers, 1990.

\bibitem{Itogi} V.I.~Arnol'd, V.A.~Vasil'ev, V.V.~Goryunov, O.V.~Lyashko.
{\it Singularities. Local and Global Theory}.
Encyclopaedia of Math.\ Sci. {\bf 6},
Springer-Varlag, Berlin-Heidelberg-New York, 1993.

\bibitem{Beer} M.~Beer, {\it Endliche Bestimmtheit und universelle
Entfaltungen von Keimen mit Gruppenoperation}. Diplomarbeit, Regensburg 1976.

\bibitem{LakshtanovMinlos} E.L.~Lakshtanov, R.A.~Minlos,
{\it Spectrum of transfer-matrix two-particle bound states},
Func.\ Anal.\ i ego Pril.\ {\bf 38(3)} (2003), 52-69.

\bibitem{MMinfP} V.A.~Malyshev, R.A.~Minlos,
{\it Linear Infinite-Particle Operators}.
Amer.\ Math.\ Soc.\ Publ., 1995.

\bibitem{MamatovMinlos} SH.S.~Mamatov, R.A.~Minlos,
{\it Bound states of two-particle cluster operator},
Theor.\ and Math.\ Phys. {\bf 79(2)} (1989), 163-182.


\bibitem{Sinay} R.A.~Minlos, Ya.G.~Sinay,
{\it Study of the spectrum of stochastic operators appearing in the lattice
models of gas},
Theor.\ and Math.\ Phys. {\bf 2} (1970), 230-243.

\bibitem{Poenaru} V.~Po\`enaru,
{\it Singularit\'es $C^\infty$ en pr\'esence de sym\'etrie},
Lecture Notes in Math.~{\bf 510}, Springer-Verlag,
Berlin-Heidelberg-New York, 1976.

\bibitem{Peter} P.~Slodowy,
{\it Einige Bemerkungen zur Entfaltung symmetrischer Funktionen},
Math.~Z.\ {\bf 158} (1978), 157-170.

\bibitem{vanderWaerden} B.L.~van der Waerden. {\it Algebra}.
Springer-Verlag, Berlin-Heidelberg-New York, 1967.

\bibitem{Wass} G.~Wassermann,
{\it Classification of singularities with compact abelian symmetry},
Singularities, Banach Center Publications {\bf 20} (1988), 475-498.

\end{thebibliography}
\end{document}